\documentclass[11pt,twoside,reqno]{amsarticlevadim} 

\usepackage[utf8]{inputenc}
\usepackage{graphicx}
\usepackage[all]{xy}
\usepackage{amssymb}
\usepackage{amsmath} 
\usepackage{amsthm}
\usepackage{verbatim}
\usepackage{latexsym}
\usepackage{mathrsfs}

\usepackage{mathrsfs}

\newcommand{\rusi}{\tilde o}
\newcommand{\rusch}{\tilde n}
\newcommand{\ruszh}{\tilde f}
\newenvironment{myItemize}{ 
  \begin{list}{\raisebox{2.2pt}{$\centerdot$}}{%
      \setlength\leftmargin{18pt}
      \setlength\labelwidth{20pt}
    }
  }{
  \end{list}
}

\newenvironment{myDescription}
{\list{}{\labelwidth10pt\leftmargin15pt
    }}
{\endlist}

\newenvironment{myEnumerate}{ 
  \begin{list}{\labelenumi}{%
      \usecounter{enumi}
      \setlength\leftmargin{25pt}
      \setlength\labelwidth{20pt}
    }
  }{
  \end{list}
}

\newenvironment{myAlphanumerate}{ 
  \begin{list}{\labelenumii}{%
      \usecounter{enumii}
      \setlength\leftmargin{25pt}
      \setlength\labelwidth{20pt}
    }
  }{
  \end{list}
}

\makeindex

\renewcommand{\ll}{\mathscr{L}}

\renewcommand{\Cup}{\bigcup}
\renewcommand{\Cap}{\bigcap}

\renewcommand{\a}{\alpha}
\renewcommand{\b}{\beta}
\newcommand{\g}{\gamma}
\newcommand{\e}{\varepsilon}
\renewcommand{\d}{\delta}
\renewcommand{\k}{\kappa}
\renewcommand{\l}{\lambda}
\newcommand{\f}{\varphi}
\renewcommand{\o}{\omega}
\newcommand{\s}{\sigma}
\renewcommand{\th}{^{\text{th}}}

\newcommand{\esm}{\preccurlyeq} 

\newcommand{\es}{\varnothing}

\newcommand{\A}{\mathcal{A}}
\newcommand{\B}{\mathcal{B}}
\newcommand{\C}{\mathcal{C}}
\newcommand{\E}{\mathcal{E}}

\newcommand{\Q}{\mathbb{Q}} 
\newcommand{\R}{\mathbb{R}} 

\newcommand{\cat}{{}^\frown}  
\newcommand{\wins}{\uparrow}

\renewcommand{\P}{\mathbb{P}}             
\newcommand{\forces}{\Vdash}
\newcommand{\la}{\langle} 
\newcommand{\ra}{\rangle}
\newcommand{\Po}{\mathcal{P}}             
\renewcommand{\S}{\mathcal{S}} 
\newcommand{\NS}{{\operatorname{NS}}}
\newcommand{\cf}{\operatorname{cf}}

\newcommand{\Succ}{\operatorname{Succ}} 
\newcommand{\Sk}{{\operatorname{Sk}}} 
\newcommand{\card}{\operatorname{card}} 
\newcommand{\sprt}{\operatorname{sprt}} 
\newcommand{\reg}{\operatorname{reg}} 
\newcommand{\id}{\operatorname{id}}
\newcommand{\CUB}{\operatorname{CUB}}
\newcommand{\Th}{\operatorname{Th}}
\newcommand{\cl}{{\operatorname{cl}}}
\newcommand{\height}{{\operatorname{ht}}}
\newcommand{\pr}{\operatorname{pr}}
\newcommand{\rest}{\!\restriction\!}
\newcommand{\restl}{\restriction}  
\newcommand{\dom}{\operatorname{dom}}
\newcommand{\ran}{\operatorname{ran}}
\newcommand{\OTP}{\operatorname{OTP}}
\newcommand{\tp}{\operatorname{tp}}  
\newcommand{\stp}{\operatorname{stp}}  
\renewcommand{\lg}{\operatorname{length}}  
\renewcommand{\le}{\leqslant}  
\renewcommand{\ge}{\geqslant}  
\newcommand{\less}{\lessdot}  
\newcommand{\sd}{\,\triangle\,}             

\newcommand{\Sii}{{\Sigma_1^1}}
\newcommand{\Pii}{{\Pi_1^1}}
\newcommand{\Dii}{{\Delta_1^1}}
\newcommand{\I}{\mathcal{I}}
\newcommand{\J}{\mathcal{J}}
\newcommand{\EF}{\operatorname{EF}} 

\newcommand{\PlOne}{\,{\textrm{\bf I}}}
\newcommand{\PlTwo}{\textrm{\bf I\hspace{-1pt}I}}

\newcommand{\Land}{\bigwedge}
\newcommand{\Lor}{\bigvee}
\DeclareMathOperator*{\Exists}{\raisebox{-6pt}{\huge$\exists$}}
\DeclareMathOperator*{\Forall}{\raisebox{-6pt}{\huge$\forall$}}

\swapnumbers
\newtheorem*{Thm*}{Theorem}
\newtheorem{Thm}{Theorem}
\newtheorem{Lemma}[Thm]{Lemma}
\newtheorem{Cor}[Thm]{Corollary}
\newtheorem*{Fact}{Fact}
\newtheorem*{Convention}{Convention}

\newenvironment{claim}[1]{\text{ }\vspace{7pt}\newline\noindent\textbf{Claim #1.}}{\hspace*{\fill}}
\newenvironment{claim*}{\vspace{7pt}\noindent\textbf{Claim.}}{}

\theoremstyle{definition}
\newtheorem{Def}[Thm]{Definition}

\newtheorem{Example}[Thm]{Example}

\theoremstyle{remark}
\newtheorem*{Remark}{Remark}
\newtheorem{RemarkN}[Thm]{Remark}
\newtheorem*{Open}{Open Problem}

\newcommand{\proofvpara}{\text{}}

\newenvironment{proofVOf}[1] {\hfill\vspace{5pt}\\\noindent \textbf{Proof of #1.}\ignorespaces\renewcommand{\proofvpara}{\text{#1}}}
{\nopagebreak\hspace*{\fill}\mbox{$\square_{\,\proofvpara}$}\\\vspace{-8pt}}

\newenvironment{proofVOf2}[2] {\hfill\vspace{5pt}\\\noindent \textbf{Proof of #1.}\ignorespaces\renewcommand{\proofvpara}{\text{#2}}}
{\nopagebreak\hspace*{\fill}\mbox{$\square_{\,\proofvpara}$}\\\vspace{-8pt}}

\newenvironment{proofV}[1] {\hfill\\\noindent \textbf{Proof.}\ignorespaces \renewcommand{\proofvpara}{\text{#1}}}
{\nopagebreak\hspace*{\fill}\mbox{$\square_{\,\proofvpara}$}\\}

\author{Sy-David Friedman\\ Kurt Gödel Research Center\\University of Vienna\\\text{ }\\ 
  Tapani Hyttinen and Vadim Kulikov\textsuperscript{\tiny 1}\\ Department of Mathematics and Statistics\\ University of Helsinki\\\text{ }\\ }

\title{\mbox{$\!\!\!\!\!\!\!$Generalized~Descriptive~Set~Theory} and Classification~Theory\newline}
\date{November 8, 2011 \\(file updated \today)}

\begin{document}

\setcounter{page}{-1}
\maketitle
\footnotetext[1]{EDIT 2025: The name of the author has changed to Vadim Weinstein}

\newpage\thispagestyle{empty}

\begin{abstract}
  Descriptive set theory is mainly concerned with studying subsets of the space of all countable
  binary sequences. In this paper we study the generalization where countable is replaced by uncountable. 
  We explore properties of generalized Baire and Cantor spaces, equivalence relations and their Borel reducibility. 
  The study shows that the descriptive
  set theory looks very different in this generalized setting compared to the classical, countable case.
  We also draw the connection between the stability theoretic complexity of first-order theories 
  and the descriptive set theoretic complexity of their isomorphism relations. Our results suggest 
  that Borel reducibility on uncountable structures 
  is a model theoretically natural way to compare the complexity of isomorphism relations.
\end{abstract}

\newpage
{\noindent{\Large{\bf Acknowledgement}}}
\newline
\vspace{10pt}

The authors wish to thank
the John Templeton Foundation for its generous support
through its project Myriad Aspects of Infinity (ID \#13152).
The authors wish to thank also 
Mittag-Leffler Institute (the Royal Swedish Academy of Sciences).

The second and the third authors wish to thank the Academy of Finland for 
its support through its grant number 1123110. 

The third author wants to express his gratitude to the
Research Foundation of the University of
Helsinki and the Finnish National Graduate School in Mathematics and its
Applications for the financial support during the work.

We are grateful to Jouko Väänänen for the useful discussions and 
comments he provided on a draft of this paper.

\newpage

\thispagestyle{empty}
\tableofcontents
\thispagestyle{empty}
\chapter{History and Motivation}\label{sec:HistMotiv}
There is a long tradition in studying connections between Borel structure of Polish spaces (descriptive set theory) 
and model theory. The connection arises from the fact that any class of countable structures can be coded into a subset of the
space $2^\o$ provided all structures in the class have domain~$\o$. A survey on this topic is given in \cite{Hjo}.
Suppose $X$ and $Y$ are subsets of $2^\o$ and let $E_1$ and $E_2$ be equivalence relations on $X$ and $Y$ respectively.
If $f\colon X\to Y$ is a map such that
$E_1(x,y)\iff E_2(f(x),f(y))$, we say that $f$ is a \emph{reduction of $E_1$ to $E_2$}.
If there exists a Borel or continuous reduction, we say that $E_1$ is Borel or continuously \emph{reducible} to $E_2$, 
denoted $E_1\le_B E_2$ or $E_1\le_c E_2$. 
The mathematical meaning of this is
that $f$ \emph{classifies $E_1$-equivalence in terms of $E_2$-equivalence}. 

The benefit of various reducibility and irreducibility theorems is roughly the following.
A reducibility result, say $E_1\le_B E_2$, tells us that $E_1$ is at most as complicated as $E_2$; once you understand $E_2$, you
understand $E_1$ (modulo the reduction).
An irreducibility result, $E_1\not\le_B E_2$ tells that there is no hope in trying to classify $E_1$ in terms of
$E_2$, at least in a ``Borel way''. From the model theoretic point of view, the isomorphism relation, and the elementary equivalence relation 
(in some language) on some class of structures are the equivalence relations of main interest. But model theory in general
does not restrict itself to countable structures. Most of stability theory and Shelah's classification theory 
characterizes first-order theories in terms of their uncountable models. This leads to the generalization 
adopted in this paper.
We consider the space $2^\k$ for an uncountable cardinal $\k$ with the idea 
that models of size $\k$ are coded into elements of
that space.

This approach, to connect such uncountable descriptive set theory with model theory, began in the early 1990's. 
One of the pioneering papers was by Mekler and V\"a\"an\"anen \cite{MekVaa}. A survey on the research done in
1990's can be found in \cite{Vaa2} and a discussion of the motivational background for this work in \cite{Vaa1}.
A more recent account is given the book \cite{Vaa3}, Chapter~9.6.

Let us explain how our approach differs from the earlier ones and why it is useful.
For a first-order complete countable theory in a countable vocabulary $T$ and a cardinal $\k\ge \o$, define
$$S_T^\k=\{\eta\in 2^\k\mid \A_\eta\models T\}\text{ and }\cong_T^\k\,=\{(\eta,\xi)\in (S^\k_T)^2\mid \A_\eta\cong\A_\xi\}.$$
where $\eta\mapsto \A_\eta$ is some fixed coding of (all) structures of size $\k$.
We can now define the partial order on the set of all theories as above by 
$$T\le^\k T'\iff \,\cong_T^\k\,\le_B\, \cong_{T'}^\k.$$
As pointed out above, $T\le^\k T'$ says that $\cong_{T}^\k$ is at most as difficult to classify as $\cong_{T'}^\k$.
But does this tell us whether $T$ is a simpler theory than $T'$? Rough answer: \emph{If $\k=\o$, then no but if $\k>\o$, then yes}.

To illustrate this, let $T=\Th(\Q,\le)$ be the theory of the order of the rational numbers (DLO) and let
$T'$ be the theory of a vector space over the field of rational numbers. Without loss of generality we may assume that
they are models of the same vocabulary.

It is easy to argue that the model class defined by $T'$ is strictly simpler than that of $T$.
(For instance there are many questions about $T$, unlike $T'$, that cannot be answered in ZFC; say existence of a saturated model.)
On the other hand $\cong^\o_T\,\le_B \,\cong_{T'}^\o$ and $\cong^\o_{T'}\,\not\le_B\,\cong^\o_T$ 
because there is only one countable model of $T$ and there are infinitely many countable models of $T'$. 
But for $\k>\o$ we have
$\cong_{T}^\k\,\not\le_B\, \cong_{T'}^\k$ and $\cong_{T'}^\k\,\le_B\, \cong_{T}^\k$, since there are $2^\k$ equivalence classes
of $\cong_T^\k$ and only one equivalence class of $\cong_{T}^\k$.
Another example, introduced in Martin Koerwien's Ph.D. thesis and his article \cite{Koe} shows that there exists an $\o$-stable 
theory without DOP and without OTOP with depth 2 for which $\cong^\o_T$ is not Borel, while we show here that
for $\k>2^\o$, $\cong^\k_T$ is Borel for all classifiable shallow theories.

The results suggest that the order $\le^\k$ 
for $\k>\o$ corresponds naturally to the classification of theories in stability theory: 
the more complex a theory is from the viewpoint of stability theory, the higher it seems to sit in the
ordering $\le^\k$ and vice versa. 
Since dealing with uncountable cardinals often implies the need for various 
cardinality or set theoretic assumptions beyond ZFC, the results are not always as simple
as in the case $\k=\o$, but they tell us a lot.
For example, our results easily imply the following (modulo some mild cardinality assumptions on $\k$):
\begin{myItemize}
\item If $T$ is deep and $T'$ is shallow, then $\cong_T\,\not\le_B\,\cong_{T'}$.
\item If $T$ is unstable and $T'$ is classifiable, then $\cong_T\,\not\le_B\,\cong_{T'}$.
\end{myItemize}

\chapter{Introduction}
\section{Notations and Conventions}\label{sec:NotationsandC}
\subsection{Set Theory}
We use standard set theoretical notation:
\begin{myItemize}
\item $A\subset B$ means that
  $A$ is a subset of $B$ or is equal to $B$.
\item $A\subsetneq B$ means proper subset.
\item Union, intersection and set theoretical difference are denoted respectively by $A\cup B$,
$A\cap B$ and $A\setminus B$. For larger unions and intersections $\bigcup_{i\in I}A_i$ etc..
\item $\Po(A)$ is the power set of $A$ and $[A]^{<\k}$ is the set of subsets of $A$ of size $<\k$ 
\end{myItemize}
\label{page:FuncGraph}

Usually the Greek letters $\kappa$, $\lambda$ and $\mu$ will stand for cardinals and
$\alpha$, $\beta$ and $\gamma$ for ordinals, but this is not strict. Also $\eta,\xi,\nu$
are usually elements of $\k^\k$ or $2^\k$ and $p,q,r$ are elements of $\k^{<\k}$ or $2^{<\k}$.
$\cf(\alpha)$ is the cofinality
of $\alpha$ (the least ordinal $\beta$ for which there exists
an increasing unbounded function $f\colon \beta\to\alpha$).

By $S^\k_\l$ we mean $\{\alpha<\k\mid \cf(\alpha)=\l\}$. 
A  $\l$-\emph{cub set} is a subset of a limit ordinal (usually of cofinality $>\l$) which is unbounded and
contains suprema of all bounded increasing sequences of length $\l$. A set is \emph{cub} if it is $\l$-cub for all $\l$.
A set is \emph{stationary} if it intersects all cub sets and $\l$\emph{-stationary} if it intersects
all $\l$-cub sets. Note that $C\subset \k$ is $\l$-cub if and only if $C\cap S^\k_\l$ is $\l$-cub and
$S\subset \k$  is $\l$-stationary if and only if $S\cap S^\k_\l$ is (just) stationary. 

If $(\P,\le)$ is a forcing notion, we write $p\le q$ if $p$ and $q$ are in $\P$ and $q$ forces more than $p$.
Usually $\P$ is a set of functions equipped with inclusion and $p\le q\iff p\subset q$. 
In that case $\es$ is the weakest condition and we write $\P\forces \f$ to mean $\es\forces_{\P}\f$.

\subsection{Functions}\label{ssec:Functions}

We denote by $f(x)$ the value of $x$ under the mapping $f$ and by $f[A]$ or just $fA$ the image
of the set $A$ under $f$. Similarly $f^{-1}[A]$ or just $f^{-1}A$ indicates the inverse image of $A$.
Domain and range are denoted respectively by $\dom f$ and $\ran f$.

If it is clear from the context that $f$ has an inverse, then $f^{-1}$ denotes that inverse.
For a map $f\colon X\to Y$ \emph{injective} means the same as \emph{one-to-one} and \emph{surjective} the same as 
\emph{onto}

Suppose $f\colon X\to Y^\a$ is a function with range consisting of sequences of elements of $Y$ of length $\a$.
The projection $\pr_\b$ is a function $Y^\a\to Y$ defined by $\pr_\b((y_i)_{i<\a})=y_\b$.
For the coordinate functions of $f$ we use the notation
$f_\b=\pr_\b\circ f$ for all~$\b<\a$.

By support of a function $f$ we mean the subset of $\dom f$ in which $f$ takes non-zero values, whatever ``zero'' means
depending on the context (hopefully never unclear). The support of $f$ is denoted by~$\sprt f$.

\subsection{Model Theory}

In Section \ref{sec:CodingModelss} we fix a countable vocabulary and assume that all theories are theories in this vocabulary.
Moreover we assume that they are first-order, complete and countable. By $\tp(\bar a/A)$ we denote the complete type
of $\bar a=(a_1,\dots,a_{\lg \bar a})$ over $A$ where $\lg \bar a$ is the length of the sequence $\bar a$. 

We think of models as tuples $\A=\la \dom \A,P^\A_n\ra_{n<\o}$ where the $P_n$ are relation symbols in the vocabulary
and the $P_n^\A$ are their interpretations. If a relation $R$ has arity $n$ (a property of the vocabulary), then for its
interpretation it holds that
$R^\A\subset (\dom \A)^n$. In Section \ref{sec:CodingModelss} we adopt more conventions concerning this.

In Section \ref{ssec:SDIR} and Chapter \ref{chapter:ComplexityofIsomRel} we will use the following 
stability theoretical notions
stable, superstable, DOP, OTOP, shallow and $\k(T)$. Classifiable means superstable with no DOP nor OTOP,
the least cardinal in which $T$ is stable is denoted by~$\l(T).$

\subsection{Reductions}\label{def:Reductions}
Let $E_1\subset X^2$ and $E_2\subset Y^2$ be equivalence relations on $X$ and $Y$ respectively. A function
$f\colon X\to Y$ is \emph{a reduction} of $E_1$ to $E_2$ if for all $x,y\in X$ we have that
$xE_1y\iff f(x)E_2f(y)$. Suppose in addition that $X$ and $Y$ are topological spaces. Then
we say that $E_1$ is \emph{continuously reducible to} $E_2$, if there exists a continuous
reduction from $E_1$ to $E_2$ and we say that $E_1$ is \emph{Borel reducible to} $E_2$ if there is
a Borel reduction. For the definition of Borel adopted in this paper, see Definition~\ref{def:Borel}.
We denote the fact that $E_1$ is continuously reducible to $E_2$ by $E_1\le_c E_2$ and
respectively Borel reducibility by $E_1\le_B E_2$.

We say that relations $E_2$ and $E_1$ are (Borel) \emph{bireducible} to each other if $E_2\le_B E_1$ and $E_1\le_B E_2$.

\section{Ground Work}
\subsection{Trees and Topologies}

Throughout the paper $\k$ is assumed to be an uncountable regular cardinal which satisfies
$$\k^{<\k}=\k\eqno(*)$$\label{kappapoten}
(For justification of this, see below.)
We look at the space $\k^\k$, i.e.
the functions from $\k$ to $\k$ and the space formed by the initial segments $\k^{<\k}$.
It is useful to think of $\k^{<\k}$ as a tree ordered by inclusion and of $\k^\k$ as a topological space
of the branches of $\k^{<\k}$; the topology is defined below. Occasionally we work in $2^\k$ and $2^{<\k}$
instead of $\k^\k$ and $\k^{<\k}$. 

\label{page:FuncGraph}
\begin{Def}
  A \emph{tree} $t$ is a partial order with a root in which the sets $\{x\in t\mid x<y\}$ are well
  ordered for each $y\in t$. A \emph{branch} in a tree is a maximal linear suborder.
  
  A tree is called \emph{a $\k\l$-tree}, if there are no branches of length $\lambda$ or higher and no element has
  $\ge \k$ immediate successors. If $t$ and $t'$ are trees, we write $t\le t'$ to mean that there exists an
  order preserving map $f\colon t\to t'$, $a<_t b \Rightarrow f(a)<_{t'}f(b)$.
\end{Def}

\begin{Convention}
  Unless otherwise said, by a tree $t\subset (\k^{<\k})^n$ we mean a tree with domain being a downward closed subset of
  $$(\k^{<\k})^n\cap \{(p_0,\dots, p_{n-1})\mid \dom p_0=\cdots=\dom p_{n-1}\}$$ 
  ordered as follows: $(p_0,\dots,p_{n-1})<(q_0,\dots,q_{n-1})$ if $p_i\subset q_i$ for all
  $i\in \{0,\dots,{n-1}\}$. It is always a $\k^+,\k+1$-tree. 
\end{Convention}

\begin{Example}\label{example:Trees}
  Let $\a<\k^+$ be an ordinal and let $t_\a$
  be the tree of descending sequences in $\a$ ordered by end extension. The root is the empty sequence.
  It is a $\k^+\o$-tree. Such $t_\a$ can be embedded into $\k^{<\o}$, but note that not all subtrees of $\k^{<\o}$ are
  $\k^+\o$-trees (there are also $\k^+,\o+1$-trees).
\end{Example}

In fact the trees $\k^{<\b}$, $\b\le\k$ and $t_\a$ are universal in the following sense:

\begin{Fact}[$\k^{<\k}=\k$]\label{thm:Embedkappapluskappatree}
   Assume that $t$ is a $\k^+,\b+1$-tree, $\b\le\k$ and $t'$ is $\k^+\o$-tree. Then 
   \begin{myEnumerate}
   \item there is an embedding $f\colon t\to \k^{<\b}$,
   \item and a strictly order preserving map $f\colon t'\to t_\a$ for some $\a<\k^+$ (in fact there is also such an embedding~$f$).\qed
   \end{myEnumerate}
\end{Fact}

Define the topology on $\k^\k$ as follows. For each $p\in \k^{<\k}$ define the basic open set
$$N_p=\{\eta\in \k^\k\mid \eta\rest\dom(p)=p\}.$$
Open sets are precisely the empty set and the sets of the form $\Cup X$, where $X$ is a collection of basic open sets.
Similarly for $2^\k$.

There are many justifications for the assumption $(*)$ which will be most apparent after seeing the proofs
of our theorems. The crucial points can be summarized as follows: if $(*)$ does not hold, then
\begin{myItemize}
\item the space $\k^\k$ does not have a dense subset of size $\k$,
\item there are open subsets of $\k^\k$ 
  that are not $\k$-unions of basic open sets which makes controlling Borel sets difficult
  (see Definition \ref{def:Borel} on page~\pageref{def:Borel}).
\item Vaught's generalization of the Lopez-Escobar theorem (Theorem \ref{thm:BorelIsLkk}) fails, see Remark \ref{separationremark} 
  on page~\pageref{separationremark}.
\item The model theoretic machinery we are using often needs this cardinality assumption (see e.g. Theorem \ref{thm:SepMkk} and 
  proof of Theorem~\ref{thm:NotDiiOrNotBorelList}).
\end{myItemize}
Initially the motivation to assume $(*)$ was simplicity. Many statements concerning the space $\k^{<\k}$ are independent
of ZFC and using $(*)$ we wanted to make the scope of such statements neater. 
In the statements of (important) theorems we mention the assumption explicitly.

Because the intersection of less than $\k$ basic open sets is either empty or a basic open set, we get the following.

\begin{Fact}[$\k^{<\k}=\k$] The following hold for a topological space $P\in \{2^\k,\k^\k\}$:
  \begin{myEnumerate}
    \item The intersection of less than $\k$  basic open sets is either empty or a basic open set,
    \item The intersection of less than $\k$ open sets is open,
    \item Basic open sets are closed,
    \item $|\{A\subset P\mid A\text{ is basic open}\}|=\k$,
    \item $|\{A\subset P\mid A\text{ is open}\}|=2^\k$.
  \end{myEnumerate}
\end{Fact}

In the space $\k^\k\times\k^\k=(\k^\k)^2$ we define the ordinary product topology. 

\begin{Def}\label{def:SiiPiiDii}
  A set $Z\subset \k^\k$ is $\Sigma_1^1$ if
  it is a projection of a closed set $C\subset (\k^\k)^2$. A set is $\Pi_1^1$ if it is the complement
  of a $\Sigma_1^1$ set. A set is $\Delta_1^1$ if it is both $\Sigma_1^1$ and $\Pi_1^1$.   
\end{Def}

As in standard descriptive set theory ($\k=\o$), we have the following:

\begin{Thm}\label{thm:Phomeo}
  For $n<\o$ the spaces $(\k^\k)^n$ and $\k^\k$ are homeomorphic. \qed
\end{Thm}
\begin{Remark}
  This standard theorem can be found for example in Jech's book \cite{Jech}.
  Applying this theorem we can extend the concepts of Definition \ref{def:SiiPiiDii} to subsets
  of $(\k^\k)^n$. For instance a subset $A$ of $(\k^\k)^n$ is $\Sii$ if for a homeomorphism 
  $h\colon (\k^\k)^n\to\k^\k$, $h[A]$ is $\Sii$ according to Definition \ref{def:SiiPiiDii}. 
\end{Remark}

\subsection{Ehrenfeucht-Fra\"iss\'e Games}

We will need Ehrenfeucht-Fra\"iss\'e games in various connections.
It serves also as a 
way of coding isomorphisms.

\begin{Def}[Ehrenfeucht-Fra\"iss\'e games]\label{def:EF1}
  Let $t$ be a tree, $\k$ a cardinal and $\A$ and $\B$ structures with domains $A$ and $B$ respectively.
  Note that $t$ might be an ordinal.
  The game $\EF^\k_t(\A,\B)$  
  is played by players $\PlOne$ and $\PlTwo$ as follows. 
  Player $\PlOne$ chooses subsets of $A\cup B$ and climbs up 
  the tree $t$ and player $\PlTwo$ chooses partial functions $A\to B$ as follows. Suppose a sequence 
  $$(X_i,p_i,f_i)_{i<\g}$$
  has been played 
  (if $\g=0$, then the sequence is empty).
  Player $\PlOne$ picks a set $X_\g\subset A\cup B$ of cardinality strictly less than $\k$
  such that $X_\delta\subset X_\g$ for all ordinals $\delta<\gamma$. Then player $\PlOne$ picks 
  a $p_\g\in t$ which is $<_t$-above all $p_\d$ where $\d<\g$. 
  Then player $\PlTwo$ chooses a partial function $f_\g\colon A\to B$ such that
  $X_\g\cap A\subset \dom f_\g$, $X_\g\cap B\subset \ran f_\g$, $|\dom f_\g|<\k$ and
  $f_\delta\subset f_\g$ for all ordinals $\delta<\gamma$. The game ends when player $\PlOne$ cannot go
  up the tree anymore, i.e. $(p_i)_{i<\g}$ is a branch. Player $\PlTwo$
  wins if 
  $$f=\Cup_{i<\g}f_i$$
  is a partial isomorphism. Otherwise player $\PlOne$ wins. 

  A \emph{strategy} of player $\PlTwo$ in $\EF^\k_t(\A,\B)$ is a function 
  $$\sigma\colon ([A\cup B]^{<\k}\times t)^{<\height(t)}\to \Cup_{I\in [A]^{<\k}}B^{I},$$
  where $[R]^{<\k}$ is the set of subsets of $R$ of size $<\k$ and $\height(t)$ is the \emph{height} of
  the tree, i.e. 
  $$\height(t)=\sup\{\a\mid \a\text{ is an ordinal and there is an order preserving embedding }\a\to t\}.$$
  A strategy of $\PlOne$ is similarly a function  
  $$\tau\colon \Big(\Cup_{I\in [A]^{<\k}}B^{I}\Big)^{<\height(t)} \to [A\cup B]^{<\k}\times t.$$
  We say that a strategy $\tau$ of player $\PlOne$ \emph{beats} strategy $\sigma$ of player $\PlTwo$
  if the play $\tau * \sigma$ is a win for $\PlOne$. The play $\tau * \s$ is just the play where $\PlOne$ uses
  $\tau$ and $\PlTwo$ uses $\s$. Similarly $\sigma$ beats $\tau$ if $\tau * \sigma$ is a win for $\PlTwo$.
  We say that a strategy is a \emph{winning strategy} if it beats all opponents strategies.

  The notation $X\uparrow \EF^\k_t(\A,\B)$ means that player $X$ has a winning strategy in $\EF^\k_t(\A,\B)$
\end{Def}

\begin{Remark}
  By our convention $\dom\A=\dom\B=\k$, so while player $\PlOne$ picks a subset of $\dom\A\cup\dom \B$ he actually
  just picks a subset of $\k$, but as a small analysis shows, this does not alter the game.
\end{Remark}

Consider the game $\EF^\k_t(\A,\B)$, where $|\A|=|\B|=\k$, $|t|\le\k$ and $\height(t)\le \k$. The set of strategies can
be identified with $\k^\k$, for example as follows. 
The moves of player $\PlOne$ are 
members of $[A\cup B]^{<\k}\times t$ and the moves of 
player $\PlTwo$ are members of $\Cup_{I\in [A]^{<\k}}B^I$.
By our convention $\dom\A=\dom\B=A=B=\k$, so these become 
$V=[\k]^{<\k}\times t$ and $U=\Cup_{I\in [\k]^{<\k}}\k^I$. By our cardinality assumption $\k^{<\k}=\k$, 
these sets are of cardinality~$\k$.

Let 
\begin{eqnarray*}
f\colon U\to \k\\
g\colon U^{<\k}\to \k\\
h\colon V\to \k\\
k\colon V^{<\k}\to \k
\end{eqnarray*}
be bijections. Let us assume that $\tau\colon U^{<\k}\to V$ is a strategy of player $\PlOne$ (there cannot be more than $\k$ moves
in the game because we assumed $\height(t)\le\k$). Let 
$\nu_\tau \colon \k\to\k$ be defined by
$$\nu_\tau=h\circ \tau\circ g^{-1}$$ 
and if $\s\colon V^{<\k}\to U$ is a strategy of player $\PlTwo$, let $\nu_\s$ be defined by
$$\nu_\s=f\circ\s\circ k^{-1}.$$
We say that $\nu_\tau$ \emph{codes} $\tau$.

\begin{Thm}[$\k^{<\k}=\k$]\label{thm:EFClosed}
  Let $\l\le\k$ be a cardinal. The set 
  $$C=\{(\nu,\eta,\xi)\in (\k^\k)^3 \mid \nu 
    \text{ codes a w.s. of }\PlTwo\text{ in }\EF^\k_\l(\A_\eta,\A_\xi)\}\subset (\k^\k)^3$$
  is closed. If $\l<\k$, then also the corresponding set for player $\PlOne$
  $$D=\{(\nu,\eta,\xi)\in (\k^\k)^3 \mid \nu \text{ codes a 
    w.s. of }\PlOne\text{ in }\EF^\k_\l(\A_\eta,\A_\xi)\}\subset (\k^\k)^3$$
  is closed.
\end{Thm}
\begin{Remark}
  Compare to Theorem \ref{thm:CodedClosed}.
\end{Remark}
\begin{proof}
  Assuming $(\nu_{_0},\eta_{_0},\xi_{_0})\notin C$, we will show that there is an open neighbourhood $U$ of 
  $(\nu_{_0},\eta_{_0},\xi_{_0})$
  such that $U\subset (\k^\k)^3\setminus C$. Denote the strategy that $\nu_{_0}$ codes by $\sigma_{_0}$.
  By the assumption there is a strategy $\tau$ of $\PlOne$ which beats $\sigma_{_0}$.
  Consider the game in which $\PlOne$ uses $\tau$ and $\PlTwo$ uses $\sigma_{_0}$. 
  
  Denote the $\g\th$ move in this game
  by $(X_\g,h_\g)$ where $X_\g\subset A_{\eta_{_0}}\cup A_{\xi_{_0}}$ and 
  $h_\g\colon A_{\eta_{_0}}\to A_{\xi_{_0}}$ are the moves of the players. 
  Since player $\PlOne$ wins this game, there is $\a<\l$ for which 
  $h_\a$ is not a partial isomorphism between $\A_{\eta_{_0}}$ and
  $\A_{\xi_{_0}}$. Let
  $$\e=\sup (X_\a\cup\dom h_\a\cup\ran h_\a)$$
  (Recall $\dom \A_\eta=A_\eta=\k$ for any $\eta$ by convention.)
  Let $\pi$ be the coding function defined in Definition \ref{def:CodingOfModels} on 
  page~\pageref{def:CodingOfModels}. Let 
  $$\b_1=\pi[\e^{<\o}]+1.$$
  The idea is that $\eta_{_0}\rest \b_1$ and $\xi_{_0}\rest\b_1$ decide 
  the models $\A_{\eta_{_0}}$ and $\A_{\xi_{_0}}$ as far as
  the game has been played. Clearly $\b_1<\k$.

  Up to this point, player $\PlTwo$ has applied her strategy $\s_{_0}$ precisely to the 
  sequences of the moves made by her opponent, namely to $S=\{(X_\g)_{\g<\b}\mid \b<\a\}\subset \dom \s_{_0}$.
  We can translate this set to represent a subset of the domain of $\nu_{_0}$:
  $S'=k[S]$, where $k$ is as defined before the statement of the present theorem. 
  Let $\b_2=(\sup S')+1$ and let 
  $$\b=\max\{\b_1,\b_2\}.$$
  Thus $\eta_{_0}\rest \b$, $\xi_{_0}\rest\b$ and $\nu_{_0}\rest\b$ decide the moves $(h_\g)_{\g<\a}$ and 
  the winner. 
  
  Now
  \begin{eqnarray*}
    U&=&\{(\nu,\eta,\xi)\mid \nu\rest \b=\nu_{_0}\rest\b\land \eta\rest \b=\eta_{_0}\rest\b\land \xi\rest \b=\xi_{_0}\rest\b\}\\
    &=&N_{\nu_{_0}\restriction\b}\times N_{\eta_{_0}\restriction\b}\times N_{\xi_{_0}\restriction\b}.  
  \end{eqnarray*}
  is the desired neighbourhood. Indeed, if $(\nu,\eta,\xi)\in U$ and $\nu$ codes a strategy $\sigma$, 
  then $\tau$ beats $\sigma$ on the structures
  $\A_\eta,\A_\xi$, since the first $\a$ moves are exactly as in the corresponding game of the triple 
  $(\nu_{_0},\eta_{_0},\xi_{_0})$.
  
  Let us now turn to $D$. The proof is similar. Assume that $(\nu_{_0},\eta_{_0},\xi_{_0})\notin D$ and $\nu_{_0}$
  codes strategy $\tau_{_0}$ of player $\PlOne$. 
  Then there
  is a strategy of $\PlTwo$, which beats $\tau_{_0}$. 
  Let $\b<\k$ be, as before, an ordinal such that
  all moves have occurred before $\b$ and the relations of the substructures generated
  by the moves are decided by $\eta_{_0}\rest\b,\xi_{_0}\rest\b$ as well as the strategy $\tau_{_0}$. 
  Unlike for player $\PlOne$, the win of $\PlTwo$ is determined always only in the end of the game, so
  $\b$ can be $\ge\l$. This is why we made the assumption
  $\l<\k$, by which
  we can always have $\b<\k$ and so 
  \begin{eqnarray*}
    U&=&\{(\nu,\eta,\xi)\mid \nu\rest \b=\nu_{_0}\rest\b\land \eta\rest \b=\eta_{_0}\rest\b\land \xi\rest \b=\xi_{_0}\rest\b\}\\
    &=&N_{\nu_{_0}\restriction\b}\times N_{\eta_{_0}\restriction\b}\times N_{\xi_{_0}\restriction\b}.  
  \end{eqnarray*}
  is an open neighbourhood of $(\nu_{_0},\eta_{_0},\xi_{_0})$ in the complement of~$D$.
\end{proof}

Let us list some theorems concerning Ehrenfeucht-Fra\"iss\'e games which we will use in the proofs.

\begin{Def}\label{def:ScottHeight}
  Let $T$ be a theory and $\A$ a model of $T$ of size $\k$. 
  The $L_{\infty\k}$-\emph{Scott height} of $\A$ is
  $$\sup\{\a\mid \exists\B\models T(\A\not\cong\B\land \PlTwo\uparrow \EF^{\k}_{t_\a}(\A,\B))\},$$
  if the supremum exists and $\infty$ otherwise,
  where $t_\a$ is as in Example \ref{example:Trees} and the subsequent Fact.
\end{Def}

\begin{Remark}
  Sometimes the Scott height is defined in terms of quantifier ranks, but this gives an equivalent
  definition by Theorem~\ref{thm:QuantifierRankAndScottRankTheorem} below.
\end{Remark}

\begin{Def}
  The \emph{quantifier rank} $R(\f)$ of a formula $\f\in L_{\infty\infty}$ 
  is an ordinal defined by induction on the length of $\f$ as follows.
  If $\f$ quantifier free, then $R(\f)=0$. If $\f=\exists\bar x \psi(\bar x)$,
  then $R(\f)=R(\psi(\bar x))+1$. If $\f=\lnot \psi$, then $R(\f)=R(\psi)$.
  If $\f=\Land_{\a<\l}\psi_\a$, then $R(\f)=\sup\{R(\psi_{\a}\mid \a<\l)\}$
\end{Def}

\begin{Thm}\label{thm:QuantifierRankAndScottRankTheorem}
  Models $\A$ and $\B$ satisfy the same $L_{\infty\k}$-sentences of quantifier rank
  $<\a$ if and only if $\PlTwo\uparrow\EF^\k_{t_\a}(\A,\B)$.\qed
\end{Thm}

The following theorem is a well known generalization of a theorem of Karp \cite{Karp}:

\begin{Thm}\label{thm:Karp}
  Models $\A$ and $\B$ are $L_{\infty\k}$-equivalent if and only if $\PlTwo\uparrow \EF^\k_\o(\A,\B)$. \qed
\end{Thm}

\begin{RemarkN}\label{rem:LkkLinfk}
  Models $\A$ and $\B$ of size $\k$ are $L_{\k^+\k}$-equivalent if and only if they are $L_{\infty\k}$-equivalent.
  For an extensive and detailed survey on this and related topics, see~\cite{Vaa3}.
\end{RemarkN}

\subsection{Coding Models} \label{sec:CodingModelss}

There are various degrees of generality to which the content of this text is applicable. 
Many of the results generalize to vocabularies with infinitary relations or to uncountable vocabularies, but not all. 
We find it reasonable though to fix the used vocabulary to make the presentation clearer.

Models can be coded to models with just one binary predicate. 
Function symbols often make situations unnecessarily complicated from the point of view of this paper.

Thus our approach is, without great loss of generality, to 
fix our attention to models with finitary relation symbols of all finite arities.

Let us fix $L$ to be the countable relational vocabulary consisting of the relations $P_n$, $n<\o$, \label{page:DefRelations}
$L=\{P_n\mid n<\o\}$,
where each $P_n$ is an $n$-ary relation: the interpretation of $P_n$ is a set consisting of $n$-tuples. 
We can assume without loss of generality that the domain of each $L$-structure of size $\k$ 
is $\k$, i.e. $\dom \A=\k$. If we restrict our attention to these models, then
the set of all $L$-models has the same cardinality as $\k^\k$.

We will next present the way we code the structures and the isomorphisms between them
into the elements of $\k^\k$ (or equivalently -- as will be seen -- to $2^\k$).

\begin{Def}\label{def:CodingOfModels}
  Let $\pi$ be a bijection $\pi\colon \k^{<\o}\to \k$. 
  If $\eta\in \k^\k$, define the structure $\A_\eta$ to have $\dom(\A_\eta)=\k$ and 
  if $(a_1,\dots a_n)\in \dom(\A_\eta)^n$, then
  $$(a_1,\dots, a_n)\in P_n^{\A_\eta}\iff \eta(\pi(a_1,\dots, a_n))>0.$$
  In that way the rule $\eta \mapsto \A_\eta$ defines a surjective (onto) function from $\k^\k$
  to the set of all $L$-structures with domain $\k$. We say that $\eta$ \emph{codes} $\A_\eta$. \label{page:CodingModels}  
\end{Def}

\begin{Remark}\label{page:DClosed}
  Define the equivalence relation on $\k^\k$ by $\eta\sim \xi\iff \sprt\eta=\sprt\xi$, where $\sprt$
  means support, see Section \ref{ssec:Functions} on page \pageref{ssec:Functions}. 
  Now we have $\eta\sim \xi \iff \A_\eta=\A_\xi$, i.e. the identity map $\k\to\k$ is
  an isomorphism between $\A_\eta$ and $\A_\xi$ when $\eta\sim \xi$ and vice versa. 
  On the other hand $\k^\k/\sim\,\cong 2^{\k}$, so the coding can be seen also as a bijection between
  models and the space $2^{\k}$. 

  The distinction will make little difference, but it is convenient to work with both spaces depending on context.
  To illustrate the insignificance of the choice between $\k^\k$ and $2^\k$, note that $\sim$ is a closed equivalence
  relation and identity on $2^\k$ is bireducible with $\sim$ on $\k^\k$ (see Definition~\ref{def:Reductions}).
\end{Remark}

\subsection{Coding Partial Isomorphisms}
Let $\xi,\eta\in\k^\k$ and let $p$ be a bijection $\k\to\k\times\k$.
Let $\nu\in \k^\a$, $\a\le \k$. The idea is that for $\b<\a$, $p_1(\nu(\b))$ is the image of $\b$ under a partial isomorphism 
and $p_2(\nu(\b))$ is the inverse image of $\b$. 
That is, for a $\nu\in \k^\a$, define a relation $F_\nu\subset \k\times\k$:
$$(\b,\g)\in F_\nu\iff \big(\b<\a\land p_1(\nu(\b))=\g\big)\lor \big(\g<\a\land p_2(\nu(\g))=\b\big)$$
If $\nu$ happens to be such that $F_\nu$ is a partial isomorphism $\A_\xi\to\A_\eta$, then we say 
that $\nu$ \emph{codes a partial isomorphism between} $\A_\xi$ and $\A_\eta$, this 
isomorphism being determined by $F_\nu$. 
If $\a=\k$ and $\nu$ codes a partial isomorphism, then $F_\nu$ is an isomorphism and we say that $\nu$ \emph{codes an isomorphism}. 

\begin{Thm}\label{thm:CodedClosed}
  The set
  $$C=\{(\nu,\eta,\xi)\in (\k^\k)^3\mid \nu\text{ codes an isomorphism between }\A_\eta\text{ and }\A_\xi\}\label{page:Closed}$$
  is a closed set. 
\end{Thm}
\begin{proof}
  Suppose that $(\nu,\eta,\xi)\notin C$ i.e. $\nu$ does not code an isomorphism $\A_\eta\cong\A_\xi$.
  Then (at least) one of the following holds:
  \begin{myEnumerate}
  \item $F_\nu$ is not a function,
  \item $F_\nu$ is not one-to-one,
  \item $F_\nu$ does not preserve relations of $\A_\eta$, $\A_\xi$.
  \end{myEnumerate}
  (Note that $F_\nu$ is always onto if it is a function and $\dom\nu=\k$.)
  If (1), (2) or (3) holds for $\nu$, then
  respectively (1), (2) or (3) holds for any triple $(\nu',\eta',\xi')$
  where $\nu'\in N_{\nu\restriction \g}$, $\eta'\in N_{\eta\restriction \g}$ and $\xi'\in  N_{\xi\restriction \g}$,
  so it is sufficient to check that (1), (2) or (3) holds
  for $\nu\rest\g$ for some $\g<\k$, because 
  
  Let us check the above in the case that (3) holds. The other cases are left to the reader. 
  Suppose (3) holds. There is $(a_0,\dots,a_{n-1})\in (\dom \A_\eta)^n=\k^n$ such that
  $(a_0,\dots,a_{n-1})\in P_n$ and $(a_0,\dots,a_{n-1})\in P_n^{\A_\eta}$ and   
  $(F_\nu(a_0),\dots,F_{\nu}(a_{n-1}))\notin P_n^{\A_\xi}$.
  Let $\b$ be greater than 
  $$\max(\{\pi(a_0,\dots,a_{n-1}),\pi(F_{\nu}(a_0),\dots,F_{\nu}(a_{n-1}))\}\cup \{a_0,\dots a_{n-1},F_{\nu}(a_0),\dots,F_\nu(a_{n-1})\})$$
  Then it is easy to verify that any $(\eta',\xi',\nu')\in N_{\eta\restl\b}\times N_{\xi\restl\b}\times N_{\nu\restl \b}$
  satisfies (3) as well.
\end{proof}

\begin{Cor}\label{cor:IsomisS}
  The set $\{(\eta,\xi)\in (\k^\k)^2\mid \A_\eta\cong \A_\xi\}$ is $\Sii$.
\end{Cor}
\begin{proof}
  It is the projection of the set $C$ of Theorem \ref{thm:CodedClosed}.
\end{proof}

\section{Generalized Borel Sets}
\begin{Def} \label{def:Borel}
  We have already discussed $\Dii$ sets which generalize Borel subsets of Polish space in one way. 
  Let us see how else can we generalize usual Borel sets to our setting. 
  \begin{myItemize}
  \item \cite{Halko, MekVaa} The collection of $\l$-\emph{Borel} subsets of $\k^\k$ is the smallest set, which contains the basic open sets of $\k^\k$ and
    is closed under complementation and under taking intersections of size $\l$. 
    Since we consider only $\k$-Borel sets, we write Borel $=\k$-Borel.
  \item The collection $\Dii=\Sii\cap\Pii$.
  \item \cite{Halko, MekVaa} The collection of \emph{Borel*} subsets of $\k^\k$. A set $A$ is 
    Borel* if there exists a $\k^+\k$-tree $t$ in which each increasing sequence of 
    limit order type has a unique supremum and a function
    $$h\colon \{\text{branches of }t\}\to \{\text{basic open sets of }\k^\k\}$$ 
    such that $\eta\in A\iff$ player $\PlTwo$ has a winning strategy in the game $G(t,h,\eta)$. The game 
    $G(t,h,\eta)$ is defined as follows. At the first round player $\PlOne$ picks a minimal element of the tree, on successive
    rounds he picks an immediate successor of the last move played by player $\PlTwo$
    and if there is no last move, he chooses an immediate successor of the supremum of all previous moves. 
    Player $\PlTwo$ always picks an immediate successor of the Player $\PlOne$'s choice. 
    The game ends when the players cannot go up the tree anymore, i.e. have chosen a branch $b$. 
    Player $\PlTwo$ wins, if $\eta\in h(b)$. Otherwise $\PlOne$ wins. 
    
    A \emph{dual} of a Borel* set $B$ is the set 
    $$B^d=\{\xi\mid \PlOne\uparrow G(t,h,\xi)\}$$ 
    where $t$ and $h$ satisfy the equation $B=\{\xi\mid \PlTwo\uparrow G(t,h,\xi)\}$. The dual is not unique.
  \end{myItemize}
\end{Def}

\begin{Remark}
  Suppose that $t$ is a $\k^+\k$ tree and $h\colon \{\text{branches of }t\}\to \text{Borel}^*$ is a labeling function
  taking values in Borel* sets instead of basic open sets. Then $\{\eta\mid \PlTwo\uparrow G(t,h,\eta)\}$
  is a Borel* set.   

  Thus if we change the basic open sets to Borel* sets in the definition of Borel*, we get Borel*. 
\end{Remark}

\begin{RemarkN}\label{remark:Borel}
  Blackwell \cite{Bla} defined Borel* sets in the case $\k=\o$ and showed that in fact Borel=Borel*. When 
  $\k$ is uncountable it is not the case.
  But it is easily seen that if $t$ is a $\k^+\o$-tree, then the Borel* set coded by $t$ (with some labeling $h$)
  is a Borel set, and vice versa: each Borel set is a Borel* set coded by a $\k^+\o$-tree.
  We will use this characterization of Borel.
\end{RemarkN}

It was first explicitly proved in \cite{MekVaa} that these are indeed generalizations:

\begin{Thm}[\cite{MekVaa}, $\k^{<\k}=\k$]\label{thm:BorelDiiBorelStar}
  Borel $\subset\Dii\subset$ Borel* $\subset\Sii$, 
\end{Thm}
\begin{proof}(Sketch)
  If $A$ is Borel*, then it is $\Sii$, intuitively, because
  $\eta\in A$ if and only if \emph{there exists} a winning strategy of player $\PlTwo$ in
  $G(t,h,\eta)$ where $(t,h)$ is a tree that codes $A$ (here one needs the assumption $\k^{<\k}=\k$
  to be able to code the strategies into the elements of $\k^\k$). By Remark \ref{remark:Borel} above if $A$ is Borel, then there is also
  such a tree. Since Borel $\subset$ Borel* by Remark \ref{remark:Borel} and Borel is closed under taking complements, Borel sets are $\Dii$. 

  The fact that $\Dii$ sets are Borel* is a more complicated issue; it follows from a separation theorem proved in \cite{MekVaa}. 
  The separation theorem says that any two disjoint $\Sii$ sets can be separated by Borel* sets. It is proved in \cite{MekVaa} for $\k=\o_1$,
  but the proof generalizes to any $\k$ \mbox{(with $\k^{<\k}=\k$)}. 
\end{proof}

Additionally we have the following results:

\begin{Thm}
  \begin{myEnumerate}
  \item ${\rm Borel}\subsetneq \Dii$.
  \item $\Dii\subsetneq \Sii$.
  \item If $V=L$, then ${\rm Borel}^*=\Sii$. 
  \item It is consistent that $\Dii\subsetneq{\rm Borel}^*$. 
  \end{myEnumerate}
\end{Thm}
\begin{proof}(Sketch)
  \begin{myEnumerate}
  \item  The following universal Borel set is not Borel itself, but is~$\Dii$:
    $$B=\{(\eta,\xi)\in 2^\k\times 2^\k\mid \eta\text{ is in the set coded by }(t_\xi,h_\xi)\},$$
    where $\xi\mapsto (t_\xi,h_\xi)$ is a continuous coding of 
    ($\k^+\o$-tree, labeling)-pairs in such a way that for all $\k^+\o$-trees $t\subset \k^{<\o}$ and labelings $h$ there is
    $\xi$ with $(t_\xi,h_\xi)=(t,h)$.
    It is not Borel since if it were, then the diagonal's complement
    $$D=\{\eta\mid (\eta,\eta)\notin B\}$$
    would be a Borel set which it is not, since it cannot be coded by any $(t_\xi,h_\xi)$. On the other hand its
    complement $C=(2^\k)^2\setminus B$ is $\Sii$, because $(\eta,\xi)\in C$ if and only if
    \emph{there exists} a winning strategy of player $\PlOne$ in the Borel-game $G(t_\xi,h_\xi,\eta)$ and the latter can be coded to a Borel set.
    It is left to the reader to verify that when $\k>\o$, then the set
    $$F=\{(\eta,\xi,\nu)\mid \nu\text{ codes a w.s. for }\PlOne\text{ in }G(t_\xi,h_\xi,\eta)\}$$
    is closed. 
    
    The existence of an isomorphism relation which is $\Dii$ but not 
    Borel follows from Theorems \ref{thm:ClasNotShalIsNotBorel} and~\ref{thm:Dii}.
  \item Similarly as above (and similarly as in the case $\k=\o$), take a universal $\Sii$-set $A\subset 2^{\k}\times 2^{\k}$ 
    with the property that if $B\subset 2^\k$ is any $\Sii$-set, then there is $\eta\in 2^\k$ such that 
    $B\times \{\eta\}\subset A$. This set can be constructed as in the case $\k=\o$, see \cite{Jech}. The
    diagonal $\{\eta\mid (\eta,\eta)\in A\}$ is $\Sii$ but not~$\Pii$.
  \item Suppose $V=L$ and $A\subset 2^\k$ is $\Sii$. There exists a formula $\f(x,\xi)$ with parameter $\xi\in 2^{\k}$ 
    which is $\Sigma_1$ in the Levy hierarchy
    (see \cite{Jech}) and for all $\eta\in 2^\k$ we have 
    $$\eta\in A\iff L\models \f(\eta,\xi)$$
    Now we have that 
    $\eta \in A$ if and only if the set 
    $$\big\{\a<\k \mid \exists\b\big(\eta\rest\a,\xi\rest\a\in L_\beta,\ L_{\b}\models\big({\rm ZF}^-\land (\alpha \text{ is a cardinal})\land 
      \f(\eta\rest\alpha,\xi\rest\a)\big)\big)\big\}$$ 
    contains an $\o$-cub set.
    
    But the $\omega$-cub filter is Borel* so $A$ is also Borel*. 
  \item This follows from the clauses \eqref{thm:CUBSETItem1}, \eqref{thm:CUBSETItem5} and \eqref{thm:CUBSETItem6} of 
    Theorem \ref{thm:CUBSET} below.\qedhere
  \end{myEnumerate}
\end{proof}

\begin{Open}
  Is it consistent that Borel* is a proper subclass of $\Sii$, or even equals $\Dii$? Is it consistent that
  all the inclusions are proper at the same time: $\Dii\subsetneq {\rm Borel}^*\subsetneq \Sii$?
\end{Open}

\begin{Thm}\label{thm:Proj}
  For a set $S\subset \k^\k$ the following are equivalent.
  \begin{myEnumerate}
  \item $S$ is $\Sii$,
  \item $S$ is a projection of a Borel set,
  \item $S$ is a projection of a $\Sii$ set,
  \item $S$ is a continuous image of a closed set.
  \end{myEnumerate} 
\end{Thm}
\begin{proof}
  Let us go in the order.
  \begin{myDescription}
  \item[$(1)\Rightarrow(2)$] Closed sets are Borel.
  \item[$(2)\Rightarrow(3)$] The same proof as in the standard case $\k=\o$ gives that Borel sets are $\Sii$ (see for instance \cite{Jech}).
  \item[$(3)\Rightarrow(4)$] Let $A\subset \k^\k\times \k^\k$ be a $\Sii$ set which is the projection of $A$, 
    $S=\pr_0 A$. Then let $C\subset \k^\k\times \k^\k\times \k^\k$ be a closed set such that $\pr_1 C = A$.
    Here $\pr_0\colon \k^\k\times \k^\k\to \k^\k$ and $\pr_1\colon \k^\k\times \k^\k\times \k^\k\to \k^\k\times \k^\k$ are the obvious projections. 
    Let $f\colon \k^\k\times \k^\k\times \k^\k\to \k^\k$
    be a homeomorphism. Then $S$ is the image of the closed set $f[C]$ under the continuous map $\pr_0\circ\pr_1\circ f^{-1}$.
  \item[$(4)\Rightarrow(1)$] The image of a closed set under a continuous map $f$ is the projection of the graph of $f$ 
    restricted to that closed set. It is a basic topological fact that a graph of a continuous partial function with closed domain
    is closed (provided the range is Hausdorff).
  \end{myDescription}
\end{proof}

\begin{Thm}[\cite{MekVaa}]\label{thm:BorelStarClosed}
  Borel* sets are closed under unions and intersections of size $\k$.\qed
\end{Thm}

\begin{Def}
  A Borel* set $B$ is \emph{determined} if there exists a tree $t$ and a labeling function $h$ such that
  the corresponding game $G(t,h,\eta)$ is determined for all $\eta\in \k^\k$ and 
  $$B=\{\eta\mid \PlTwo\text{ has a winning strategy in }G(t,h,\eta)\}.$$
\end{Def}

\begin{Thm}[\cite{MekVaa}]\label{thm:DiiIsBorelStar}
  $\Dii$ sets are exactly the determined Borel* sets.\qed
\end{Thm}

\chapter{Borel Sets, $\Dii$ Sets and Infinitary Logic}
\section{The Language $L_{\k^+\k}$ and Borel Sets}

The interest in the class of Borel sets is explained by the fact that
the Borel sets are relatively simple yet at the same time this class
includes many interesting definable sets. We prove Vaught's theorem (Theorem \ref{thm:BorelIsLkk}),
which equates ``invariant'' Borel sets with those definable in the
infinitary language $L_{\kappa^+\kappa}$.
Recall that models $\A$ and $\B$ of size $\k$ are $L_{\k^+\k}$-equivalent if and only if they are $L_{\infty\k}$-equivalent.
 Vaught proved his theorem for the case $\kappa = \omega_1$ assuming CH in
\cite{Vau}, but the proof works for arbitrary $\kappa$ assuming $\kappa^{<\kappa} = \kappa$.

\begin{Def}\label{def:Permutation}
  Denote by $S_\k$ the set of all permutations of $\k$. 
  If $u\in \k^{<\k}$, denote 
  $$\bar u=\{p\in S_\k\mid p^{-1}\rest\dom u = u\}.$$
  Note that $\bar \es=S_\k$ and if $u\in \k^\a$ is not injective, then $\bar u=\es$.
  
  A permutation $p\colon \k\to\k$ acts on $2^{\k}$ by
  $$p\eta=\xi\iff p\colon \A_\eta\to\A_\xi\text{ is an isomorphism}.$$
  The map $\eta\mapsto p\eta$ is well defined for every $p$ and it is easy to 
  check that it defines an action of the permutation group $S_\k$ on the space $2^{\k}$.
  We say that a set $A\subset 2^{\k}$ is \emph{closed under permutations} if it is a union
  of orbits of this action.
\end{Def}

\begin{Thm}[\cite{Vau}, $\k^{<\k}=\k$]\label{thm:BorelIsLkk}
  A set $B\subset \k^\k$ is Borel and closed under permutations if and only if there is
  a sentence $\f$ in $L_{\k^+\k}$ such that $B=\{\eta\mid \A_\eta\models \f\}$.
\end{Thm}
\begin{proof}
  Let $\f$ be a sentence in $L_{\k^+\k}$. 
  Then $\{\eta\in 2^{\k}\mid \A_\eta\models \f\}$
  is closed under permutations, because if $\eta=p\xi$, then
  $\A_\eta\cong \A_{\xi}$ and $\A_\eta\models\f\iff \A_\xi\models \f$ for every sentence $\f$.
  If $\f$ is a formula with parameters $(a_{i})_{i<\a}\in \k^\a$,
  one easily verifies by induction on the complexity of $\f$ that
  the set 
  $$\{\eta\in 2^{\k}\mid \A_\eta\models\f((a_{i})_{i<\a})\}$$ 
  is Borel. This of course implies that for every sentence $\f$ the set $\{\eta\mid \A_\eta\models \f\}$
  is Borel.

  The converse is less trivial.
  Note that the set of permutations $S_\k\subset \k^\k$ is Borel, since
  $$S_\k = \Cap_{\b<\k}\Cup_{\a<\k}\underbrace{\{\eta\mid \eta(\a)=\b\}}_{\text{open}}\cap  
  \Cap_{\a<\b<\k}\underbrace{\{\eta\mid \eta(\a)\ne\eta(\b)\}}_{\text{open}}.\eqno(\cdot)$$\label{page:PermBorel}
  For a set $A\subset \k^\k$ and $u\in \k^{<\k}$, define 
  $$A^{*u}=\big\{\eta\in 2^{\k}\mid \{p\in \bar u\mid p\eta\in A\} \text{ is co-meager in }\bar u\big\}.$$
  From now on in this section we will write ``$\{p\in \bar u\mid p\eta\in A\}$ is co-meager'', when we really mean ``co-meager in $\bar u$''.

  Let us show that the set 
  $$Z=\{A\subset 2^{\k}\mid A\text{ is Borel}, A^{*u} \text{ is }L_{\k^+\k}\text{-definable for all }u\in \k^{<\k}\}$$  
  contains all the basic open sets, is closed under intersections of size $\k$ and under
  complementation in the three steps (\emph{a}),(\emph{b}) and (\emph{c}) below. This implies that $Z$ is the collection of all Borel sets. 
  We will additionally keep track
  of the fact that the formula, which defines $A^{*u}$ depends only on $A$ and $\dom u$, i.e. 
  for each $\b<\k$ and Borel set $A$ 
  there exists $\f=\f^A_\b$ such that for all $u\in \k^\b$ we have $A^{*u}=\{\eta \mid \A_\eta\models \f((u_i)_{i<\b})\}$.
  Setting $u=\es$, we have the intended result, because $A^{*\es}=A$ for all $A$ which are closed under permutations
  and $\f$ is a sentence (with no parameters).  
  
  If $A$ is fixed we denote $\f^A_\b=\f_\b$.
  \begin{myItemize}
  \item[(\emph{a})] Assume $q\in 2^{<\k}$ and let $N_q$ be the corresponding basic open set. Let us show that
    $N_q\in Z$. Let $u\in \k^\b$ be arbitrary. We have to find $\f^{N_q}_\b$. 
    Let $\theta$ be a quantifier free formula with $\a$ parameters such that:
    $$N_q=\{\eta\in 2^{\k}\mid \A_\eta\models \theta((\g)_{\g<\a})\}.$$
    Here $(\g)_{\g<\a}$ denotes both an initial segment of $\k$ as well as an $\a$-tuple
    of the structure.
    Suppose $\a\le\b$. We have
    $p\in \bar u\Rightarrow  u\subset p^{-1}$, so 
    \begin{eqnarray*}
      \eta\in N_q^{*u}&\iff& \{p\in \bar u\mid p\eta\in N_q\} \text{ is co-meager}\\
      &\iff& \{p\in \bar u\mid \A_{p\eta}\models \theta((\g)_{\g<\a})\} \text{ is co-meager}\\
      &\iff& \{p\in \bar u\mid \A_{\eta}\models \theta((p^{-1}(\g))_{\g<\a})\} \text{ is co-meager}\\
      &\iff& \{p\in \bar u\mid \underbrace{\A_{\eta}\models \theta((u_{\g})_{\g<\a})}_{\text{independent of }p}\}\text{ is co-meager}\\
      &\iff& \A_{\eta}\models \theta((u_{\g})_{\g<\a}).
    \end{eqnarray*}
    Then $\f_\b=\theta$. 

    Assume then that $\a>\b$. By the above, we still have
    $$\eta\in N_q^{*u}\iff E=\big\{p\in \bar u\mid \A_{\eta}\models \theta\big( (p^{-1}(\g))_{\g<\a}\big)\big\} \text{ is co-meager}$$
    Assume that $w=(w_{\g})_{\g<\a}\in \k^{\a}$ is an arbitrary sequence with no repetition and such that
    $u\subset w$. Since $\bar w$ is an open subset of $\bar u$ and $E$ is co-meager, there is $p\in \bar w\cap E$.
    Because $p\in E$, we have $\A_\eta\models\theta\big( (p^{-1}(\g))_{\g<\a}\big)$. On the other hand $p\in \bar w$, so
    we have $w\subset p^{-1}$, i.e. $w_\g=w(\g)=p^{-1}(\g)$ for $\g<\a$. Hence
    $$\A_\eta\models \theta((w_{\g})_{\g<\a}).\eqno(\star)$$ 
    On the other hand, if
    for every injective $w\in\k^\a$, $w\supset u$, we have $(\star)$, then 
    in fact $E=\bar u$ and is trivially co-meager. Therefore we have an equivalence:
    $$\eta\in N_q^{*u}\iff(\forall w\supset u)(w\in \k^\a\land w\text{ inj.}\Rightarrow\A_\eta\models \theta((w_{\g})_{\g<\a})).$$ 
    But the latter can be expressed
    in the language $L_{\k^+k}$ by the formula $\f_\b((w_i)_{i<\b})$:
    $$\Land_{i<j<\b}(w_i\ne w_j)\land \Big(\Forall_{\b\le i<\a}w_i\Big)\Big(\Land_{i<j<\a}(w_i\ne w_j)\rightarrow\theta((w_i)_{i<\a})\Big)$$
    $\theta$ was defined to be a formula defining $N_q$ with parameters. It is clear thus that $\theta$
    is independent of $u$. Furthermore the formulas constructed above from $\theta$ depend only
    on $\b=\dom u$ and on $\theta$. Hence the formulas defining $N_q^{*u}$ and $N_q^{*v}$ for 
    $\dom u=\dom v$ are the same modulo parameters.
  \item[(\emph{b})] For each $i< \k$ let $A_i\in Z$. We want to show that $\Cap_{i<\k}A_i\in Z$. 
    Assume that $u\in \k^{<\k}$ is arbitrary. It suffices to show that 
    $$\Cap_{i<\k}(A_i^{*u})=\Big(\Cap_{i<\k}A_i\Big)^{*u},$$
    because then $\f^{\cap_iA_i}_\b$ is just the $\k$-conjunction of the formulas $\f^{A_i}_\b$ which exist 
    by the induction hypothesis. Clearly the resulting formula depends again only on $\dom u$ if the previous did.
    Note that a $\k$-intersection of co-meager sets is co-meager. Now
    \begin{eqnarray*}
      &&     \eta\in \Cap_{i<\k}(A_i^{*u})\\
      &\iff& (\forall i<\k) (\{p\in \bar u\mid p\eta\in A_i\} \text{ is co-meager})\\
      &\iff& (\forall i<\k)(\forall i<\k)(\{p\in \bar u\mid p\eta\in A_i\} \text{ is co-meager})\\
      &\iff& \Cap_{i<\k}\{p\in \bar u\mid p\eta\in A_i\} \text{ is co-meager}\\
      &\iff& \{p\in \bar u\mid p\eta\in \Cap_{i<\k}A_i\} \text{ is co-meager}\\
      &\iff& \eta\in \Big(\Cap_{i<\k}A_i\Big)^{*u}.
    \end{eqnarray*}
  \item[(\emph{c})] Assume that $A\in Z$ i.e. that $A^{*u}$ is definable for any $u$. Let
    $\f_{\dom u}$ be the formula, which defines $A^{*u}$. Let now $u\in \k^{<\k}$ be arbitrary and
    let us show that $(A^c)^{*u}$ is definable. We will show that
    $$(A^c)^{*u}=\Cap_{v\supset u} \left(A^{*v}\right)^c$$
    i.e. for all $\eta$ 
    \begin{equation}
      \label{eq:BP}
      \eta\in (A^c)^{*u}\iff \forall v\supset  u (\eta\notin A^{*v}).  
    \end{equation}
    Granted this, one can write the formula ``$\forall v\supset u \lnot\f_{\dom u}((v_i)_{i<\dom v})$'', 
    which is not of course the real $\f^{A^c}_\b$ which we will write in the end of the proof.
    
    To prove \eqref{eq:BP}
    we have to show first that for all $\eta\in \k^\k$ the set $B=\{p\in \bar u\mid p\eta\in A\}$ has the Property of Baire
    (P.B.), see Section \ref{sec:PB}.
    
    The set of all permutations $S_\k\subset \k^\k$ is Borel by $(\cdot)$ on page \pageref{page:PermBorel}.
    The set $\bar u$ is an intersection of $S_\k$ with an open set. Again the 
    set $\{p\in \bar u\mid p\eta\in A\}$ is the intersection of $\bar u$ and 
    the inverse image of $A$ under the continuous map $(p\mapsto p\eta)$, so is Borel and
    so has the Property of Baire.
    
    We can now turn to proving the equivalence \eqref{eq:BP}. First ``$\Leftarrow$'': 
    \begin{eqnarray*}
      \eta\notin (A^c)^{*u}&\Rightarrow& B=\{p\in \bar u\mid p\eta\in A\} \text{ is not meager in }\bar u\\
      &\Rightarrow& \text{By P.B. of $B$ there is a non-empty open $U$ such that }U\setminus B\text{ is meager}\\
      &\Rightarrow& \text{There is non-empty }\bar v\subset \bar u \text{ such that }\bar v\setminus B\text{ is meager.}\\
      &\Rightarrow& \text{There exists }\bar v\subset \bar u\text{ such that } \{p\in \bar v\mid p\eta\in A\}=\bar v\cap B \text{ is co-meager}\\
      &\Rightarrow& \exists v\supset u(\eta\in A^{*v}).
    \end{eqnarray*}
    And then the other direction ``$\Rightarrow$'':
    \begin{eqnarray*}
      \eta\in (A^c)^{*u}&\Rightarrow& \{p\in \bar u\mid p\eta\in A\} \text{ is meager}\\
      &\Rightarrow& \text{ for all }\bar v\subset \bar u\text{ the set }\{p\in \bar v\mid p\eta\in A\}\text{ is meager}.\\
      &\Rightarrow& \forall \bar v\subset \bar u (\eta\notin A^{*v}).
    \end{eqnarray*}
  \end{myItemize}
  Let us now write the formula $\psi=\f^{A^c}_\b$ such that 
  $$\forall \bar v\subset \bar u(\eta\notin A^{*v})\iff \A_\eta\models \psi((u_i)_{i<\b}),$$
  where $\b=\dom u$: let $\psi((u_i)_{i<\b})$ be
  $$\Land_{\b\le\g<\k}\Forall_{i<\g}x_i\left(\Big[\Land_{j<\b}(x_j=u_j)\land \Land_{i<j<\g}(x_i\ne x_j)\Big]\rightarrow \lnot\f_\g((x_i)_{i<\g})\right)$$
  One can easily see, that this is equivalent to $\forall v\supset u\big(\lnot \f_{\dom v}((v_i)_{i<\dom v})\big)$
  and that $\psi$ depends only on $\dom u$ modulo parameters.
\end{proof}

\begin{RemarkN}\label{separationremark}
  If $\k^{<\k}>\k$, then the direction from right to left of the above theorem does not in general hold. Let $\la \k,\less,A\ra$ be
  a model with domain $\k$, $A\subset \k$ and $\less$ a well ordering of $\k$ of order type $\k$.
  V\"a\"an\"anen and Shelah have shown in \cite{SheVaa} (Corollary 17)
  that if $\k=\l^+$, $\k^{<\k}>\k$, $\l^{<\l}=\l$ and a forcing axiom holds (and $\o_1^L=\o_1$ if $\l=\o$)
  then there is a sentence of $L_{\k\k}$ defining the set

  \newcommand{\Stat}{\operatorname{STAT}}
  $$\Stat=\{\la \k,\less,A\ra\mid A\text{ is stationary}\}.$$
  If now $\Stat$ is Borel, then so would be the set $\CUB$ defined in Section \ref{sec:PB},
  but by Theorem \ref{thm:CUBSET} this set cannot be Borel since Borel sets have the Property of Baire by Theorem~\ref{thm:BorelNoBaire}.
\end{RemarkN}

\begin{Open}
  Does the direction left to right of Theorem \ref{thm:BorelIsLkk} hold without the assumption
  $\k^{<\k}=\k$?
\end{Open}

\section{The Language $M_{\k^+\k}$ and $\Dii$-Sets}
In this section we will present a theorem similar to Theorem \ref{thm:BorelIsLkk}. 
It is also a generalization of the known result which follows from \cite{MekVaa} and \cite{Vaa2}:

\begin{Thm}[\cite{MekVaa,Vaa2}:]
  Let $\A$ be a model of size  $\o_1$. Then the isomorphism type $I=\{\eta\mid \A_\eta\cong\A\}$ is $\Dii$ if and only if
  there is a sentence $\f$ in $M_{\k^+\k}$ such that $I=\{\eta\mid \A_\eta\models \f\}$ and 
  $2^\k\setminus I=\{\eta\mid \A_\eta\models \sim\f\}$, where $\sim\theta$ is the dual of $\theta$.  
\end{Thm}

The idea of the proof of the following Theorem
is due to Sam Coskey and Philipp Schlicht:

\begin{Thm}[$\k^{<\k}=\k$]\label{thm:DiiIsMkk}
  A set $D\subset 2^\k$ is $\Dii$ and closed under permutations if and only if there is
  a sentence $\f$ in $M_{\k^+\k}$ such that $D=\{\eta\mid \A_\eta\models \f\}$ and 
  $\k^\k\setminus D=\{\eta\mid \A_\eta\models \sim\f\}$, where $\sim\theta$ is the dual of $\theta$.  
\end{Thm}
We have to define these concepts before the proof.
\begin{Def}[Karttunen \cite{Karttunen}]
  Let $\l$ and $\k$ be cardinals. The language 
  $M_{\l\k}$ is then defined to be the set of pairs $(t,\ll)$ of a tree $t$ and a labeling function $\ll.$
  The tree $t$ is a $\l\k$-tree where
  the limits of increasing sequences of $t$ exist and are unique.
  The labeling $\ll$ is a function satisfying the following conditions:
  \begin{myEnumerate}
  \item  $\ll\colon t\to a\cup \bar a\cup\{\Land,\Lor\}\cup \{\exists x_i\mid i<\k\}\cup \{\forall x_i\mid i<\k\}$ where $a$ is the set of atomic 
    formulas and $\bar a$ is the set of negated atomic formulas.
  \item If $x\in t$ has no successors, then $\ll(t)\in a\cup \bar a$.
  \item If $x\in t$ has exactly one immediate successor then $\ll(t)$ is either $\exists x_i$ or $\forall x_i$ for some $i<\k$.
  \item Otherwise $\ll(t)\in\{\Lor,\Land\}$.
  \item If $x<y$, $\ll(x)\in \{\exists x_i,\forall x_i\}$ and $\ll(y)\in \{\exists x_j,\forall x_j\}$, then $i\ne j$.
  \end{myEnumerate}
\end{Def}

\begin{Def}
  Truth for $M_{\l\k}$ is defined in terms of a semantic game. Let $(t,\ll)$ be the pair which
  corresponds to a particular sentence $\f$ and let $\A$ be a model. The semantic game $S(\f,\A)=S(t,\ll,\A)$ for $M_{\l\k}$ 
  is played by players $\PlOne$ and $\PlTwo$
  as follows. At the first move the players are at the root and later in the game at some other element of $t$. 
  Let us suppose that they are at the element $x\in t$. If $\ll(x)=\Lor$, then Player $\PlTwo$ chooses a successor of $x$ and the players move to that
  chosen element. If $\ll(x)=\Land$, then player $\PlOne$ chooses a successor of $x$ and the players move to that chosen element.
  If $\ll(x)=\forall x_i$ then player $\PlOne$ picks an element $a_i\in\A$ and if $\ll(x)=\exists x_i$ then player $\PlTwo$ picks 
  an element $a_i$
  and they move to the immediate successor of $x$. If they come to a limit, they move to the unique supremum. If $x$ is
  a maximal element of $t$, then they plug the elements $a_i$ in place of the corresponding free variables in the atomic formula
  $\ll(x)$. Player $\PlTwo$ wins if this atomic 
  formula is true in $\A$ with these interpretations. Otherwise player $\PlOne$ wins.

  We define $\A\models\f$ if and only if $\PlTwo$ has a winning strategy in the semantic game.

  Given a sentence $\f$, the sentence $\sim\f$ is defined by modifying the labeling function as follows. The atomic formulas are replaced by
  their negations, the symbols $\Lor$ and $\Land$ switch places and the quantifiers $\forall$ and $\exists$ switch places.
  A sentence $\f\in M_{\l\k}$ is \emph{determined} if for all models $\A$ either $\A\models\f$ or $\A\models\sim\f$.
\end{Def}

Now the statement of Theorem \ref{thm:DiiIsMkk} makes sense. 
Theorem \ref{thm:DiiIsMkk} concerns a sentence $\f$ whose dual defines the complement of the set defined by $\f$ among the models of size $\k$, 
so it is determined in that model class.
Before the proof let us 
recall a separation theorem for $M_{\k^+\k}$, Theorem 3.9 from \cite{Tuuri}:
\begin{Thm}\label{thm:SepMkk}
  Assume $\k^{<\k}=\l$ and let $\exists R\f$ and $\exists S\psi$ be two $\Sii$ sentences where $\f$ and $\psi$ are in
  $M_{\k^+\k}$ and $\exists R$ and $\exists S$ are second order quantifiers. If 
  $\exists R\f\land \exists S\psi$ does not have a model, then there is a sentence $\theta\in M_{\l^+\l}$ 
  such that
  for all models $\A$
  $$\A\models\exists R\f\Rightarrow \A\models\theta\text{ and }\A\models\exists S\psi\Rightarrow \A\models \sim\theta\eqno\qed$$
\end{Thm}

\begin{Def}\label{def:sigmaOfTree}
  For a tree $t$, let $\s t$ be the tree of downward closed linear subsets of $t$ ordered by inclusion.
\end{Def}

\begin{proofVOf}{Theorem \ref{thm:DiiIsMkk}}
  Let us first show that if $\f$ is an arbitrary sentence of $M_{\k^+\k}$, then $D_{\f}=\{\eta\mid\A_\eta\models \f\}$
  is $\Sii$. The proof has the same idea as the proof of Theorem \ref{thm:BorelDiiBorelStar} that Borel* $\subset\Sii$.
  Note that this implies that if $\sim\f$ defines the complement of $D_\f$ in $2^\k$, then $D_\f$ is $\Dii$. 

  A strategy in the semantic game $S(\f,\A_\eta)=S(t,\ll,\A_\eta)$ is a function 
  $$\upsilon\colon \s t\times (\dom\A_\eta)^{<\k}\to t\cup (t\times\dom\A_\eta).$$ 
  This is because the previous moves always form an initial segment of a branch of the tree together with the 
  sequence of constants picked by the players
  from $\dom\A_\eta$ at the quantifier moves, and a move consists either of going to some node of the tree or
  going to a node of the tree together with choosing an element from $\dom \A_\eta$. 
  By the convention that $\dom\A_\eta=\k$, a strategy becomes a function 
  $$\upsilon\colon \s t\times \k^{<\k}\to t\cup (t\times\k),$$ 
  
  Because $t$ is a $\k^+\k$-tree, there are fewer than $\k$ moves in a play (there are no branches of length $\k$ and the players go up the tree
  on each move).
  Let 
  $$f\colon \s t\times \k^{<\k}\to \k$$ 
  be any bijection and let 
  $$g\colon t\cup (t\times \k)\to \k$$ 
  be another bijection. Let 
  $F$ be the bijection 
  $$F\colon (t\cup (t\times\k))^{\s t\times \k^{<\k}}\to \k^\k$$
  defined by $F(\upsilon)=g\circ\upsilon\circ f^{-1}$. Let 
  $$C=\{(\eta,\xi)\mid F^{-1}(\xi) \text{ is a winning strategy of }\PlTwo\text{ in }S(t,\ll,\A_\eta)\}.$$
  Clearly $D_\f$ is the projection of $C$. Let us show that $C$ is closed. Consider an element $(\eta,\xi)$ in the complement
  of $C$. We shall show that there is an open neighbourhood of $(\eta,\xi)$ outside $C$. 
  Denote $\upsilon=F^{-1}(\xi)$. Since $\upsilon$ is not a winning strategy there is a strategy $\tau$ of $\PlOne$ that beats $\upsilon$.
  There are $\a+1<\k$ moves in the play $\tau*\upsilon$ (by definition all branches have
  successor order type). Assume that $b=(x_i)_{i\le\a}$ is the chosen branch of the tree and 
  $(c_i)_{i<\a}$ the constants picked by the players.
  Let $\b<\k$ be an ordinal with the properties $\{f((x_{i})_{i<\g},(c_i)_{i<\g})\mid \g\le\a+1\}\subset\b$ and
  $$\eta'\in N_{\eta\restriction\b}\rightarrow \A_{\eta'}\not\models \ll(x_\a)((c_{i})_{i<\a}).\eqno(\star)$$
  Such $\b$ exists, since $|\{f((x_{i})_{i<\g},(c_i)_{i<\g})\mid \g\le\a+1\}|<\k$ and $\ll(x_\a)$ is a (possibly negated) atomic formula 
  which is not true in $\A_{\eta}$, because $\PlTwo$ lost the game $\tau*\upsilon$ and because already a 
  fragment of size $<\k$ of $\A_\eta$ decides this.
  Now if $(\eta',\xi')\in N_{\eta\restriction\b}\times N_{\xi\restriction \b}$
  and $\upsilon'=F^{-1}(\xi')$, then $\upsilon*\tau$ is the same play as $\tau*\upsilon'$. 
  So $\A_{\eta'}\not\models \ll(x_\a)((c_i)_{i<\a})$ by $(\star)$ and  $(\eta',\xi')$ is not in $C$ and  
  $$N_{\eta\restriction \b}\times N_{\xi\restriction \b}$$
  is the intended open neighbourhood of $(\eta,\xi)$ outside $C$. 
  This completes the ``if''-part of the proof.

  Now for a given $A\in \Dii$ which is closed under permutations we want to find a sentence $\f\in M_{\k^+\k}$ 
  such that $A=\{\eta\mid\A_{\eta}\models \f\}$ and $2^\k\setminus A=\{\eta\mid\A_{\eta}\models \sim\f\}$. 
  By our assumption $\k^{<\k}=\k$ and 
  Theorems \ref{thm:DiiIsBorelStar} and \ref{thm:SepMkk}, it is enough to show
  that for a given Borel* set $B$ which is closed under permutations, there is a sentence $\exists R\psi$ 
  which is $\Sii$ over $M_{\k^+\k}$ (as in the formulation of Theorem \ref{thm:SepMkk}), such that
  $B=\{\eta\mid \A_\eta\models \exists R\psi\}$.

  The sentence ``$R$ is a well ordering of the 
  universe of order type $\k$'', is definable by the formula $\theta=\theta(R)$ of $L_{\k^+\k}\subset M_{\k^+\k}$:
  \begin{eqnarray}
    &&"R\text{ is a linear ordering on the universe}"\nonumber\\
    &\land& \Big(\Forall_{i<\o}x_i\Big) \Big(\Lor_{i<\o}\lnot R(x_{i+1},x_i)\Big)\nonumber\\
    &\land& \forall x \Lor_{\a<\k}\Exists_{i<\a}y_i\left[\big(\forall y(R(y,x)\rightarrow \Lor_{i<\a}y_i=y)\big)\right]\label{eq:WO}
  \end{eqnarray}
  (We assume $\k>\o$, so the infinite quantification is allowed. The second row says that there are no descending sequences of length $\o$ and
  the third row says that the initial segments are of size less than $\k$. This ensures that $\theta(R)$ says that $R$ is a well ordering
  of order type $\k$). 

  Let $t$ and $h$ be the tree and the labeling function corresponding to $B$. Define the tree $t^\star$ as follows.
  \begin{myEnumerate}
  \item Assume that $b$ is a branch of $t$ with $h(b)=N_{\xi\restriction \a}$ for some $\xi\in\k^\k$ and $\a<\k$. Then
    attach a sequence of order type $\a^*$ on top of $b$ where 
    $$\a^*=\Cup_{s\in\pi^{-1}[\a]} \ran s,$$
    where $\pi$ is the bijection $\k^{<\o}\to\k$ used in the coding, see Definition \ref{def:CodingOfModels} 
    on page \pageref{def:CodingOfModels}.
  \item Do this to each branch of $t$ and add a root $r$ to the resulting tree.
  \end{myEnumerate}
  After doing this, the resulting tree is $t^\star$. Clearly it is a $\k^+\k$-tree, because $t$ is. 
  Next, define the labeling function $\ll$.
  If $x\in t$ then either $\ll(x)=\Land$ or $\ll(x)=\Lor$ depending on 
  whether it is player $\PlOne$'s move or player $\PlTwo$'s move:
  formally let $n<\o$ be such that $\OTP(\{y\in t^\star\mid y\le x\})=\a+n$ where $\a$ is a limit ordinal or $0$; then if 
  $n$ is odd, put $\ll(x)=\Land$ and otherwise $\ll(x)=\Lor$. 
  If $x=r$ is the root, then $\ll(x)=\Land$. Otherwise, if $x$ is not maximal, define
  $$\b=\OTP\{y\in t^\star\setminus (t\cup \{r\})\mid y\le x\}$$ 
  and set $\ll(x)=\exists x_{\b}$. 

  Next we will define the labeling of the maximal nodes of $t^\star$. By definition
  these should be atomic formulas or negated atomic formulas, but it is 
  clear that they can be replaced without loss of generality 
  by any formula of $M_{\k^+\k}$; this fact will make the proof simpler.
  Assume that $x$ is maximal in $t^\star$. $\ll(x)$ will depend only on $h(b)$ where $b$ is the unique branch
  of $t$ leading to $x$. Let us define $\ll(x)$ to be the formula
  of the form
  $\theta\land\Theta_b((x_i)_{i<\a^*})$, where $\theta$ is defined above and $\Theta_b$ is defined below.
  The idea is that
  $$\A_\eta\models \Theta_b((a_{\g})_{\g<\a^*})\}\iff \eta\in h(b)\text{ and }\forall \g<\a^*(a_\g=\g).$$   
  Let us define such a $\Theta_b$.
  Suppose that $\xi$ and $\a$ are such that $h(b)=N_{\xi\restriction\a}$. 
  Define for $s\in \pi^{-1}[\a]$ the formula $A^s_b$ as follows:
  $$A^s_b=
  \begin{cases}
    P_{\dom s}, &\text{ if }\A_\xi\models P_{\dom s}((s(i))_{i\in \dom s})\\
    \lnot P_{\dom s}, &\text{ if }\A_\xi\not \models P_{\dom s}((s(i))_{i\in \dom s})
  \end{cases}
  $$
  Then define
  \begin{eqnarray*}
    \psi_0((x_i)_{i<\a^*})&=&\Land_{i<\a^*}\big[\forall y(R(y,x_i)\leftrightarrow \Lor_{j<i}(y=x_j))\big]\\
    \psi_1((x_i)_{i<\a^*})&=&\Land_{s\in \pi^{-1}[\a]} A^s_b((x_{s(i)})_{i\in \dom s}),\\
    \Theta_b&=&\psi_0\land\psi_1.
  \end{eqnarray*}
  The disjunction over the empty set is considered false.
  \begin{claim}{1}
    Suppose for all $\eta$, $R$ is the standard order relation on $\k$.
    Then
    $$(\A_\eta,R)\models \Theta_b((a_\g)_{\g<\a^*})\iff \eta\in h(b)\land \forall\g<\a^*(\a_\g=\g).$$
  \end{claim}
  \begin{proofVOf}{Claim 1}
    Suppose $\A_\eta\models \Theta((a_\g)_{\g<\a^*})$. Then by $\A_\eta\models \psi_0((a_\g)_{\g<\a^*})$ 
    we have that $(a_\g)_{\g<\a^*}$ is an initial segment of $\dom \A_\eta$ with respect to $R$.
    But $(\dom\A_\eta,R)=(\k,<)$, so $\forall\g<\a^*(\a_\g=\g)$. Assume that $\b<\a$ and $\eta(\b)=1$
    and denote $s=\pi^{-1}(\b)$. Then $\A_\eta\models P_{\dom s}((s(i))_{i\in \dom s})$. Since $\Theta$ is true
    in $\A_\eta$ as well, we must have $A^s_b=P_{\dom s}$ which by definition means that $\A_\xi\models P_{\dom s}((s(i))_{i\in \dom s})$
    and hence $\xi(\b)=\xi(\pi(s))=1$. In the same way one shows that if $\eta(\b)=0$, then $\xi(\b)=0$ for all $\b<\a$. Hence
    $\eta\rest\a=\xi\rest\a$.
    
    Assume then that $a_\g=\g$ for all $\g<\a^*$ and that $\eta\in N_{\xi\restriction \a}$.
    Then $\A_\eta$ trivially satisfies $\psi_0$. Suppose that $s\in\pi^{-1}[\a]$ is such that
    $\A_{\xi}\models P_{\dom s}((s(i))_{i\in \dom s})$.
    Then $\xi(\pi(s))=1$ and since $\pi(s)<\a$, also $\eta(\pi(s))=1$, so $\A_{\eta}\models P_{\dom s}((s(i))_{i\in \dom s})$.
    Similarly one shows that if $$\A_{\xi}\not\models P_{\dom s}((s(i))_{i\in \dom s}),$$ then
    $\A_\eta\not\models P_{\dom s}((s(i))_{i\in \dom s})$. This shows that $\A_\eta\models A^s_b((s(i))_{i\in \dom s})$ for all $s$.
    Hence $\A_\eta$ satisfies $\psi_1$, so we have $\A_\eta\models \Theta$.
  \end{proofVOf}
  \begin{claim}{2}
    $t$, $h$, $t^{\star}$ and $\ll$ are such that for all $\eta\in \k^\k$
    $$\PlTwo\uparrow G(t,h,\eta)\iff\exists R\subset (\dom \A_\eta)^2\ \ \PlTwo\uparrow S(t^\star,\ll,\A_\eta).$$
    \vspace{-10pt}
  \end{claim}
  \begin{proofVOf}{Claim 2}
    Suppose $\s$ is a winning strategy of $\PlTwo$ in $G(t,h,\eta)$. Let $R$ be the well ordering of $\dom \A_\eta$
    such that $(\dom\A_\eta,R)=(\k,<)$.
    Consider the game $S(t^\star,\ll,\A_\eta)$. On the first move
    the players are at the root and player $\PlOne$ chooses where to go next. They go to
    to a minimal element of $t$. From here on $\PlTwo$ uses $\s$ as long as they
    are in $t$. Let us see what happens if they got to a maximal element of $t$, i.e. they picked a branch $b$
    from $t$. Since $\s$ is a winning strategy of $\PlTwo$ in $G(t,h,\eta)$, we have
    $\eta\in h(b)$ and $h(b)=N_{\xi\restriction \a}$ for some $\xi$ and $\a$. For the next $\a$ moves the players
    climb up the tower defined in item (1) of the definition of $t^\star$. All labels are of the form $\exists x_\b$,
    so player $\PlTwo$ has to pick constants from $\A_\eta$. She picks them as follows: 
    for the variable $x_\b$ she picks $\b\in\k=\dom\A_\eta$. She wins now if
    $\A_\eta\models\Theta((\b)_{\b<\a^*})$ and $\A_\eta\models \theta$. But $\eta\in h(b)$, so by Claim~1 the former holds and
    the latter holds because we chose $R$ to be a well ordering of order type $\k$.

    Let us assume that there is no winning strategy of $\PlTwo$ in $G(t,h,\eta)$. 
    Let $R$ be an arbitrary relation on $\dom\A_\eta$. Here we shall finally use the fact that
    $B$ is closed under permutations. Suppose $R$ is not a well ordering of the universe of order type
    $\k$. Then after the players reached the final node of $t^\star$, player $\PlOne$ chooses to go to $\theta$
    and player $\PlTwo$ loses. So we can assume that $R$ is a well ordering of the universe of order type
    $\k$. Let $p\colon \k\to\k$ be a bijection such that $p(\a)$ is the $\a\th$ element of $\k$ with respect to $R$.
    Now $p$ is a permutation and $\{\eta\mid \A_{p\eta}\in B\}=B$ since $B$ is closed under permutations.
    So by our assumption that $\eta\notin B$ (i.e. $\PlTwo\not\uparrow G(t,h,\eta)$), we also have $p\eta\notin B$,
    i.e. player $\PlTwo$ has no winning strategy in $G(t,h,p\eta)$ either.

    Suppose $\s$ is any strategy of $\PlTwo$ in 
    $S(t^\star,\ll,\A_\eta)$. Player $\PlOne$ imagines that $\s$ is a strategy in $G(t,h,p\eta)$ and picks a strategy $\tau$
    that beats it. In the game $S(t^\star,\ll,\A_\eta)$,
    as long as the players are still in $t$, player $\PlOne$ uses 
    $\tau$ that would beat $\s$ if they were playing $G(t,h,p\eta)$ instead of $S(t^\star,\ll,\eta)$.
    Suppose they picked a branch $b$ of $t$. Now
    $p\eta\notin h(b)$. If $\PlTwo$ wants to satisfy $\psi_0$ of the definition of $\Theta_b$, she
    is forced to pick the constants $(a_i)_{i<\a^*}$ such that $a_i$ is the $i\th$ element of $\dom \A_\eta$
    with respect to $R$. Suppose that $\A_\eta\models \psi_1((a_i)_{i<\a^*})$ (recall $\Theta_b=\psi_0\land \psi_1$). But then 
    $\A_{p\eta}\models \psi_1((\g)_{\g<\a^*})$ and also $\A_{p\eta}\models \psi_0((\g)_{\g<\a^*})$, so by
    Claim 1 we should have $p\eta\in h(b)$ which is a contradiction. 
  \end{proofVOf}
\vspace{-20pt}
\end{proofVOf}

\chapter{Generalizations From Classical Descriptive Set Theory}

\section{Simple Generalizations}

\subsection{The Identity Relation}

Denote by $\id$ the equivalence relation $\{(\eta,\xi)\in (2^\k)^2\mid \eta=\xi\}$.
With respect to our choice of topology, the natural generalization of the
equivalence relation
$$E_0=\{(\eta,\xi)\in 2^\o\times 2^\o\mid \exists n<\o\forall m>n(\eta(m)=\xi(m))\}$$
is equivalence modulo sets of size $<\k$:
$$E_0^{<\k}=\{(\eta,\xi)\in 2^\k\times 2^\k\mid \exists \a<\k\forall \b>\a(\eta(\b)=\xi(\b))\},$$
although the equivalences modulo sets of size $<\l$ for $\l<\k$ can also be studied:
$$E_0^{<\l}=\{(\eta,\xi)\in 2^\k\times 2^\k\mid \exists A\subset \k [|A|<\l\land \forall \b\notin A(\eta(\b)=\xi(\b))]\},$$
but for $\l<\k$ these turn out to be bireducible with $\id$ (see below). Similarly one can define $E_0^{<\l}$ on
$\k^\k$ instead of $2^\k$.

It makes no difference whether we define these relations on $2^\k$ or $\k^\k$ since they become bireducible to each other:

\begin{Thm}\label{thm:Bireducible}
  Let $\l\le \k$ be a cardinal and let $E_0^{<\l}(P)$ denote the equivalence relation
  $E_0^{<\l}$ on $P\in \{2^\k,\k^\k\}$ (notation defined above). Then
  $$E^{<\l}_0(2^\k)\le_c E^{<\l}_0(\k^\k)\text{ and }E^{<\l}_0(\k^\k)\le_c E^{<\l}_0(2^\k).$$
  Note that when $\l=1$, we have $E^{<1}_0(P)=\id_{P}$.
\end{Thm}
\begin{proof}
  In this proof we think of functions $\eta,\xi\in \k^\k$ as 
  graphs $\eta=\{(\a,\eta(\a))\mid \a<\k\}$.
  Fix a bijection $h\colon \k\to\k\times\k$. Let $f\colon 2^\k\to\k^\k$ be the inclusion, 
  $f(\eta)(\a)=\eta(\a)$. Then $f$ is easily seen to be a continuous reduction $E^{<\l}_0(2^\k)\le_c E^{<\l}_0(\k^\k)$.
  Define $g\colon \k^\k\to 2^\k$ as follows. For $\eta\in \k^\k$ let
  $g(\eta)(\a)=1$ if $h(\a)\in \eta$ and $g(\eta)(\a)=0$ otherwise.
  Let us show that $g$ is a continuous reduction $E^{<\l}_0(\k^\k)\le_c E^{<\l}_0(2^\k)$.
  Suppose $\eta,\xi\in \k$ are $E^{<\l}_0(\k^\k)$-equivalent. Then clearly $|\eta\sd \xi|<\l$.
  On the other hand
  $$I=\{\a\mid g(\eta)(\a)\ne g(\xi)(\a)\}=\{\a\mid h(\a)\in \eta \sd \xi\}$$
  and because $h$ is a bijection, we have that $|I|<\l$.

  Suppose $\eta$ and $\xi$ are not $E^{<\l}_0(\k^\k)$-equivalent. But then $|\eta\sd \xi|\ge\l$
  and the argument above shows that also $|I|\ge \l$, so $g(\eta)(\a)$ is not $E^{<\l}_0(2^\k)$-equivalent
  to $g(\xi)(\a)$.

  $g$ is easily seen to be continuous.
\end{proof}

We will need the following Lemma which is a straightforward generalization from the case~$\k=\o$:
\begin{Lemma}\label{lem:BaireCont}
  Borel functions are continuous on a co-meager set.
\end{Lemma}
\begin{proof}
  For each $\eta\in\k^{<\k}$ let $V_\eta$ be an open subset of $\k^\k$ such that $V_\eta\sd f^{-1}N_\eta$ is meager.
  Let 
  $$D=\k^\k\setminus \Cup_{\eta\in \k^{<\k}}V_\eta\sd f^{-1}N_\eta.$$
  Then $D$ is as intended. Clearly it is co-meager, since we took away only a $\k$-union of meager sets.
  Let $\xi\in \k^{<\k}$ be arbitrary. The set $D\cap f^{-1}N_\xi$ is open in $D$ since 
   $D\cap f^{-1}N_\xi=D\cap V_\xi$ and so $f\rest D$ is continuous. 
\end{proof}

\begin{Thm}[$\k^{<\k}=\k$]\label{thm:EzeroTheorem}
  $E_0^{<\l}$ is an equivalence relation on $2^\k$ for all $\l\le \k$ and
  \begin{myEnumerate}
  \item $E_0^{<\l}$ is Borel.
  \item $E_0^{<\k}\not\le_B\id$.
  \item If $\l\le \k$, then $\id\le_c E_0^{<\l}$.
  \item If $\l<\k$, then $E_0^{<\l}\le_c\id$.
  \end{myEnumerate}
\end{Thm}
\begin{proof}
  $E_0^{<\l}$ is clearly reflexive and symmetric. Suppose 
  $\eta E_0^{<\l}\xi$ and $\xi E_0^{<\l}\zeta$. Denote $\eta=\eta^{-1}\{1\}$
  and similarly for $\eta,\zeta$.
  Then $|\eta\sd\xi|<\l$ and $|\xi\sd\zeta|<\l$;
  but $\eta\sd\zeta\subset (\eta\sd\xi)\cup(\xi\sd\zeta)$.
  Thus $E_0^{<\l}$ is indeed an equivalence relation. 
  \begin{myEnumerate}
  \item $\displaystyle E_0^{<\l}=\Cup_{A\in [\k]^{<\l}}\Cap_{\a\notin A}\underbrace{\{(\eta,\xi)\mid \eta(\a)=\xi(\a)\}}_{\text{open}}$.
  \item Assume there were a Borel reduction $f\colon 2^\k\to 2^\k$ witnessing $E_0\le_B \id$.
    By Lemma \ref{lem:BaireCont} there are dense open sets $(D_i)_{i<\k}$ such that $f\rest \Cap_{i<\k}D_i$ is continuous.
    If $p,q\in 2^\a$ for some $\a$ and $\xi\in N_p$, let us denote $\xi^{(p/q)}=q\cat (\xi\rest(\k\setminus\a))$, and
    if $A\subset N_p$, denote  
    $$A^{(p/q)}=\{\eta^{(p/q)}\mid \eta\in A\}.$$
    Let $C$ is be the collection of sets, each of which is of the form
    $$\Cup_{q\in 2^\a}[D_i\cap N_p]^{(p/q)}$$
    for some $\a<\k$ and some $p\in 2^\a$. It is easy to see that each such set is dense and open, so $C$
    is a collection of dense open sets.
    By the assumption $\k^{<\k}=\k$, $C$ has size $\k$.
    Also $C$ contains the sets $D_i$ for all $i<\k$, 
    (taking $\a=0$). Denote $D=\Cap_{i<\k}D_i$. Let $\eta\in \Cap C$, $\xi=f(\eta)$
    and $\xi'\ne\xi$, $\xi'\in \ran (f\rest D)$. 
    Now $\xi$ and $\xi'$ have disjoint open neighbourhoods $V$ and $V'$ respectively. 
    Let $\a$ and $p,q\in 2^\a$ be such that $\eta\in N_p$ and such that 
    $D\cap N_p\subset f^{-1}[V]$ and $D\cap N_q\subset f^{-1}[V']$. These $p$ and $q$ exist by the continuity of $f$ on $D$.
    Since $\eta\in \Cap C$ and $\eta\in N_p$, we have 
    $$\eta\in [D_i\cap N_q]^{(q/p)}$$
    for all $i<\k$, which is equivalent to 
    $$\eta^{(p/q)}\in [D_i\cap N_q]$$
    for all $i<\k$, i.e. $\eta^{(p/q)}$ is in $D\cap N_q$. 
    On the other hand (since $D_i\in C$ for all $i<\k$ and because $\eta\in N_p$), 
    we have $\eta\in D\cap N_p$. This implies that $f(\eta)\in V$ and $f(\eta^{(p/q)})\in V'$ which is a contradiction, because $V$
    and $V'$ are disjoint and $(\eta,\eta^{p/q})\in E_0$.
  \item Let $(A_i)_{i<\k}$ be a partition of $\k$ into pieces of size $\k$:
    if $i\ne j$ then $A_i\cap A_j=\es$, $\Cup_{i<\k}A_i=\k$ and $|A_i|=\k$. Obtain
    such a collection for instance by taking a bijection $h\colon\k\to\k\times\k$ and defining
    $A_i=h^{-1}[\k\times \{i\}]$.
    Let $f\colon 2^\k\to 2^\k$ be defined by $f(\eta)(\a)=\eta(i)\iff \a\in A_i$.
    Now if $\eta=\xi$, then clearly $f(\eta)=f(\xi)$ and so $f(\eta) E^{<\l}_0 f(\xi)$. If $\eta\ne\xi$,
    then there exists $i$ such that $\eta(i)\ne \xi(i)$ and we have that
    $$A_i\subset \{\a\mid f(\eta)(\a)\ne f(\xi)(\a)\}$$
    and $A_i$ is of size $\k\ge \l$.
  \item Let $P=\k^{<\k}\setminus \k^{<\l}$.
    Let $f\colon P\to \k$ be a bijection. It induces a bijection $g\colon 2^P\to 2^{\k}$. 
    Let us construct a map $h\colon 2^\k\to 2^{P}$ such that 
    $g\circ h$ is a reduction $E^{<\l}_0\to \id_{2^\k}$.
    Let us denote by $E^{<\l}(\a)$ the equivalence relation on $2^\a$ such that
    two subsets $X,Y$ of $\a$ are $E^{<\l}(\a)$-equivalent if and only if 
    $|X\sd Y|<\l$.
    
    For each $\a$ in $\l<\a<\k$ let
    $h_{\a}$ be any reduction of $E^{<\l}(\a)$ to $\id_{2^\a}$. This exists because both
    equivalence relations have $2^\alpha$ many classes. Now reduce
    $E^{<\l}_0$ to $\id_{\k^{<\k}}$ by 
    $f(A) = (h_\alpha(A \cap \alpha) \mid \lambda \le \alpha < \kappa)$. If $A$, $B$ are $E^{<\l}_0$-equivalent,
    then $f(A) = f(B)$. Otherwise $f_\alpha(A \cap \alpha)$
    differs from $f_\alpha(B \cap \alpha)$ for large enough $\alpha < \kappa$
    because $\lambda$ is less than $\kappa$ and $\k$ is regular. Continuity of $h$ is easy to check.    
    \qedhere
  \end{myEnumerate}
\end{proof}

\section{On the Silver Dichotomy}

To begin with, let us define the Silver Dichotomy and the Perfect Set Property:
\begin{Def}\label{def:SDPSP}
  Let $\C\in \{{\rm Borel},\Dii,{\rm Borel}^*,\Sii,\Pii\}$.

  By \emph{the Silver Dichotomy}, or more specifically, 
  \emph{$\k$-SD for $\C$} we mean the statement that there are no
  equivalence relations $E$ in the class $\C$ such that
  $E\subset 2^{\k}\times 2^{\k}$ and $E$ has more than $\k$ equivalence classes
  such that $\id\not\le_B E$, $\id=\id_{2^\k}$. 

  Similarly the \emph{Perfect Set Property}, or \emph{$\k$-PSP for $\C$}, means that
  each member $A$ of $\C$ has either size $\le\k$ or 
  there is a Borel injection $2^\k\to A$. 
  Using Lemma \ref{lem:BaireCont} it is not hard to see that this definition is equivalent
  to the game definition given in~\cite{MekVaa}.
\end{Def}

\subsection{The Silver Dichotomy for Isomorphism Relations}\label{ssec:SDIR}

Although the Silver Dichotomy for Borel sets is not provable from ZFC for $\k>\o$ (see Theorem \ref{thm:NoSilver} on page \pageref{thm:NoSilver}),
it holds when the equivalence relation is an isomorphism relation, if $\k>\o$ is an inaccessible cardinal:

\begin{Thm}\label{thm:VaughtConjecture}
  Assume that $\k$ is inaccessible.
  If the number of equivalence classes of $\cong_T$ is greater than $\k$,
  then $\id \le_c\,\,\cong_T$.
\end{Thm}
\begin{proof}
  Suppose that there are more than $\k$ equivalence classes of $\cong_T$. We will show that then $\id_{2^\k}\le_c \ \cong_T$.
  If $T$ is not classifiable, then as was done in \cite{Shelah3}, we can construct a tree $t(S)$ for each $S\subset S^\k_\o$ and Ehrenfeucht-Mostowski-type
  models $M(t(S))$ over these trees such that if $S\sd S'$ is stationary, then $M(t(S))\not\cong M(t(S'))$. Now it is easy to
  construct a reduction $f\colon\id_{2^\k}\le_c E_{S^\k_\o}$ (see notation defined in Section \ref{sec:NotationsandC}), 
  so then $\eta\mapsto M(t(f(\eta)))$ is 
  a reduction $\id\le_c\ \cong_T$. 
  
  Assume now that $T$ is classifiable. By $\l(T)$ we denote the least cardinal in which $T$ is stable.
  By \cite{Shelah2} Theorem XIII.4.8 (this is also mentioned in \cite{HarHruLas} Theorem 2.5), assuming that $\cong_T$ has more than $\k$ equivalence
  classes, it has depth at least $2$ and so 
  there are: a $\l(T)^+$-saturated model $\B\models T$, $|\B|=\l(T)$, and a $\l(T)^+$-saturated 
  elementary submodel $\A\esm\B$ and $a\notin \B$ such that
  $\tp(a/\B)$ is orthogonal to $\A$. Let $f\colon\k\to\k$ be strictly increasing and 
  such that for all $\a<\k$, $f(\a)=\mu^+$, for some $\mu$ with the properties
  $\l(T)<\mu<\k$, $\cf(\mu)=\mu$ and $\mu^{2^\o}=\mu$.
  For each $\eta\in 2^\k$ with $\eta^{-1}\{1\}$ is unbounded we will construct a model $\A_\eta$. As above, it will be enough to show that
  $\A_\eta\not\cong\A_\xi$ whenever $\eta^{-1}\{1\}\sd \xi^{-1}\{1\}$ is $\l$-stationary where $\l=\l(T)^+$.
  Fix $\eta\in 2^\k$ and let~$\l=\l(T)^+$.
  
  For each $\a\in \eta^{-1}\{1\}$ choose $\B_\a\supset \A$ such that 
  \begin{enumerate}
  \item\label{item1} $\exists\pi_\a\colon\B\cong\B_\a$, $\pi_{\a}\rest\A=\id_{\A}$.
  \item\label{item2} $\B_\a\downarrow_\A\Cup\{\B_\b\mid \b\in \eta^{-1}\{1\},\b\ne\a\}$
  \end{enumerate}
  Note that \ref{item2} implies that if $\a\ne\b$, then $\B_\a\cap \B_\b=\A$.
  For each $\a\in \eta^{-1}\{1\}$ and $i<f(\a)$ choose tuples $a^{\a}_i$ with the properties
  \begin{enumerate}\setcounter{enumi}{2}
  \item\label{item3} $\tp(a^\a_i/\B_\a)=\pi_\a(\tp(a/\B))$
  \item\label{item4} $a^\a_i\downarrow_{\B_\a}\Cup\{a^\a_j\mid j<f(\a),j\ne i\}$
  \end{enumerate}
  Let $\A_\eta$ be $F_{\l}^s$-primary over 
  $$S_\eta=\Cup\{B_\a\mid a<\eta^{-1}\{1\}\}\cup \Cup\{a^{\a}_i\mid \a<\eta^{-1}\{1\},i<f(\a)\}.$$
  
  It remains to show that if $S^\k_{\l}\cap\eta^{-1}\{1\}\sd \xi^{-1}\{1\}$ is stationary, then $\A_\eta\not\cong\A_\xi$.
  Without loss of generality we may assume that $S^\k_{\l}\cap \eta^{-1}\{1\}\setminus \xi^{-1}\{1\}$ is stationary.
  Let us make a counter assumption, namely that there is an isomorphism $F\colon \A_\eta\to \A_\xi$.

  Without loss of generality there exist singletons $b^\eta_i$ and sets $B^\eta_i$, $i<\k$ of size $<\l$ such that
  $\A_\eta=S_\eta\cup\Cup_{i<\k}b^\eta_i$ and $(S_\eta,(b^\eta_i,B^\eta_i)_{i<\k})$ is an $F^s_{\l}$-construction.

  Let us find an ordinal $\a<\k$ and sets $C\subset \A_\eta$ and $D\subset \A_\xi$ with the properties listed below:
  \begin{myAlphanumerate}
  \item \label{ai1} $\a\in \eta^{-1}\{1\}\setminus \xi^{-1}\{1\}$
  \item \label{aai2} $D=F[C]$
  \item \label{ai3} $\forall \b\in(\a+1)\cap \eta^{-1}\{1\}(\B_\b\subset C)$ and $\forall \b\in(\a+1)\cap \xi^{-1}\{1\}(\B_\b\subset D)$, 
  \item \label{ai4} for all $i<f(\a)$, $\forall \b\in \a\cap \eta^{-1}\{1\} (a^{\b}_i\in C)$ and $\forall \b\in \a\cap \xi^{-1}\{1\} (a^{\b}_i\in D)$,
  \item \label{ai5} $|C|=|D|<f(\a)$,
  \item \label{ai6} For all $\b$, if $\B_\b\cap C\setminus \A\ne \es$, then $\B_\b\subset C$ and if 
    $\B_\b\cap D\setminus \A\ne \es$, then $\B_\b\subset D$,
  \item \label{ai7} $C$ and $D$ are $\l$-saturated,
  \item \label{ai8} if $b^\eta_i\in C$, then 
    $B^\eta_i\subset [S_\eta\cup \Cup\{b^\eta_i\mid j< i\}]\cap C$ and 
    if $b^{\xi}_i\in D$, then 
    $B^\xi_i\subset [S_\xi\cup \Cup\{b^\xi_i\mid j< i\}]\cap D$.
  \end{myAlphanumerate}
  This is possible, because $\eta^{-1}\{1\}\setminus \xi^{-1}\{1\}$ is stationary and we can close under the
  properties~\mbox{(b)--(h)}.

  Now $\A_\eta$ is $F^s_{\l}$-primary over $C\cup S_\eta$ and $\A_\xi$ is $F^s_{\l}$-primary over $D\cup S_\eta$ and
  thus $\A_\eta$ is $F_{\l}^s$-atomic over $C\cup S_\eta$ and $\A_\xi$ is $F_{\l}^s$-atomic over $D\cup S_\xi$. Let 
  $$I_{\a}=\{a^\a_i\mid i<f(\a)\}.$$
  Now $|I_\a\setminus C|=f(\a)$, because $|C|<f(\a)$, and so $I_\a\setminus C\ne\es$. Let $c\in I_\a\setminus C$
  and let $A\subset S_\xi\setminus D$ and $B\subset D$ be such that $\tp(F(c)/A\cup B)\vdash \tp(F(c)/D\cup S_\xi)$
  and $|A\cup B|<\l$.
  Since $\a\notin \xi^{-1}\{1\}$, we can find (just take disjoint copies) a sequence $(A_i)_{i<f(\a)^+}$ such that $A_i\subset I_\a\cap\A_\xi$, 
  $\tp(A_i/D)=\tp(A/D)$ and $A_i\downarrow_D \Cup\{A_j\mid j\ne i, j<f(\a)^+\}$
  
  Now we can find $(d_i)_{i<f(\a)^+}$, such that 
  $$\tp(d_i\cat A_i\cat B_i/\es)=\tp(F(c)\cat A\cat B/\es).$$
  Then it is a Morley sequence over $D$ and for all $i<f(\a)^+$, 
  $$\tp(d_i/D)=\tp(F(c)/D),$$
  which implies 
  $$\tp(F^{-1}(d_i)/C)=\tp(c/C),$$
  for some $i$, since  for some $i$ we have $c=a^\a_i$.
  Since by (c), $\B_\a\subset C$, the above implies that
  $$\tp(F^{-1}(d_i)/\B_\a)=\tp(a^\a_i/\B_\a)$$
  which by the definition of $a^\a_i$, item \ref{item3} implies
  $$\tp(F^{-1}(d_i)/\B_\a)=\pi_{\a}(\tp(a/\B)).$$
  Thus the sequence $(F^{-1}(d_i))_{i<f(\a)^+}$
  witnesses that the dimension of $\pi_{\a}(\tp(a/\B))$ in $\A_\eta$ is greater than 
  $f(\a)$. Denote that sequence by $J$. Since $\pi_{\a}(\tp(a/\B))$ is orthogonal to $\A$, 
  we can find $J'\subset J$ such that
  $|J'|=f(\a)^+$ and $J'$ is a Morley sequence over $S_\eta$. Since $f(\a)^+>\l$, this
  contradicts Theorem 4.9(2) of Chapter IV of~\cite{Shelah2}. 
\end{proof}

\begin{Open}
  Under what conditions on $\kappa$ does the conclusion of Theorem \ref{thm:VaughtConjecture} hold?
\end{Open}

\subsection{Theories Bireducible With $\id$}

\begin{Thm}
Assume $\k^{<\k}=\k=\aleph_{\a}>\o$, $\k$ is not weakly inaccessible
and $\l =\vert\a +\o \vert$. Then the following are equivalent.
\begin{myEnumerate}
\item There is $\g <\o_{1}$ such that $\beth_{\g}(\l )\ge\k$.
\item There is a complete countable $T$ such that $\id\le_{B}\cong_{T}$ and $\cong_{T}\le_{B}\id$.  
\end{myEnumerate}
\end{Thm}
\begin{proof}
(2)$\Rightarrow$(1): Suppose that (1) is not true.
Notice that then $\k >2^{\o}$.
Then every shallow classifiable theory has $<\k$ many models of power $\k$
(see \cite{HarHruLas}, item 6. of the Theorem which is on the first page of the article.)
and thus $\id\not\le_{B}\cong_{T}$. On the other hand
if $T$ is not classifiable and shallow,
$\cong_{T}$ is not Borel by Theorem \ref{thm:ClasNotShalIsNotBorel} and thus it is not Borel reducible
to~$\id$ by Fact~\ref{fact:ReductionNote}.

(1)$\Rightarrow$(2): Since $\cf(\k )>\o$, (1) implies that there is
$\a =\b +1<\o_{1}$ such that
$\beth_{\a}(\l )=\k$. But then there is an $L^{*}$-theory
$T^{*}$ which has exactly $\k$ many models in cardinality $\k$
(up to isomorphism, use \cite{HarHruLas}, Theorem 6.1 items 2. and 8.).
But then it has exactly $\k$ many models of cardinality $\le\k$,
let $\A_{i}$, $i<\k$, list these. Such a theory must be classifiable and
shallow.
Let $L$ be the vocabulary we get from $L^{*}$ by adding one binary
relation symbol $E$. Let
$\A$ be an $L$-structure in which $E$ is an equivalence relation
with infinitely many equivalence classes such that
for every equivalence class $a/E$, $(\A\rest a/E)\rest L^{*}$
is a model of $T^{*}$. Let~$T=Th(\A )$.

We show first that identity on $\{\eta\in 2^{\k}
\vert\ \eta (0)=1\}$ reduces to $\cong_{T}$.
For all $\eta\in 2^{\k}$, let $\B_{\eta}$ be a model of $T$
of power $\k$ such that if $\eta (i)=0$, then the number of equivalence classes
isomorphic to $\B_{i}$ is countable and otherwise the number
is $\k$. Clearly we can code $\B_{\eta}$ as $\xi_{\eta}\in 2^{\k}$
so that $\eta\mapsto\xi_{\eta}$ is the required Borel reduction.

We show then that $\cong_{T}$ Borel reduces to identity on
$$X=\{ \eta \colon\k\rightarrow (\k +1)\}.$$
Since $T^{*}$ is classifiable and shallow, for all $\d ,i<\k$ the set 
$$\{\eta\in X\vert\ (\A_{\eta}\rest \d /E)\rest L^{*}\cong\A_{i}\}$$
is Borel. But then for all cardinals $\theta\le\k$ and $i<\k$,
the set
$$\{\eta\in X\mid \card(\{ \d/E\mid\ \d <\k,\ (\A_{\eta}\rest \d /E)\rest L^{*}\cong\A_{i}\}) =\theta\}$$
is Borel. But then
$\eta\mapsto\xi_{\eta}$ is the required reduction when\\
  \begin{equation*}
    \xi_{\eta}(i)=\vert\{ \d/E\mid\ \d <\k,\ (\A_{\eta}\rest \d /E)\rest L^{*}\cong\A_{i}\}\vert.      \qedhere
  \end{equation*}
\end{proof}

\subsection{Failures of Silver's Dichotomy}\label{sec:Silver}

There are well-known dichotomy theorems for Borel equivalence relations on $2^\o$. Two of them are:

\begin{Thm}[Silver, \cite{Silver}]\label{thm:ClassSilver}
  Let $E\subset 2^\o\times 2^\o$ be a $\Pii$ equivalence relation. If $E$ has uncountably many equivalence classes, 
  then $\id_{2^\o}\le_B E$.\qed
\end{Thm}

\begin{Thm}[Generalized Glimm-Effros dichotomy, \cite{HaKeLo}]
  Let $E\subset 2^\o\times 2^\o$ be a Borel equivalence relation. Then either
  $E\le_B \id_{2^\o}$ or else $E_0 \le_c E$.\qed
\end{Thm}

As in the case $\k=\o$ we have the following also for uncountable $\k$ (see Definition~\ref{def:SDPSP}):

\begin{Thm}
  If $\k$-SD for $\Pii$ holds, then the $\k$-PSP holds for $\Sii$-sets.
  More generally, if $\C\in \{{\rm Borel},\Dii,{\rm Borel}^*,\Sii,\Pii\}$, then
  $\k$-SD for $\C$ implies $\k$-PSP for $\C'$, where 
  elements in $\C'$ are all the complements of those in $\C$.
\end{Thm}
\begin{proof}
  Let us prove this for $\C=\Pii$, the other cases are similar.
  Suppose we have a $\Sii$-set $A$. Let 
  $$E=\{(\eta,\xi)\mid \eta=\xi\text{ or } ((\eta\notin A)\land(\xi\notin A))\}.$$
  Now $E=\id \cup (2^\k\setminus A)^2$.
  Since $A$ is $\Sii$, $(2^\k\setminus A)^2$ is $\Pii$
  and because $\id$ is Borel, also $E$ is $\Pii$. Obviously $|A|$ is the 
  number of equivalence classes of $E$ provided $A$ is infinite. 
  Then suppose $|A|>\k$. Then there are more than $\k$  equivalence classes of $E$, so 
  by $\k$-SD for $\Pii$, there is a reduction $f\colon \id\le E$. 
  This reduction in fact witnesses the PSP of $A$.
\end{proof}

The idea of using Kurepa trees for this purpose arose already in the paper \cite{MekVaa} by Mekler and V\"a\"an\"anen.

\begin{Def}
  If $t\subset 2^{<\k}$ is a tree, a \emph{path} through $t$ is a branch of length $\k$.
  A \emph{$\k$-Kurepa tree} is a tree $K\subset 2^{<\k}$ which satisfies the following:
  \begin{myAlphanumerate}
  \item $K$ has more than $\k$ paths,
  \item $K$ is downward closed,
  \item for all $\a<\k$, the levels are small:  $|\{p\in K\mid \dom p=\a\}|\le |\a+\o|.$    
  \end{myAlphanumerate}
\end{Def}

\begin{Thm}\label{thm:NoSilver}
  Assume one of the following: 
  \begin{myEnumerate}
  \item $\k$ is regular but not strongly inaccessible and there exists a $\k$-Kurepa tree $K\subset 2^{<\k}$,
  \item $\k$ is regular (might be strongly inaccessible), $2^\k>\k^+$ and there exists a tree $K\subset 2^{<\k}$ 
    with more than $\k$ but less than $2^\k$ branches.
  \end{myEnumerate}
  Then the Silver Dichotomy for $\k$ does not hold. In fact there an equivalence 
  relation
  $E\subset 2^{\k}\times 2^{\k}$ which is the union of a closed and an open set, 
  has more than $\k$ equivalence classes but $\id_{2^{\k}}\not\le_B E$.
\end{Thm}
\begin{proof}
  Let us break the proof according to the assumptions (1) and (2). So first let us consider the case where $\k$ is
  not strongly inaccessible and there is a $\k$-Kurepa tree.
    \vspace{5pt}
    
  \noindent {\bf(1)}:
    Let us carry out the proof in the case $\k=\o_1$. It should be obvious then how to generalize it to any
    $\k$ not strongly inaccessible. So let $K\subset 2^{<\o_1}$ be an $\o_1$-Kurepa tree.
    Let $P$ be the collection of all paths of $K$. For $b\in P$, denote $b=\{b_\a\mid \a<\o_1\}$ where
    $b_\a$ is an element of $K$ with domain $\a$.
    
    Let 
    $$C=\{\eta\in 2^{\o_1}\mid \eta=\Cup_{\a<\o_1}b_\a, b\in P\}.$$
    Clearly $C$ is closed.
    
    Let $E=\{(\eta,\xi)\mid (\eta\notin C\land \xi\notin C)\lor (\eta\in C\land \eta=\xi)\}$. In words,
    $E$ is the equivalence relation whose equivalence classes are the complement of $C$ and 
    the singletons formed by the elements of $C$. $E$ is the union of the open set
    $\{(\eta,\xi)\mid \eta\notin C\land \xi\notin C\}$ and the closed set
    $\{(\eta,\xi)\mid \eta\in C\land \eta=\xi\}=\{(\eta,\eta)\mid \eta\in C\}$.
    The number of equivalence classes equals the number of paths of $K$, so
    there are more than $\o_1$ of them by the definition of Kurepa tree. 
    
    Let us show that $\id_{2^{\o_1}}$ is not embeddable
    to $E$. Suppose that $f\colon 2^{\o_1}\to 2^{\o_1}$ is a Borel reduction. We will show that then $K$ must have
    a level of size $\ge \o_1$ which contradicts the definition of Kurepa tree. By Lemma \ref{lem:BaireCont}
    there is a co-meager set $D$ on which $f\rest D$ is continuous. There is at most one $\eta\in 2^{\o_1}$
    whose image $f(\eta)$ is outside $C$, so without loss of generality $f[D]\subset C$. Let $p$ be an arbitrary
    element of $K$ such that $f^{-1}[N_p]\ne\es$. By continuity there is a $q\in 2^{<\o_1}$ 
    with $f[N_{q}\cap D]\subset N_p$. 
    Since 
    $D$ is co-meager, there are $\eta$ and $\xi$ such that $\eta\ne\xi$, $q\subset \eta$ and $q\subset \xi$.
    Let $\a_1<\o_1$ and $p_{0}$ and $p_{1}$ be extensions of $p$ with the properties 
    $p_{0}\subset f(\eta)$, $p_{1}\subset f(\xi)$, 
    $\a_1=\dom p_{0}=\dom p_{1}$, $f^{-1}[N_{p_0}]\ne\es\ne f^{-1}[N_{p_1}]$
    and $N_{p_{0}}\cap N_{p_{1}}=\es$. Note that $p_{0}$ and $p_{1}$ are in $K$. Then, again by continuity, 
    there are $q_{0}$ and $q_{1}$ such that
    $f[N_{q_{0}}\cap D]\subset N_{p_{0}}$ and $f[N_{q_{1}}\cap D]\subset N_{p_{1}}$. Continue in the same manner to
    obtain $\a_n$ and $p_s\in K$ for each $n<\o$ and $s\in 2^{<\o}$ so that $s\subset s'\iff p_{s}\subset p_{s'}$ and 
    $\a_n=\dom p_s\iff n=\dom s$. 
    Let $\a=\sup_{n<\o}\a_n$. Now clearly the $\a$'s level of $K$
    contains continuum many elements: by (b) in the definition of Kurepa tree
    it contains all the elements of the form $\Cup_{n<\o}p_{\eta\restl n}$
    for $\eta\in 2^\o$ and $2^\o\ge \o_1$.
    
    If $\k$ is arbitrary regular not strongly inaccessible cardinal, then the proof is the same, only instead of 
    $\o$ steps one has to do $\l$ steps where $\l$ is the least cardinal satisfying $2^\l\ge \k$.
    \vspace{5pt}
    
    \noindent {\bf(2)}: The argument is even simpler. Define the equivalence relation $E$ exactly as above. Now 
    $E$ is again closed and has as many equivalence classes as is the number of paths in $K$. Thus the number of 
    equivalence classes is $>\k$ but $\id$ cannot be reduced to $E$ since there are less than $2^\k$ equivalence classes.
\end{proof}

\begin{Remark}
  Some related results:
  \begin{myEnumerate}
  \item In $L$, the PSP fails for closed sets for all uncountable regular $\k$. 
    This is because ``weak Kurepa trees'' exist (see the proof sketch of (3) below for the definition of ``weak Kurepa tree'').
  \item (P. Schlicht)
    In Silver's model where an inaccessible 
    $\kappa$ is made into $\omega_2$ by Levy collapsing each ordinal below to $\omega_1$ with 
    countable conditions, every $\Sigma^1_1$ subset $X$ of $2^{\omega_1}$ obeys the PSP.
  \item Supercompactness does not imply the PSP for closed sets.
  \end{myEnumerate}
  \begin{proof}[Sketch of a proof of item (3)]
    Suppose $\k$ is supercompact and by a reverse Easton iteration add
    to each inaccessible $\alpha$ a ``weak Kurepa tree'', i.e., a
    tree $T_\alpha$ with $\alpha^+$ branches whose $\beta\th$ level has size $\beta$
    for stationary many $\beta < \alpha$. The forcing at stage $\a$ is $\a$-closed and
    the set of branches through $T_\kappa$ is a closed set with no perfect
    subset. If $j\colon V \to M$ witnesses $\l$-supercompactness ($\lambda >\kappa$) and $G$ is the generic then we can find $G^*$ which is
    $j(P)$-generic over $M$ containing $j[G]$: Up to $\l$ we copy $G$, between
    $\lambda$ and $j(\kappa)$ we build $G^*$ using $\lambda^+$ closure of the
    forcing and of the model $M$, and at $j(\kappa)$ we form a master
    condition out of $j[G(\kappa)]$ and build a generic below it, again
    using $\lambda^+$ closure.    
  \end{proof}
\end{Remark}

\begin{Cor}
  The consistency of the Silver Dichotomy for Borel sets on $\o_1$ with CH implies the 
  consistency of a strongly inaccessible cardinal.
  In fact, if there is no equivalence relation witnessing the failure of the Silver Dichotomy for $\o_1$, 
  then $\o_2$ is inaccessible in $L$.
\end{Cor}
\begin{proof}
  By a result of Silver, if there are no $\o_1$-Kurepa trees, then $\o_2$ is 
  inaccessible in $L$, see Exercise 27.5 in Part III of \cite{Jech}.
\end{proof}


\begin{Open}
  Is the Silver Dichotomy for uncountable $\kappa$ consistent?
\end{Open}

\section{Regularity Properties and Definability of the CUB Filter}\label{sec:PB}

In the standard descriptive theory ($\k=\o$), the notions of Borel, $\Dii$ and Borel* coincide 
and one of the most important observations in the theory is that such sets have the Property of Baire and that the
$\Sii$-sets obey the Perfect Set Property.
In the case $\k>\o$ the situation is more complicated as the following shows. 
It was already pointed out in the previous section that Borel $\subsetneq \Dii$. 
In this section we focus on the cub filter 
$$\CUB=\{\eta\in 2^{\k}\mid \eta^{-1}\{1\}\text{ contains a cub}\}.$$ 
The set $\CUB$ is easily seen to be $\Sii$: the set
$$\{(\eta,\xi)\mid (\eta^{-1}\{1\}\subset \xi^{-1}\{1\})\land (\eta^{-1}\{1\}\text{ is cub})\}$$
is Borel.
$\CUB$ (restricted to cofinality $\o$, see Definition \ref{def:CUBRestSet}) 
will serve (consistently) as a counterexample to $\Dii=$ Borel*, but we will show that it is also consistent that
CUB is $\Dii$. The latter implies that it is consistent that $\Dii$-sets do not have the Property of Baire
and we will also show that in a forcing extension of $L$, $\Dii$-sets all have the Property of Baire.

\begin{Def}\label{def:BaireProperty}
  A \emph{nowhere dense set} is a subset of a set whose complement is dense and open.
  Let $X\subset \k^\k$. A subset $M\subset X$ is $\k$-\emph{meager in $X$}, if $M\cap X$ is the 
  union of no more than $\k$ nowhere dense sets,
  $$M=\Cup_{i<\k}N_i.$$
  We usually drop the prefix ``$\k$-''.

  Clearly $\k$-meager sets form a $\k$-complete ideal. A \emph{co-meager} set is a set whose complement
  is meager.
  
  A subset $A\subset X$ \emph{has the Property of Baire} or shorter \emph{P.B.}, if there exists an open $U\subset X$
  such that the symmetric difference $U\sd A$ is meager.
\end{Def}

Halko showed in \cite{Halko} that 
\begin{Thm}[\cite{Halko}]\label{thm:BorelNoBaire}
  Borel sets have the Property of Baire.\qed  
\end{Thm}
(The same proof as when $\k=\o$ works.)
This is independent of the assumption $\k^{<\k}=\k$.
Borel* sets do not in general have the Property of Baire.

\begin{Def}[\cite{MekShe,MekVaa,HytRau}]\label{def:Canary}
  A $\k^+\k$-tree $t$ is a $\k\l$-\emph{canary tree} if for all stationary $S\subset S^{\k}_\l$ it holds that
  if $\P$ does not add subsets of $\k$ of size less than $\k$ and $\P$ kills the stationarity of $S$, then $\P$
  adds a $\k$-branch to $t$.
\end{Def}
  
\begin{Remark}
  Hyttinen and Rautila \cite{HytRau} use the notation \emph{$\k$-canary tree} for our \emph{$\k^+\k$-canary tree}.
\end{Remark}
  
It was shown by Mekler and Shelah \cite{MekShe} and Hyttinen and Rautila \cite{HytRau} that it is consistent with ZFC+GCH that there is a
$\k^+\k$-canary tree \emph{and} it is consistent with ZFC+GCH that there are no $\k^+\k$-canary trees. 
The same proof as in \cite{MekShe, HytRau} gives the following:

\begin{Thm}\label{thm:NoCanary}
  Assume GCH and assume $\l<\k$ are regular cardinals.
  Let $\P$ be the forcing which adds $\k^+$ Cohen subsets of $\k$. Then in the forcing extension
  there are no $\k\l$-canary trees. \qed
\end{Thm}  

\begin{Def}\label{def:CUBRestSet}
  Suppose $X\subset \k$ is stationary. For each such $X$ define the set 
  $$\CUB(X)=\{\eta\in 2^\k\mid X\setminus \eta^{-1}\{1\}\text{ is non-stationary}\},$$
  so $\CUB(X)$ is ``cub in $X$''.   
\end{Def}

\begin{Thm}\label{thm:CUBSET} In the following $\k$ satisfies $\k^{<\k}=\k>\o$.
  \begin{myEnumerate}
  \item $\operatorname{CUB}(S^\k_\o)$ is Borel*.\label{thm:CUBSETItem1}
  \item For all regular $\l<\k$, $\operatorname{CUB}(S^\k_\l)$ is not $\Dii$ in the forcing extension after
    adding $\k^+$ Cohen subsets of $\k$.\label{thm:CUBSETItem2}
  \item If $V=L$, then for every stationary $S\subset \k$, the set $\CUB(S)$ is not $\Dii$.\label{thm:CUBSETItem7}
  \item Assume GCH and that $\k$ is not a successor of a singular cardinal. For any stationary set $Z\subset\k$ there exists a forcing 
    notion $\P$ which has the $\k^+$-c.c., does not add bounded subsets of $\k$ and preserves GCH and 
    stationary subsets of $\k\setminus Z$
    such that $\operatorname{CUB}(\k\setminus Z)$ is $\Dii$ in the forcing extension. \label{thm:CUBSETItem4}
  \item Let the assumptions for $\k$ be as in~\eqref{thm:CUBSETItem4}. 
    For all regular $\l<\k$, 
    $\operatorname{CUB}(S^\k_\l)$ is $\Dii$ in a forcing extension 
    as in~\eqref{thm:CUBSETItem4}. \label{thm:CUBSETItem3}
  \item $\operatorname{CUB}(X)$ does not have the Property of Baire for stationary $X\subset \k$. 
    (Proved by Halko and Shelah in \cite{HalShe} for $X=\k$)\label{thm:CUBSETItem5}
  \item It is consistent that all $\Dii$-sets have the Property of Baire. (Independently known to P. L\"ucke and P. Schlicht.)\label{thm:CUBSETItem6}
  \end{myEnumerate}
\end{Thm}
\vspace{-15pt}
\begin{proofVOf}{Theorem \ref{thm:CUBSET}}
  \begin{proofVOf}{item \eqref{thm:CUBSETItem1}}
    Let $t=[\k]^{<\o}$ (increasing functions ordered by end extension) and for all branches $b\subset t$
    $$h(b)=\{\xi\in 2^{\k}\mid \xi(\sup_{n<\o}b(n))\ne 0\}.$$
    Now if $\k\setminus\xi^{-1}\{0\}$ contains an $\o$-cub set $C$, then player $\PlTwo$ has a winning strategy in 
    $G(t,h,\xi)$: for her $n\th$ move she picks an element $x\in t$ with domain $2n+2$ such that
    $x(2n+1)$ is in $C$. Suppose the players picked a branch $b$ in this way. 
    Then the condition $\xi(b(2n+1))\ne 0$ holds for all $n<\o$ and because $C$ is cub  
    outside $\xi^{-1}\{0\}$, we have $\xi(\sup_{n<\o}b(n))\ne 0$. 
    
    Suppose on the contrary that $S=\xi^{-1}\{0\}$ is stationary. 
    Let $\s$ be any strategy of player $\PlTwo$. Let $C_{\s}$ be the set of ordinals closed under this strategy.
    It is a cub set, so there is an $\a\in C_{\s}\cap S$. Player $\PlOne$ can now easily play
    towards this ordinal to force $\xi(b(\o))=0$, so $\s$ cannot be a winning strategy.
  \end{proofVOf}
  \begin{proofVOf}{item \eqref{thm:CUBSETItem2}}
    It is not hard to see that $\operatorname{CUB}^\k_\l$ is $\Dii$ if and only if there exists a 
    $\k\l$-canary tree. This fact is proved in detail 
    in \cite{MekVaa} in the case $\k=\o_1$, $\l=\o$ and the proof generalizes easily to any regular uncountable $\k$
    along with the assumption $\k^{<\k}=\k$. So the statement follows from Theorem~\ref{thm:NoCanary}.
  \end{proofVOf}
  \begin{proofVOf}{item \eqref{thm:CUBSETItem7}}
    Suppose that $\f$ is $\Sigma_1$ and for simplicity assume that $\f$ has no parameters. 
    Then for $x \subset \k$ we have:

    \begin{claim*}
      $\f(x)$ holds if and only if  
      the set $A$ of those $\a$ for which there exists $\b>\a$ such that
      $$L_\beta \models \big(ZF^-\land (\omega < \a \text{ is regular}) \land ((S \cap \a)\text{ is stationary })\land \f(x \cap \a)\big)$$
      contains $C \cap S$ for some cub set~$C$.
    \end{claim*}
    \begin{proofVOf2}{the Claim}{Claim}
      ``$\Rightarrow$''. If $\f(x)$ holds then choose a continuous chain $(M_i \mid i < \k)$ of
        elementary submodels of some large $ZF^-$ model $L_\theta$ so that $x$ and $S$
        belong to $M_0$ and the intersection of each $M_i$ with $\k$ is an ordinal
        $\a_i$ less than $\k$. Let $C$ be the set of $\a_i$'s, cub in $\k$. Then
        any $\a$ in $C \cap S$ belongs to $A$ by condensation.

      ``$\Leftarrow$''. If $\f(x)$ fails then let $C$ be any cub in $\k$ and let $D$ be the
        cub of $\a < \k$ such that $H(\a)$ is the Skolem Hull in some large
        $L_\theta$ of $\a$ together with $\{\k, S, C\}$ contains no ordinals in the
        interval $[\a,\k)$. Let $\a$ be the least element of $S \cap \lim(D)$.  Then
        $\a$ does not belong to $A$: If $L_\beta$ satisfies $\f(x \cap \a)$ then $\beta$
        must be greater than $\bar\b$ where $\overline{H(\a)} = L_{\bar\b}$ is the
        transitive collapse of $H(\a)$, because $\f(x \cap \a)$ fails in
        $\overline{H(\a)}$. But as $\lim(D) \cap \a$ is an element of $L_{\bar\b + 2}$ and
        is disjoint from $S$, it follows that either $\a$ is singular in $L_\b$ or
        $S \cap \a$ is not stationary in $L_{\bar\b + 2}$ and hence not in
        $L_\b$. Of course $\a$ does belong to $C$ so we have shown that $A$ does not
        contain $S \cap C$ for an arbitrary cub $C$ in $\k$.
    \end{proofVOf2}

    It follows from the above that any $\Sigma_1$ subset of $2^\k$ is
    $\Delta_1$ over $(L_\k^+, \CUB(S))$ and therefore if $\CUB(S)$ were $\Delta_1$ then any $\Sigma_1$ subset of $2^\k$ would be
    $\Delta_1$, a contradiction.
  \end{proofVOf}
  \begin{proofVOf}{item \eqref{thm:CUBSETItem4}}
    If $X\subset 2^\k$ is $\Dii$, then $\{\eta\in X\mid \eta^{-1}\{1\}\subset \k\setminus Z\}$ is $\Dii$, so it is sufficient
    to show that we can force a set $E\subset Z$ which has the claimed property.
    So we force a set $E\subset Z$ such that $E$ is stationary but $E\cap \a$ is non-stationary in $\a$ for all
    $\a<\k$ and $\k\setminus E$ is fat. A set is \emph{fat} if its intersection with any cub set contains closed increasing sequences
    of all order types $<\k$. 
    
    This can be easily forced with 
    $$\R= \{p\colon \a\to 2\mid \a<\k, p^{-1}\{1\}\cap\b\subset Z\text{ is non-stationary in }\b\text{ for all }\b\le \a\}$$
    ordered by end-extension. It is easy to see that for any $\R$-generic $G$ the set $E=(\cup G)^{-1}\{1\}$ satisfies
    the requirements. Also $\R$ does not 
    add bounded subsets of $\k$ and has the $\k^+$-c.c. and does not kill stationary sets.
    
    Without loss of generality assume that such $E$ exists in $V$ and that $0\in E$.

    Next let $\P_0=\{p\colon \a\to 2^{<\a}\mid \a<\k, p(\b)\in 2^{\b}, p(\b)^{-1}\{1\}\subset E\}$. This forcing adds a 
    $\diamondsuit_E$-sequence 
    $\la A_\a\mid \a\in E\ra$ (if $G$ is generic, set $A_\a=(\cup G)(\a)^{-1}\{1\}$) such that for all $B\subset E$ there is a stationary
    $S\subset E$ such that $A_\a=B\cap \a$ for all $\a\in S$. 
    This forcing $\P_0$ is $<\k$-closed and clearly has the $\k^+$-c.c., so it is easily seen that
    it does not add bounded subsets of $\k$ and does not kill stationary sets. 

    Let $\psi(G,\eta,S)$ be a formula with parameters $G\in (2^{<\k})^{\k}$ and $\eta\in 2^\k$ and a free variable $S\subset \k$ which says:
    $$\forall \a<\k(\a\in S\iff G(\a)^{-1}\{1\}=\eta^{-1}\{1\}\cap\a).$$
    If $\la G(\a)^{-1}\{1\}\ra_{\a<\k}$ happens to be a $\diamondsuit_E$-sequence, then $S$ satisfying $\psi$ is always stationary.
    Thus if $G_0$ is $\P_0$-generic over $V$ and $\eta\in 2^{E}$, then $(\psi(G_0,\eta,S)\rightarrow (S\text{ is stationary}))^{V[G_0]}$. 

    For each $\eta\in 2^{E}$, let $\dot S_\eta$ be a nice $\P_0$-name for the set $S$ such that 
    $V[G_0]\models \psi(G_0,\eta,S)$
    where $G_0$ is $\P_0$-generic over $V$.
    By the definitions, $\P_0\forces$ ``$\dot S_\eta\subset \check E$ is stationary''
    and if $\eta\ne\eta'$, then $\P_0\forces$ ``$\dot S_\eta\cap \dot S_{\eta'}$ is bounded''.

    Let us enumerate $E=\{\b_i\mid i<\k\}$ such that $i<j\Rightarrow \b_i<\b_j$ and 
    for $\eta\in 2^{E}$ and $\g\in \k$ define $\eta+\g$ to be the $\xi\in 2^E$ such that
    $\xi(\b_i)=1$ for all $i<\g$ and $\xi(\b_{\g+j})=\eta(\b_j)$ for $j\ge 0$. 
    Let $$F_0=\{\eta\in 2^E \mid \eta(0)=0\}^V\eqno(*)$$
    Now for all $\eta, \eta'\in F_0$ and $\a,\a'\in \k$, 
    $\eta + \a=\eta'+\a'$ implies $\eta=\eta'$ and $\a=\a'$. Let us now define the formula $\f(G,\eta,X)$ with parameters
    $G\in (2^{<\k})^\k$, $\eta\in 2^{\k}$ and a free variable $X\subset \k\setminus E$ which says:
    \begin{eqnarray*}
      (\eta(0)=0)\land\forall \a<\k[ \!\!\!\!\!\!\! && \!\!\!\!\!\!\!\!\!\!\!\! (\a\in    X\rightarrow \exists S(\psi(G,\eta+2\a,  S)\land S\text{ is non-stationary}))\\
                                    &\land&\!\!\!\!\!\!\!\!\!\! (\a\notin X\rightarrow \exists S(\psi(G,\eta\!+\!2\a\!+\!1,S)\land S\text{ is non-stationary}))].      
    \end{eqnarray*}
    
    Now, we will construct an iterated forcing $\P_{\k^+}$, starting with $\P_0$, which
    kills the stationarity of $\dot S_\eta$ for suitable $\eta\in 2^E$, such that if $G$ is $\P_{\k^+}$-generic, 
    then for all 
    $S\subset \k\setminus E$, $S$ is stationary if and only if 
    $$\exists \eta\in 2^{E}(\f(G_0,\eta,S))$$
    where $G_0=G\rest \{0\}$.
    In this model, for each $\eta\in F_0$, there will be a unique $X$ such that $\f(G_0,\eta,X)$, so let us denote this $X$
    by $X_\eta$.
    It is easy to check that the mapping $\eta\mapsto X_\eta$ defined by $\f$ is $\Sii$ 
    so in the result, also $\S=\{S\subset \k\setminus E\mid S\text{ is stationary}\}$
    is~$\Sii$. Since cub and non-stationarity are also $\Sii$, we get that $\S$ is $\Dii$, as needed.

    Let us show how to construct the iterated forcing. 
    For $S\subset \k$, we denote by $T(S)$ the partial order of all closed increasing sequences 
    contained in the complement of $S$. Clearly $T(S)$ is a forcing that kills the stationarity of $S$.
    If the complement of $S$ is fat and $S$ is non-reflecting, then $T(S)$ has all the nice properties we need, as the 
    following claims show.
    Let $f\colon \k^+\setminus \{0\}\to \k^+\times \k^+$ be a bijection such that $f_1(\g)\le\g$. 

    $\P_0$ is already defined and it has the $\k^+$-c.c. and it is $<\k$-closed.
    Suppose that $\P_i$ has been defined for $i<\a$ and $\s_i$ has been defined for
    $i<\cup\a$ such that $\s_i$ is a (nice) $\P_i$-name for a $\k^+$-c.c. partial order. Also suppose that for all $i<\cup\a$,
    $\{(\dot S_{ij},\d_{ij})\mid j<\k^+\}$ is the list of all pairs $(\dot S,\d)$ such that 
    $\dot S$ is a nice $\P_i$-name for a subset of $\check\k\setminus \check E$ and $\d<\k$, and suppose that
    $$g_{\a}\colon \{\dot S_{f(i)}\mid i<\a\}\to F_0\eqno(***)$$ 
    is an injective function, where $F_0$ is defined at $(*)$.

    If $\a$ is a limit, let
    $\P_\a$ consist of those $p\colon \a\to \Cup_{i<\a} \dom \s_i$ with \mbox{$|\sprt(p)|<\k$} (support,
    see Section \ref{ssec:Functions} on page \pageref{ssec:Functions}) such that for all $\g<\a$,
    $p\rest \g\in \P_\g$ and let $g_\a=\Cup_{i<\a}g_i$. 
    Suppose $\a$ is a successor,
    $\a=\g+1$. Let $\{(\dot S_{\g j},\d_{\g j})\mid j<\k\}$ be the the list of pairs as defined above.
    Let $(\dot S, \d) = (\dot S_{f(\g)},\delta_{f(\g)})$ where $f$ is the bijection defined above.
    If there exists $i<\g$ such that $\dot S_{f(i)}=\dot S_{f(\g)}$ (i.e. $\dot S_i$ has been already under focus),
    then let $g_\a=g_\g$. Otherwise let 
    $$g_\a=g_\g\cup \{(\dot S_{f(\g)},\eta)\}.$$ 
    where $\eta$ is some element in $F_0\setminus \ran g_\g$. Doing this, we want to make sure
    that in the end $\ran g_{\k^+}=F_0$. We omit the technical details needed to ensure that.

    Denote $\eta=g(\dot S_{f(\g)})$. 
    Let $\s_\g$ be a $\P_\g$-name such that for all $\P_{\g}$-generic $G_\g$ it holds that
    $$\P_{\g}\forces        \begin{cases}  
      \s_\g=T(\dot S_{\eta+2\d}),   &\text{ if } V[G_\g]\models [(\d_{f(\g)}\in \dot S_{f(\g)})\land (\dot S_{f(\g)}\text{ is stationary})] \\
      \s_\g=T(\dot S_{\eta+2\d+1}), &\text{ if } V[G_\g]\models [(\d_{f(\g)}\notin \dot S_{f(\g)})\land (\dot S_{f(\g)}\text{ is stationary})] \\
      \s_\g=\{\check \es\},         &\text{ otherwise.}
    \end{cases}
    $$
    Now let $\P_\a$ be the collection of sequences $p=\la\rho_i\ra_{i\le\g}$
    such that $p\rest\g=\la \rho_i\ra_{i<\g}\in \P_\g,$ $\rho_\g\in \dom \s_\g$ 
    and $p\rest\g\forces_{\P_\g}\rho_{\g}\in \s_\g$
    with the ordering defined in the usual way.

    Let $G$ be $\P_{\k^+}$-generic.
    Let us now show that  
    the extension $V[G]$ satisfies what we want, namely that $S\subset \k\setminus E$ is stationary if and only if
    there exists $\eta\in 2^E$ such that $S=X_\eta$ (Claims 3 and~4 below). 
    \vspace{-1pt}
    \begin{claim}{1}
      For $\a\le\k^+$ the forcing $\P_{\a}$ does not add bounded subsets of $\k$ and the suborder
      $$\Q_{\a}=\{p\mid p\in \P_{\a}, p=\la \check\rho_i\ra_{i<\a}\text{ where }\rho_i\in V\text{ for }i<\a\}$$
      is dense in $\P_{\a}$. 
    \end{claim}
    \begin{proofVOf}{Claim 1}
      Let us show this by induction on $\a\le \k^+$.
      For $\P_{0}$ this is already proved and the limit case is left to the reader.
      Suppose this is proved for all $\g<\a<\k^+$ and $\a=\b+1$.
      Then suppose $p\in \P_{\a}$, 
      $p=\la\rho_i\ra_{i<\a}$. Now $p\rest \b\forces \rho_\b\in \s_\b$.
      Since by the induction hypothesis $\P_{\b}$ does not add bounded subsets of $\k$ and $\Q_\b$ is dense in $\P_\b$, 
      there exists a condition $r\in \Q_\b$, $r>p\rest\b$ and a standard name $\check q$ such that $r\forces \check q=\rho_\b$.
      Now $r\cat (\check q)$ is in $\Q_{\a}$, so it is dense in $\P_{\a}$.
      To show that $\P_{\a}$ does not add bounded sets, it is enough to show that $\Q_\a$ does not. 
      Let us think of $\Q_{\a}$ as a suborder of the product $\prod_{i<\a}2^{<\k}$.
      Assume that $\tau$
      is a $\Q_\a$-name and $p\in \Q_\a$ forces that $|\tau|=\check\l<\check \k$ for some cardinal $\l$. Then let 
      $\la M_\d\ra_{\d<\k}$ be a sequence of elementary submodels of $H(\k^+)$ such that
      for all~$\d$,~$\b$
      \begin{myAlphanumerate}
      \item $|M_\d|<\k$
      \item $\d<\b\Rightarrow M_{\d}\preceq M_{\b}$, \label{proofofclaim1item2}
      \item $M_{\d}\cap \k\subset M_\d$,
      \item if $\b$ is a limit ordinal, then $M_{\b}=\Cup_{\a<\b}M_\a$,
      \item if $\k=\l^+$, then $M_\d^{<\l}\subset M_{\d}$ and if $\k$ is inaccessible, then $M_\d^{|M_\d|}\subset M_{\d+1}$,\label{proofofclaim1item5}
      \item $M_\a\in M_{\a+1}$,\label{proofofclaim1item6}
      \item $\{p,\k,\Q_{\a},\tau,\check E\}\subset M_{0}$.
      \end{myAlphanumerate}
      This (especially (e)) is possible since $\k$ is not a successor of a singular cardinal and GCH holds.
      Now the set $C=\{M_\d\cap \k\mid \d<\k\}$ is cub, so because $\k\setminus E$ is fat, there is a closed sequence $s$ of length
      $\l+1$ in $C\setminus E$. Let $(\d_i)_{i\le\l}$ be the sequence such that $s=\la M_{\d_i}\cap \k\ra_{i\le\l}$.
      For $q\in \Q_\a$, let
      $$m(q)=\inf_{\g\in\sprt q}\ran q(\g).\eqno(\star)$$

      Let $p_0=p$ and for all $i<\g$ let $p_{i+1}\in M_{\d_{i+1}}\setminus M_{\d_i}$ be
      such that $p_{i}<p_{i+1}$, $p_{i+1}$ decides $i+1$ first values of $\tau$ (think of $\tau$ as a name for a function
      $\l\to\k$ and that $p_{i}$ decides the first $i$ values of that function) and 
      $m(p_{i+1})\ge M_{\d_i}\cap \k$. 
      This $p_{i+1}$ can be found because clearly $p_i\in M_{\d_{i+1}}$ and $M_{\d_{i+1}}$ is an elementary submodel. 
      If $i$ is a limit, $i<\l$, then 
      let $p_i$ be an upper bound of $\{p_j\mid j<i\}$ which can be found in $M_{\d_{i+1}}$  
      by the assumptions  (f), (e) and (b), and because
      $M_{\d_{i}}\cap\k\notin E$. Finally let $p_\l$ be an upper bound of $\la p_{i}\ra_{i<\l}$ which exists because 
      for all $\a\in \Cup_{i<\l}\sprt p_i$
      $\sup_{i<\l}\ran p_i(\a)=M_{\d_{\l}}\cap\k$ is not in $E$ and the forcing is closed under such sequences. So
      $p_{\l}$ decides the whole $\tau$. This completes the proof of the claim.
    \end{proofVOf}
    So for simplicity, instead of $\P_{\k^+}$ let us work with $\Q_{\k^+}$.

    \begin{claim}{2} Let $G$ be $\P_{\k^+}$-generic over $V$.
      Suppose $S\subset \k$, $S\in V[G]$ and $\dot S$ is a nice name for a subset of $\k$ such that $\dot S_G=S$.
      Then let $\g$ be the smallest ordinal with $S\in V[G_\g]$. If 
      $(S\subset \k\setminus E\text{ is stationary})^{V[G_\g]}$, then $S$ is stationary in $V[G]$.
      If $\dot S=\dot S_\eta$ for some $\eta\in V$ and $V[G_{\g}]\models \s_{\g}\ne T((\dot S_\eta)_{G_{\g}\restl\{0\}})$ 
      for all $\g<\k^+$, then $S$ is stationary in $V[G]$. 
    \end{claim}
    \begin{proofVOf}{Claim 2}
      Recall, $\s_\g$ is as in the construction of $\P_{\k^+}$.
      Suppose first that $S\subset \k\setminus E$ is a stationary set in $V[G_\g]$ for some $\g<\k^+$. 
      Let us show that $S$ is stationary in $V[G]$. Note
      that $V[G]=V[G_\g][G^{\g}]$ where $G^\g=G\rest\{\a\mid \a\ge\g\}$. Let us show this 
      in the case $\g=0$ and $S\in V$, the other cases being similar. 
      Let $\dot C$ be a name and $p$ a condition which forces that $\dot C$ is cub. 
      Let us show that then $p\forces \check S\cap\dot C\ne\check\es$.
      For $q\in \Q_{\k^+}$ let $m(q)$ be defined as in $(\star)$ above.

      Like in the proof of Claim 1, construct a continuous increasing 
      sequence $\la M_{\a}\ra_{\a<\k}$ of elementary submodels
      of $H(\k^{++})$ such that $\{p,\k,\P_{\k^+},\check S,\dot C\}\subset M_{0}$
      and $M_{\a}\cap\k$ is an ordinal. Since
      $\{M_\a\cap \k\mid \a<\k, M_\a\cap\k=\a\}$ is cub, 
      there exists $\a\in S$ such that $M_\a\cap \k=\a$ and because $E$ does not reflect to $\a$
      there exists a cub sequence
      $$c\subset \{M_\b\cap \k\mid \b<\a, M_\b\cap\k=\b\}\setminus E,$$ 
      $c=\la c_i\ra_{i<\cf(\a)}$.
      Now, similarly as in the proof of Claim 1, we can choose an increasing $\la p_i\ra_{i\le\cf(\a)}$
      such that $p_0=p$,
      $p_i\in \Q_{\k^+}$ for all $i$, $p_{i+1}\forces \check\b\in \dot C$ for some 
      $c_i\le\b\le c_{i+1}$, $p_{i+1}\in M_{c_{i+1}}\setminus M_{c_i}$ and
      $m(p_{i+1})\ge c_i$. If $i$ is a limit, let $p_i$ be again 
      an upper bound of $\{ p_j\mid {j<i}\}$ in $M_{c_i}$. 
      Since the limits are not in $E$, the upper bounds exist. Finally
      $p_{\cf(\a)}\forces \a\in \dot C$, which implies $p_{\cf(\a)}\forces \check S\cap \dot C\ne\es$,
      because $\a$ was chosen from $S$.    
      
      Assume then that $\dot S=\dot S_\eta$ for some $\eta\in V$ such that 
      $$V[G_{\g}]\models \s_\g\ne T((\dot S_\eta)_{G_{\g}\restl \{0\}})$$
      for all $\g<\k^+$.
      To prove that $(\dot S_\eta)_{G}$ is stationary in $V[G]$, we carry the same argument as the above, a little modified. 
      Let us work in $V[G_0]$ and let $p_0$ force that 
      $$\forall \g< \k^+(\sigma_\g\ne T(S_\eta)).$$
      (This $p_0$ exists for example because there is at most one $\g$ such that $\s_\g=T(S_\eta)$)
      Build the sequences $c$, $\la M_{c_i}\ra_{i<\cf(\a)}$ and 
      $\la p_i\ra_{i<\cf(\a)}$ in the same fashion as above, except that assume additionally
      that the functions $g_{\k^+}$ and $f$, defined along with $\P_{\k^+}$, are in $M_{c_0}$. 
      
      At the successor steps one has to choose $p_{i+1}$ such that for each $\g\in \sprt p_i$,
      $p_{i+1}$ decides $\sigma_\g$. This is possible, since there are only three choices for
      $\sigma_\g$, namely $\{\es\}$, $T(S_{\xi+2\a+1})$ or $T(S_{\xi+2\a})$ where $\xi$ and $\a$ are justified by
      the functions $g_{\k^+}$ and $f$. For all $\g\in\sprt p_i$ let us denote by $\xi_\g$ the function such that
      $p_{i+1}\rest\g\forces \sigma_\g=T(S_{\xi_\g})$. Clearly $\eta\ne \xi_\g$ for all $\g\in\sprt p_i$.
      Further demand that
      $m(p_{i+1})>\sup (S_{\eta}\cap S_{\xi_{\g}})$ for all $\g\in \sprt p_i$.
      It is possible to find such $p_{i+1}$ from $M_{i+1}$ because $M_{i+1}$ is an elementary submodel
      and such can be found in $H(\k^{++})$ since $\xi_\g\ne \eta$ and by the definitions $S_\eta\cap S_{\xi_\g}$ is bounded.
    \end{proofVOf}
    \begin{claim}{3}
      In $V[G]$ the following holds: if $S\subset \k\setminus E$ 
      is stationary, then there exists $\eta\in 2^{E}$ with $\eta(0)=0$ such that $S=X_\eta$.
    \end{claim}
    \begin{proofVOf}{Claim 3}
      Recall the function $g_{\k^+}$ from the construction of $\P_{\k^+}$ (defined at $(***)$ and the paragraph below that). 
      Let $\eta=g_{\k^+}(\dot S)$ where
      $\dot S$ is a nice name $\dot S\in V$ such that $\dot S_G=S$.
      If $\a\in S$, then there is the smallest $\g$ such that $\dot S=S_{f(\g)}$ and $\a=\delta_{f(\g)}$ (where $f$ is as in the 
      definition of $\P_{\k^+}$). This stage
      $\g$ is the only stage where it is possible that $V[G_\g]\models \sigma_\g=T(S_{\eta+2\a+1})$, but since $V[G_\g]\models \check \a\in \dot S$,
      by the definition of 
      $\P_{\k^+}$ it is not the case, so
      the stationarity  of $S_{\eta+2\a+1}$ has not been killed by Claim~2. On the other hand the stationarity of $S_{\eta+2\a}$
      is killed at this level $\g$ of the construction, so $\a\in X_\eta$ by the 
      definitions of $\f$ and $X_\eta$. Similarly if $\a\notin S$, we conclude that
      $\a\notin X_\eta$.
    \end{proofVOf}
    \begin{claim}{4}
      In $V[G]$ the following holds: if $S\subset \k\setminus E$ is not stationary, then for all $\eta\in 2^{E}$ with $\eta(0)=0$ we have
      $S\ne X_\eta$.      
    \end{claim}
    \begin{proofVOf}{Claim 4}
      It is sufficient to show that $X_\eta$ is stationary for all $\eta\in 2^{E}$ with $\eta(0)=0$.
      Suppose first that $\eta\in F_0\subset V$. Then since $g_{\k^+}$ is a surjection onto $F_0$ (see~$(***)$), 
      there exists a name $\dot S$ such that $S=\dot S_G$ is stationary,
      $S\subset \k\setminus E$ and $g_{\k^+}(S)=\eta$. 
      Now the same argument as in the proof of Claim 3 implies that $X_\eta=S$, so $X_\eta$ is stationary by Claim 2.

      If $\eta\notin F_0$, then by the definition of $\eta\mapsto X_\eta$ it is sufficient to show that the $\diamondsuit$-sequence 
      added by $\P_0$ guesses in $V[G]$ every new set on a stationary set.  

      Suppose that $\tau$ and $\dot C$ are nice $\P_{\k^+}$-names for subsets of $\check\k$ and let $p$ be a condition
      forcing that $\dot C$ is cub. We want to find $\g$ and $q>p$ such that 
      $$q\forces ((\cup\dot G_0)(\check \g)^{-1}\{1\}=\tau\cap\check\g)\land (\check\g\in\dot C)$$
      where $\dot G_0=\dot G\rest \{0\}$ is the name for the $\P_0$-generic.
      To do that let
      $p_0\ge p$ be such that $p_0\forces \tau\notin \check{\Po(\k)^{V}}$. 
      
      Similarly as in the proofs above define a suitable sequence $\la M_{i}\ra_{i<\l}$ of elementary submodels, 
      of length $\l<\k$, where
      $\l$ is a cofinality of a point in $E$, such that $\sup_{i<\l}(M_{i}\cap\k)=\a\in E$ and $M_i\cap\k\notin E$
      for all $i<\l$. 
      Assume also that $p_0\in M_0$.
      Suppose $p_i\in M_{i}$ is defined. Let $p_{i+1}>p_i$ be an element of $M_{i+1}\setminus M_i$
      satisfying the following:
      \begin{myEnumerate}
      \item $p_{i+1}$ decides $\sigma_\b$ for all $\b\in\sprt p_i$,
      \item for all $\b\in \sprt p_i$ there is $\b'\in M_{i+1}$ such that $p_{i+1}\forces \b'\in\tau\sd\xi_\b$,
        where $\xi_\b$ is defined as in the proof of Claim 2 and $p_{i+1}$ decides what it~is,
      \item $p_{i+1}$ decides $\tau$ up to $M_i\cap \k$,
      \item $p_{i+1}\forces \d\in \dot C$ for some $\d\in M_{i+1}\setminus M_i$,
      \item $m(p_{i+1})>M_i\cap\k$, ($m(p)$ is defined at $(\star)$),
      \end{myEnumerate}
      Item (1) is possible for the same reason as in the proof of Claim 2 and (2) is possible since
      $p_{i}\forces \forall\eta\in \check{\Po(\k)^V}(\tau\ne S_{\check \eta}).$
      
      Since $M_i\cap\k\notin E$ for $i<\l$, 
      this ensures that the sequence $p_0\le p_1\le \dots$ closes under 
      limits $<\l$. Let $p_\l=\Cup_{i<\l}p_i$ and let us define $q\supset p_\l$ as follows: 
      $\sprt q=\sprt p_\l$, for $\d\in \sprt p_\l\setminus \{0\}$ 
      let $\dom q=\a+1$, $p_\l(\d)\subset q(\d)$, $q(\a)=1$
      and $q(0)(\a)=\tau\cap \g$ ($\tau$~means here what have been decided by $\{p_i\mid i<\l\}$).
      Now $q$ is a condition in the forcing notion. 
      
      Now certainly, if $q\in G$, then 
      in the extension $\tau_G\cap \a=(\cup G_0)(\a)^{-1}\{1\}$ and $\a\in C$, so we finish. 
    \end{proofVOf}
  \end{proofVOf}

  \begin{proofVOf}{item \eqref{thm:CUBSETItem3}}  
    If $\k=\l^+$, this follows from the result of Mekler and Shelah \cite{MekShe} 
    and Hyttinen and Rautila \cite{HytRau} that the existence of a $\k\l$-canary tree is consistent.  
    For arbitrary $\l<\k$ the result follows from the item \eqref{thm:CUBSETItem4} of this theorem proved above 
    (take $Z=\k\setminus S^\k_{\l}$).
  \end{proofVOf}
  \begin{proofVOf}{item \eqref{thm:CUBSETItem5}}
    For $X=\k$ this was proved by Halko and Shelah in \cite{HalShe}, Theorem~4.2.
    For $X$ any stationary subset of $\k$ the proof is similar. It is sufficient to show that 
    $2^\k\setminus \CUB(X)$ is not meager in any open set. Suppose $U$ is an open set and $(D_\a)_{\a<\k}$
    is a set of dense open sets and let us show that 
    $$(2^\k\setminus \CUB(X))\cap U\cap\Cap_{\a<\k} D_{\a}\ne\es.$$ Let $p\in 2^{<\k}$
    be such that $N_{p}\subset U$. Let $p_0\ge p$ be such that $p_0\in D_0$. Suppose $p_\b$ are defined for
    $\b<\a+1$. Let $p_{\a+1}$ be such that $p_{\a+1}\ge p_\a$, $p_{\a+1}\in D_{\a+1}$. Suppose $p_{\b}$ is defined for
    $\b<\a$ and $\a$ is a limit ordinal. Let $p_\a$ be any element of $2^{<\k}$ such that $p_\a > \Cup_{\b<\a}p_\b$,
    $p_\a(\sup\limits_{\b<\a}\dom p_\b)=0$ and $p_\a\in D_\a$. Let $\eta=\Cup_{\a<\k}p_\a$. 
    The complement of $\eta^{-1}\{1\}$ contains a cub, so $X\setminus \eta^{-1}\{1\}$ is stationary whence
    $\eta\notin \CUB(X)$ and so $\eta\in 2^\k\setminus \CUB(X)$. Also clearly $\eta\in U\cap \Cap_{\a<\k}D_\a$.
  \end{proofVOf}
  \begin{proofVOf}{item \eqref{thm:CUBSETItem6}}
    Our proof is different from that given by L\"ucke and Schlicht.
    Suppose $\k^{<\k}=\k >\o$. 
    We will show that in a
    generic extension of $V$ all $\Dii$-sets have the Property of Baire. Let
    $$\P=\{p\mid p\text{ is a function,} |p|<\k, \dom p\subset \k\times\k^+, \ran p\subset \{0,1\}\}$$
    with the ordering $p<q\iff p\subset q$
    and let $G$ be $\P$-generic over $V$. Suppose that $X\subset 2^{\k}$ is a $\Dii$-set in $V[G]$. It is sufficient to show
    that for every $r\in 2^{<\k}$ there is $q\supset r$ such that either $N_q\setminus X$ or $N_q\cap X$ is co-meager.
    So let $r\in 2^{<\k}$ be arbitrary.
    
    Now suppose that $\la p_i\ra_{i<\k}$ and $\la q_i\ra_{i<\k}$ are sequences in $V[G]$ such that $p_i,q_i\in (2^{<\k})^{2}$
    for all $i<\k$ and $X$ is the projection of 
    $$C_0=(2^{\k})^{2}\setminus \Cup_{i<\k}N_{p_i}$$
    and $2^{\k}\setminus X$ is the projection of 
    $$C_1=(2^{\k})^{2}\setminus \Cup_{i<\k}N_{q_i}.$$
    (By $N_{p_i}$ we mean $N_{p_i^1}\times N_{p_i^2}$ where $p_i=(p_i^1,p_i^2)$.)
    Since these sequences have size $\k$, there exists $\a_{1}<\k^+$ such that they are already in
    $V[G_{\a_1}]$ where $G_{\a_1}=\{p\in G\mid \dom p\subset \k\times \a_1\}$.
    More generally, for $E\subset \P$ and $A\subset \k^+$,
    we will denote $E_A=\{p\in E\mid \dom p\subset \k\times A\}$ and
    if $p\in \P$, similarly $p_{A}=p\rest (\k\times A)$.

    Let $\a_2\ge \a_1$ be such that $r\in G_{\{\a_2\}}$ (identifying $\k\times \{\a_2\}$ with $\k$). 
    This is possible since $G$ is generic.
    Let $x=G_{\{\a_2\}}$. Since in $V[G]$, $x\in X$ or $x\in 2^\k\setminus X$, there are $\a_3>\a_2$, 
    $p\in G_{\a_3}$, $p_{\{\a_2\}}\supset r$ and a name $\tau$ such that $p$ forces 
    that $(x,\tau)\notin N_{p_i}$ for all $i<\k$ or $(x,\tau)\notin N_{q_i}$ for all $i<\k$.
    Without loss of generality assume
    that $p$ forces that $(x,\tau)\notin N_{p_i}$ for all $i<\k$. Also we can assume that
    $\tau$ is a $\P_{\a_3}$-name and that $\a_3=\a_2+2$.

    By working in $V[G_{\a_{2}}]$ we may assume that
    $\a_{2}=0$. For all $q\in \P_{\{ 1\}}$,
    $p_{\{ 1\}}\subseteq q$ and $i<\k$, let
    $D_{i,q}$ be the set of all
    $s\in \P_{\{ 0\}}$ such that $p_{\{ 0\}}\subseteq s$,
    $\dom(s)\ge \dom(p_{i}^{1})$
    and there is $q'\in \P_{\{ 1\}}$ such that
    $q\subseteq q'$ and $s\cup q'$ decides $\tau\rest \dom(p_{i}^{2})$.
    Clearly each $D_{i,q}$ is dense above $p_{\{0\}}$ in $\P_{\{ 0\}}$
    and thus it is enough to show that
    if $y\in 2^{\k}$ is such that for all
    $i <\k$ and $q$ as above there is $\a <\k$ such that
    $y\rest\a\in D_{i,q}$, then
    $y\in X$.
    
    So let $y$ be such. Then we can find $z\in 2^{\k}$ such that
    for all $i<\k$ and $q$ as above there are $\a ,\b <\k$
    such that $\a\ge \dom(p_{i}^{1})$ and
    $y\rest\a \cup z\rest\b$ 
    decides $t=\tau\rest \dom(p_{i}^{2})$. By the choise of $p$,
    $(y\rest \dom(p_{i}^{1}),t)\ne p_{i}$.
    Thus letting $\tau^{*}$ be the function as decided by $y$ and $z$,
    $(y,\tau^{*})\in C_{0}$ and thus $y\in X$.
  \end{proofVOf}
\end{proofVOf}

\begin{Remark}[$\cf(\k)=\k>\o$]
  There are some more results and strengthenings of the results in Theorem \ref{thm:CUBSET}:
  \begin{myEnumerate}
  \item (Independently known by S. Coskey and P. Schlicht) 
    If $V=L$ then there is a $\Delta^1_1$ wellorder of $\Po(\kappa)$ and this implies that there is a $\Delta^1_1$ set wihtout the Baire Property.
  \item Suppose that $\omega < \k < \l$, $\k$ regular and $\l$
    inaccessible. Then after turning $\l$ into $\k^+$ by collapsing
    each ordinal less than $\l$ to $\k$ using conditions of size $<
    \k$, the Baire Property holds for $\Dii$ subsets of $\k^\k$.
  \end{myEnumerate}
\end{Remark}

\begin{Cor} \label{cor:NSIsNotBorel}
  For a regular $\l<\k$
  let $\NS_\l$ denote the equivalence relation on $2^\k$ such that $\eta\NS_\l\xi$ if and only if
  $\eta^{-1}\{1\}\sd\xi^{-1}\{1\}\text{ is not }\l\text{-stationary}.$
  Then $\NS_\lambda$ is not Borel and it is not $\Delta^1_1$ in $L$ or in the forcing extensions after 
  adding $\kappa^+$ Cohen subsets of $\kappa$.
\end{Cor}
\begin{proof}
  Define a map $f\colon 2^\k\to (2^\k)^2$ by
  $\eta\mapsto (\es,\k\setminus \eta)$. Suppose for a contradiction that $\NS_\l$
  is Borel. Then 
  $$\NS_\es=\NS_\l\cap \underbrace{\{(\es,\eta)\mid \eta\in 2^{\k}\}}_{\text{closed}}$$
  is Borel, and further $f^{-1}[\NS_{\es}]$ is Borel by continuity of $f$. But
  $f^{-1}[\NS_{\es}]$ equals $\operatorname{CUB}$ which is not Borel by Theorem \ref{thm:CUBSET} 
  \eqref{thm:CUBSETItem5} and Theorem \ref{thm:BorelNoBaire}.
  Similarly, using items \eqref{thm:CUBSETItem2} and \eqref{thm:CUBSETItem7} of Theorem \ref{thm:CUBSET}, one 
  can show that $\NS_\l$ is not~$\Dii$ under the stated assumptions.
\end{proof}

\section{Equivalence Modulo the Non-stationary Ideal}

In this section we will investigate the relations defined as follows:

\begin{Def}
  For $X\subset \k$, we denote by $E_X$ the relation 
  $$E_X=\{(\eta,\xi)\in 2^{\k}\times 2^{\k}\mid (\eta^{-1}\{1\}\sd \xi^{-1}\{1\})\cap X\text{ is not stationary}\}.$$
\end{Def}

The set $X$ consists usually of ordinals of fixed cofinality, i.e. $X\subset S^\k_\mu$ for some $\mu$.
These relations are easily seen to be $\Sii$. If $X\subset S^\k_\o$, then it is in fact Borel*. To see this 
use the same argument as in the proof of Theorem \ref{thm:CUBSET} \eqref{thm:CUBSETItem1} 
that the $\operatorname{CUB}^{\k}_{\o}$-set is Borel*.

\subsection{An Antichain}

\begin{Thm}\label{thm:AntiChainChain}
  Assume GCH, $\kappa^{<\k}=\k$ is uncountable and $\mu<\k$ is a regular
  cardinal such that if $\k=\l^+$, then $\mu\le\cf(\l)$. 
  Then in a cofinality 
  and GCH preserving forcing extension,
  there are stationary sets $K(A)\subset S^\k_\mu$ for each $A\subset \k$ such that
  $E_{K(A)}\not\le_B E_{K(B)}$ if and only if $A\not\subset B$.
\end{Thm}
\begin{proofV}{Theorem \ref{thm:AntiChainChain}}
  In this proof we identify functions $\eta\in 2^{\le\k}$ with the sets $\eta^{-1}\{1\}$: for example
  we write $\eta\cap \xi$ to mean $\eta^{-1}\{1\}\cap \xi^{-1}\{1\}$.

  The embedding will look as follows. Let $(S_{i})_{i<\k}$ be pairwise disjoint stationary subsets of 
  $$\lim S^\k_\mu = \{\a\in S^\k_\mu\mid \a\text{ is a limit of ordinals in }S^\k_\mu\}.$$
  Let
  $$K(A)=E_{\displaystyle\mathop{\cup}_{\tiny{\a\in A}}S_\a}.\eqno(*)$$

  If $X_1\subset X_2\subset \k$, then $E_{X_1}\le_B E_{X_2}$, because $f(\eta)=\eta\cap X_1$ is a reduction. This
  guarantees that 
  $$A_1\subset A_2\Rightarrow K(A_1)\le_B K(A_2).$$ 
  
  Now suppose that for all $\alpha < \kappa$ we have killed (by forcing) all reductions from 
  $K({\alpha}) = E_{S_\alpha}$ to $K(\kappa \setminus {\alpha}) =
  E_{\bigcup_{\beta \neq \alpha}S_\beta}$ for all $\alpha < \kappa$. Then if $K(A_1) \le_B K(A_2)$ it follows 
  that $A_1 \subset A_2$: Otherwise choose $\alpha \in A_1 \setminus A_2$ and we have:
  $$K({\alpha}) \le_B K(A_1) \le_B K(A_2) \le_B K(\kappa \setminus{\alpha}),$$
  contradiction. So we have:
  $$A_1 \subset A_2 \iff K(A_1) \le_B K(A_2).$$
  It is easy to obtain an antichain of length $\k$ in $\Po(\k)$ and so the result follows.

  Suppose that $f\colon E_X \le_B E_Y$ is a Borel reduction. Then $g\colon 2^{\k}\to 2^{\k}$
  defined by $g(\eta)=f(\eta)\sd f(0)$ is a Borel function with the following property:
  $$\eta\cap X\text{ is stationary }\iff g(\eta)\cap Y\text{ is stationary}.$$
  The function $g$ is Borel, so by Lemma \ref{lem:BaireCont}, page~\pageref{lem:BaireCont},
  there are dense open sets
  $D_i$ for $i<\k$ such that $g\rest D$ is continuous where $D=\Cap_{i<\k} D_i$. Note that
  $D_i$ are open so for each $i$ we can write $D_i=\Cup_{j<\k}N_{p(i,j)}$, where $(p(i,j))_{j<\k}$
  is a suitable collection of elements of $2^{<\k}$.

  Next define $Q_g\colon 2^{<\k}\times 2^{<\k}\to \{0,1\}$ by $Q_g(p,q)=1\iff N_p\cap D\subset g^{-1}[N_q]$
  and $R_g\colon \k\times\k\to 2^{<\k}$ by $R_g(i,j)=p(i,j)$ where $p(i,j)$ are as above. 

  For any $Q\colon 2^{<\k}\times 2^{<\k}\to \{0,1\}$
  define 
  $Q^*\colon 2^{\k}\to 2^{\k}$ by
  $$
  Q^*(\eta)=
  \begin{cases}
    \xi,&\text{ s.t. }\forall\a<\k\exists\b<\k Q(\eta\rest\b,\xi\rest\a)=1\text{ if such exists,}\\
    0,&\text{ otherwise}.
  \end{cases}
  $$
  And for any $R\colon \k\times\k\to 2^{<\k}$ define
  $$R^*=\Cap_{i<\k}\Cup_{j<\k}N_{R(i,j)}.$$
  
  Now clearly $R_g^*=D$ and $Q_g^*\rest D=g\rest D$, i.e. $(Q,D)$ \emph{codes} $g\rest D$ in this sense. 
  Thus we have shown that if there is a reduction
  $E_X\le_B E_Y$, then there is a pair $(Q,R)$ which satisfies the following conditions:
  \begin{myEnumerate}
  \item \label{QR1} $Q\colon (2^{<\k})^2\to \{0,1\}$ is a function.
  \item \label{QR2} $Q(\es,\es)=1$,
  \item \label{QR3} If $Q(p,q)=1$ and $p'>p$, then $Q(p',q)=1$,
  \item \label{QR4} If $Q(p,q)=1$ and $q'< q$, then $Q(p,q')=1$
  \item \label{QR5} Suppose $Q(p,q)=1$ and $\a>\dom q$. There exist $q'>q$ and $p'>p$ such that $\dom q'=\a$ and
    $Q(p',q')=1$,
  \item \label{QR6} If $Q(p,q)=Q(p,q')=1$, then $q\le q'$ or $q'<q$,
  \item \label{QR7} $R\colon \k\times\k\to 2^{<\k}$ is a function.
  \item \label{QR8} For each $i\in \k$ the set $\Cup_{j<\k}N_{R(i,j)}$ is dense.
  \item \label{QR9} For all $\eta\in R^*$, $\eta\cap X$ is stationary if and 
    only if $Q^*(\eta\cap X)\cap Y$ is stationary.
  \end{myEnumerate}
  
  Let us call a pair $(Q,R)$ which satisfies \eqref{QR1}--\eqref{QR9} \emph{a code for a reduction (from $E_X$ to $E_Y$)}. 
  Note that it is not the same as the Borel code for the graph of a reduction function as a set.
  Thus we have shown that if
  $E_X\le_B E_Y$, then there exists a code for a reduction from $E_X$ to $E_Y$.
  We will now prove the following lemma which is stated in a general enough form so we can use it also in the next section:

  \begin{Lemma}[GCH]\label{lem:KillOneCode} 
    Suppose $\mu_1$ and $\mu_2$ are regular cardinals less than $\kappa$ such that if $\k=\l^+$, then $\mu_2\le\cf(\l)$, and
    suppose $X$ is a stationary subset of $S^\kappa_{\mu_1}$, $Y$ is a subset of
    $S^\kappa_{\mu_2}$, $X \cap Y = \es$ (relevant if $\mu_1 = \mu_2$) and if
    $\mu_1 < \mu_2$ then $\alpha \cap X$ is not stationary in $\alpha$ for all
    $\alpha \in Y$.
    Suppose that $(Q,R)$ is an arbitrary pair.
    Denote by $\f$ the statement ``$(Q,R)$ is not a code for a reduction from $E_X$ to $E_Y$''. 
    Then there is a $\k^+$-c.c. $<\k$-closed forcing $\R$ such that $\R\forces \f$.
  \end{Lemma}
  \begin{Remark}
    Clearly if $\mu_1=\mu_2=\o$, then the condition $\mu_2\le \cf(\l)$ is of course true.
    We need this assumption in order to have $\nu^{<\mu_2}<\k$ for all $\nu<\k$.
  \end{Remark}
  \vspace{-10pt}
  \begin{proofVOf}{Lemma \ref{lem:KillOneCode}}
    We will show that one of the following holds: 
    \begin{myEnumerate}
    \item $\f$ already holds, i.e. $\{\es\}\forces \f$,
    \item $\P=2^{<\k}=\{p\colon\a\to 2\mid \a<\k\}\forces \f$,
    \item $\R\forces \f$,
    \end{myEnumerate}
    where $$\R=\{(p,q)\mid p,q\in 2^{\a},\a<\k, X\cap p\cap q=\es, q\text{ is }\mu_1\text{-closed}\}$$
    Above ``$q$ is $\mu_1$-closed'' means ``$q^{-1}\{1\}$ is $\mu_1$-closed'' etc., and we will use this abbreviation 
    below.
    Assuming that (1) and (2) do not hold, we will show that (3) holds. 
    
    Since (2) does not hold,
    there is a $p\in \P$ which forces $\lnot \f$ and so
    $\P_p=\{q\in \P\mid q>p\}\forces\lnot\f$. But $\P_p\cong \P$, so in fact $\P\forces \lnot \f$,
    because $\f$ has only standard names as parameters (names for elements in $V$, such as $Q$, $R$, $X$ and $Y$).
    Let $G$ be any $\P$-generic and let us denote the set $G^{-1}\{1\}$ also by $G$. 
    Let us show that $G\cap X$ is stationary. Suppose that $\dot C$ is a name and $r\in\P$ is a condition which forces that
    $\dot C$ is cub. For an arbitrary $q_0$, let us find a $q>q_0$ which forces $\dot C\cap \dot G\cap \check X\ne\es$.
    Make a counter assumption: no such $q>q_0$ exists.
    Let $q_1>q_0$ and $\a_1>\dom q_0$ be such that $q_1\forces \check\a_1\in \dot C$, $\dom q_1>\a_1$ is a successor and $q_1(\max\dom q_1)=1$.
    Then by induction on $i<\k$ let $q_{i+1}$ and $\a_{i+1}>\dom q_{i}$ be such that $q_{i+1}\forces \check\a_{i+1}\in \dot C$,
    $\dom q_{i+1}>\a_{i+1}$ is a successor and $q_{i+1}(\max\dom q_{i+1})=1$. If $j$ is a limit ordinal, let 
    $q_j=\Cup_{i<j}q_i\cup \{(\sup_{i<j}\dom q_i,1)\}$ and $\a_j=\sup_{i<j}\a_i$.
    We claim that for some $i<\k$, the condition $q_i$ is as needed, i.e.
    $$q_i\forces \dot G\cap \check X\cap \dot C\ne\es.$$
    Clearly for limit ordinals $j$, we have $\a_j=\max\dom q_j$ and $q_j(\a_j)=1$ and $\{\a_j\mid j\text{ limit}\}$ is cub. 
    Since $X$ is stationary, there exists
    a limit $j_0$ such that $\a_{j_0}\in X$. Because
    $q_0$ forces that $\dot C$ is cub, $q_j>q_i>q_0$ for all $i<j$, 
    $q_i\forces\check\a_i\in \dot C$ and $\a_j=\sup_{i<j}\a_i$, we have $q_j\forces \a_{j}\in\dot C\cap \check X$. On 
    the other hand $q_j(\a_j)=1$, so $q_j\forces \a_j\in G$ so we finish.    

    So now we have in $V[G]$ that $G\cap X$ is stationary, $G\in R^*$ (since $R^*$ is co-meager) 
    and $Q$ is a code for a reduction, so $Q^*$ has the property \eqref{QR9} and $Q^*(G\cap X)\cap Y$ is stationary.
    Denote $Z=Q^*(G\cap X)\cap Y$.
    We will now construct a forcing $\Q$ in $V[G]$ such that
    $$V[G]\models (\Q\forces \text{``}G\cap X\text{ is not stationary, but }Z\text{ is stationary''}).$$
    Then $V[G]\models (\Q\forces\f)$ and hence $\P*\Q\forces\f$. On the other hand $\Q$ will be chosen such that
    $\P*\Q$ and $\R$ give the same generic extensions. 
    So let 
    $$\Q=\{q\colon\a\to 2\mid X\cap G\cap q=\es,q\text{ is }\mu_1-\text{closed}\},\eqno(***)$$
    Clearly $\Q$ kills the stationarity of $G\cap X$. Let us show that it preserves the stationarity of
    $Z$. For that purpose it is sufficient to show that for any nice $\Q$-name $\dot C$ 
    for a subset of $\k$ and
    any $p\in \Q$, if $p\forces \text{``}\,\dot C\text{ is }\mu_2\text{-cub''}$, then 
    $p\forces (\dot C\cap \check Z\ne \check \es)$.
    
    So suppose $\dot C$ is a nice name for a subset of $\k$ and $p\in \Q$ is such that
    $$p\forces \text{``}\,\dot C\text{ is cub''}$$
    Let $\l>\k$ be a sufficiently large regular cardinal and let $N$ be an elementary submodel
    of $\la H(\l),p,\dot C,\Q,\k\ra$ which has the following properties:
    \begin{myItemize}
    \item $|N|=\mu_2$
    \item $N^{<\mu_2}\subset N$
    \item $\a=\sup (N\cap \k)\in Z$ (This is possible because $Z$ is stationary).
    \end{myItemize}
    Here we use the hypothesis that $\mu_2$ is at most $\cf(\l)$ when $\k=\l^+$.
    Now by the assumption of the theorem, $\a\setminus X$ contains a  $\mu_1$-closed unbounded sequence of length $\mu_2$,
    $\la\a_i\ra_{i<\mu_2}$.
    Let $\la D_i\ra_{i<\mu_2}$ list all the dense subsets of $\Q^N$ in $N$. Let $q_0\ge p$, $q_0\in \Q^N$ be arbitrary and
    suppose $q_i\in \Q^N$ is defined for all $i<\g$. 
    If $\g=\b+1$, then define $q_\g$ to be an extension of $q_\b$
    such that $q_\g\in D_\b$ and $\dom q_\g=\a_{i}$ for some $\a_i>\dom q_\b$. To do that, for instance,
    choose $\a_i>\dom q_\b$ and define 
    $q'\supset q_\b$ by $\dom q'=\a_i$, $q(\delta)=0$ for all $\delta\in \dom q'\setminus\dom q_\b$
    and then extend $q'$ to $q_\b$ in $D_\b$. 
    If $\g$ is a limit ordinal with $\cf(\g)\ne\mu_1$, then let
    $q_\g=\Cup_{i<\g} q_i$. If $\cf(\g)=\mu_1$, let 
    $$q_\g=\big(\Cup_{i<\g} q_i\big)\cat \la\sup_{i<\g}\dom q_i,1\ra$$

    Since $N$ is closed under taking sequences of length less than $\mu_2$,
    $q_\g\in N$. Since we required elements of $\Q$ to be $\mu_1$-closed but not $\g$-closed if $\cf(\g)\ne\mu_1$,
    $q_\g\in \Q$ when $\cf(\g)\ne\mu_1$. When $\cf(\g)=\mu_1$, the limit $\sup_{i<\g}\dom q_i$ coincides with a limit
    of a subsequence of $\la\a_i\ra_{i<\mu_2}$ of length $\mu_1$, i.e. the limit is $\a_\b$ for some $\b$ since this sequence
    is $\mu_1$-closed. So by definition $\sup_{i<\g}\dom q_i\notin X$ and again $q_\g\in \Q$.

    Then $q=\Cup_{\g<\mu}q_\g$ is a $\Q^N$-generic over $N$. Since $X\cap Y=\es$, also $(X\cap G)\cap Z=\es$ 
    and $\a\notin X\cap G$. Hence $q\cat(\a,1)$ is in $\Q$. We claim that $q\forces (\dot C\cap \check Z\ne\es)$.
    
    Because $p\forces \text{``}\,\dot C\text{ is unbounded''}$,
    also $N\models (p\forces \text{``}\,\dot C\text{ is unbounded''})$ by elementarity. Assuming that $\l$ is chosen large enough,
    we may conclude that for all $\Q^N$-generic $g$ over $N$, $N[g]\models \text{``}\dot C_g\text{ is unbounded''}$,
    thus in particular $N[g]\models \text{``}\dot C_g\text{ is unbounded in }\k\text{''}$.
    Let $G_1$ be $\Q$-generic over $V[G]$ with $q\in G_1$. Then $\dot C_{G_1}\supset\dot C_{q}$
    which is unbounded in $\a$ by the above, since $\sup(\k\cap N)=\a$. Because $\dot C_{G_1}$ is $\mu_2$-cub, $\a$ is in $\dot C_{G_1}$.

    Thus $\P*\Q\forces \f$. It follows straightforwardly from the definition of iterated forcing that $\R$ is 
    isomorphic to a dense suborder of $\P*\dot\Q$ where $\dot Q$ is a $\P$-name for a partial order such that 
    $\dot \Q_{G}$ equals $\Q$ as defined in $(***)$ for any $\P$-generic $G$.
    
    Now it remains to show that $\R$ has the $\k^{+}$-c.c. and is $<\k$-closed.
    Since $\R$ is a suborder of $\P\times \P$, which has size $\k$, it trivially has the $\k^+$-c.c. 
    Suppose $(p_i,q_i)_{i<\g}$ is an increasing sequence, $\g<\k$. Then the pair
    $$(p,q)=\Big\la\Big(\Cup_{i<\g}p_i\Big)\cat \la \a,0\ra,\Big(\Cup_{i<\g}q_i\Big)\cat \la \a,1\ra\Big\ra$$
    is an upper bound.
  \end{proofVOf}

  \begin{Remark}
    Note that the forcing used in the previous proof is equivalent to $\kappa$-Cohen forcing.
  \end{Remark}
  
  \begin{Cor}[GCH]
    Let $K\colon A\mapsto E_{\Cup_{\a\in A}S_\a}$ be as in the beginning of the proof.
    For each pair $(Q,R)$ and each $\a$ there is a $<\k$-closed, $\k^+$-c.c. forcing 
    $\R(Q,R,\a)$ such that 
    $$\R(Q,R,\a)\forces \text{``}\,(Q,R)\text{ is not a code for a reduction from }K(\{\a\})\text{ to }K(\k\setminus\{\a\})\text{''}$$
  \end{Cor}
  \begin{proof}
    By the above lemma one of the choices $\R=\{\es\}$, $\R=2^{<\k}$ or
    $$\R=\{(p,q)\mid p,q\in 2^{\b},\b<\k, S_\a\cap p\cap q=\es, q\text{ is }\mu\text{-closed}\}$$
    suffices.
  \end{proof}
  
  Start with a model satisfying GCH.
  Let $h\colon \k^+\to\k^+\times\k\times\k^+$ be a bijection such that $h_3(\a)<\a$ for $\a>0$ and
  $h_3(0)=0$. Let $\P_0=\{\es\}$. For each $\a<\k$, let $\{\sigma_{\b\a 0}\mid \b<\k^+\}$
  be the list of all $\P_0$-names for codes for a reduction from $K(\{\a\})$ to $K(\k\setminus \{\a\})$.
  Suppose $\P_i$ and 
  $\{\sigma_{\b\a i}\mid \b<\k^+\}$ are defined for all $i<\g$ and $\a<\k$, where $\g<\k^+$ 
  is a successor $\g=\b+1$, $\P_i$ is $<\k$-closed and has the 
  $\k^+$-c.c.

  Consider $\sigma_{h(\b)}$. By the above corollary, the following holds:
  \begin{eqnarray*} 
    \P_\b&\forces& \big[\exists\R\in \Po(2^{<\k}\times 2^{<\k})(\R\text{ is }<\k\text{-closed, }\k^+\text{-c.c. p.o. 
      and }\\
    &&\R\forces\text{``}\,\sigma_{h(\b)}\text{ is not a code for a reduction.''})\big]   
  \end{eqnarray*}
  So there is a $\P_{\b}$-name $\rho_\b$ such that $\P_{\b}$ forces that $\rho_\b$ is as $\R$ above.
  Define 
  $$\P_\g=\{(p_i)_{i<\g}\mid ((p_i)_{i<\b}\in \P_\b)\land ((p_i)_{i<\b}\forces p_{\b}\in \rho_\b)\}.$$
  And if $p=(p_i)_{i<\g}\in \P_\g$ and $p'=(p_i')_{i<\g}\in \P_{\g}$, then 
  $$p\le_{\P_\g} p'\iff [(p_i)_{i<\b}\le_{\P_\b} (p_i')_{i<\b}]\land  [(p_i')_{i<\b}\forces (p_\b\le_{\rho_\b} p_\b')]$$
  If $\g$ is a limit, $\g\le\k^+$, let
  $$\P_{\g}=\{(p_i)_{i<\g}\mid \forall \b(\b<\g\rightarrow (p_i)_{i<\b}\in \P_{\b})\land(|\sprt (p_i)_{i<\g}|<\k)\},$$
  where $\sprt$ means support, see Section \ref{ssec:Functions} on page \pageref{ssec:Functions}.
  For every $\a$, let $\{\sigma_{\b\a\g}\mid \b<\k^+\}$ list all $\P_\b$-names for codes for a reduction.
  It is easily seen that $\P_\g$ is $<\k$-closed and has the $\k^+$-c.c. for all $\g\le \k^+$
  
  We claim that $\P_{\k^+}$ forces
  that for all $\a$, $K(\{\a\})\not\le_B K(\k\setminus\{\a\})$ which suffices by the 
  discussion in the beginning of the proof, see $(**)$ for the notation.
  
  Let $G$ be $\P_{\k^+}$-generic and let $G_\g=\text{``}\,G\cap \P_\g\text{''}$ for every $\g<\k$.
  Then $G_\g$ is $\P_{\g}$-generic.

  Suppose that in $V[G]$, $f\colon 2^{\k}\to 2^{\k}$ is a 
  reduction $K(\{\a\})\le_B K(\k\setminus\{\a\})$ 
  and $(Q,R)$ is the corresponding code for a reduction. By \cite{Kunen} Theorem VIII.5.14, there is a $\d<\k^+$ such that
  $(Q,R)\in V[G_\d]$. Let $\d_0$ be the smallest such $\d$. 

  Now there exists $\sigma_{\g\a\d_0}$, a
  $\P_{\d_0}$-name for $(Q,R)$. By the definition of $h$, there exists a $\d>\d_0$ with $h(\d)=(\g,\a,\d_0)$.
  Thus $$\P_{\d+1}\forces \text{``}\sigma_{\g\a\d_0}\text{ is not a code for a reduction''},$$ 
  i.e. $V[G_{\d+1}]\models (Q,R)$ is not a code for a reduction. Now one of the items \eqref{QR1}--\eqref{QR9} fails for $(Q,R)$ in $V[G_{\d+1}]$. 
  We want to show that then 
  one of them fails in $V[G]$. The conditions \eqref{QR1}--\eqref{QR8} are absolute, so if one of them fails in $V[G_{\d+1}]$, then 
  we are done. Suppose \eqref{QR1}--\eqref{QR8} hold but \eqref{QR9} fails. Then  there is an 
  $\eta\in R^*$ such that $Q^*(\eta\cap S_{\{\a\}})\cap S_{\k\setminus\a}$ is stationary but $\eta\cap S_{\{\a\}}$ is not or
  vice versa. In $V[G_{\d+1}]$ define 
  $$\P^{\d+1}=\{(p_i)_{i<\k^+}\in \P_{\k^+}\mid (p_i)_{i<\d+1}\in G_{\d+1}\}.$$
  Then $\P^{\d+1}$ is $<\k$-closed. Thus it does not kill stationarity of any set. 
  So if $G^{\d+1}$ is $\P_{\d+1}$-generic over $V[G_{\d+1}]$,
  then in $V[G_{\d+1}][G^{\d+1}]$, $(Q,R)$ is not a code for a reduction. Now it remains to show that
  $V[G]=V[G_{\d+1}][G^{\d+1}]$ for some $G^{\d+1}$. In fact putting $G^{\d+1}=G$ we get $\P^{\d+1}$-generic over
  $V[G_{\d+1}]$ and of course $V[G_{\d+1}][G]=V[G]$ (since $G_{\d+1}\subset G$). 
\end{proofV}

\begin{Remark}
  The forcing constructed in the proof of Theorem \ref{thm:AntiChainChain} above, combined with the forcing
  in the proof of item \eqref{thm:CUBSETItem4} of Theorem \ref{thm:CUBSET} gives that for $\k^{<\k}=\k>\o_1$ not
  successor of a singular cardinal, we have in a forcing extension that
  $\la \Po(\k),\subset\ra$ embeds into $\la\E^{\Dii},\le_B\ra,$ i.e. the partial order of $\Dii$-equivalence relations
  under Borel reducibility.
\end{Remark}

\begin{Open}
  Can there be two equivalence relations, 
  $E_1$ and $E_2$ on $2^{\k}$, $\k>\o$ such that $E_1$ and $E_2$ are Borel and incomparable,
  i.e. $E_1\not\le_B E_2$ and $E_2\not\le_B E_1$? 
\end{Open}

\subsection{Reducibility Between Different Cofinalities}\label{ssec:RedDiffCof}

Recall the notation defined in Section \ref{sec:NotationsandC}.
In this section we will prove the following two theorems:

\begin{Thm}\label{thm:WeaklyCompact1}
  Suppose that $\k$ is a weakly compact cardinal and that $V=L$. Then
  \begin{myEnumerate}
  \item[(A)] $E_{S^\k_\l}\le_c E_{\reg(\k)}$ for any regular $\l<\k$, where $\reg(\k)=\{\l<\k\mid \l\text{ is regular}\}$,
  \item[(B)] In a forcing extension $E_{S^{\o_2}_\o}\le_c E_{S^{\o_2}_{\o_1}}.$ Similarly for $\l$, $\l^+$ and $\l^{++}$ instead of
    $\o$, $\o_1$ and $\o_2$ for any regular $\l<\k$.
  \end{myEnumerate}
\end{Thm}

\begin{Thm}\label{thm:DifferentCofinalitiesNoReduction}
  For a cardinal $\k$ which is a successor of a regular cardinal or $\k$ inaccessible, 
  there is a cofinality-preserving forcing
  extension in which for all regular $\l<\k$, the relations $E_{S^\k_\l}$ are $\le_B$-incomparable with each other.
\end{Thm}
Let us begin by proving the latter.
\begin{proofVOf}{Theorem \ref{thm:DifferentCofinalitiesNoReduction}}
  Let us show that there is a forcing extension of $L$ in which $E_{S^{\o_2}_{\o_1}}$
  and $E_{S^{\o_2}_{\o}}$ are incomparable. The general case is similar.

  We shall use Lemma \ref{lem:KillOneCode} with $\mu_1=\o$ and $\mu_2=\o_1$ and vice versa, and then a similar iteration
  as in the end of the proof of Theorem \ref{thm:AntiChainChain}. 
  First we force, like in the proof of Theorem \ref{thm:CUBSET} \eqref{thm:CUBSETItem4}, a
  stationary set $S\subset S^{\o_2}_{\o}$ such that for all $\a\in S^{\o_2}_{\o_1}$, $\a\cap S$ is non-stationary in $\a$. 
  Also for all $\a\in S^{\o_2}_\o$, $\a\cap S^{\o_2}_{\o_1}$ is non-stationary.
  
  By Lemma \ref{lem:KillOneCode}, for each code for a reduction from $E_{S}$ to $E_{S^{\o_2}_{\o_1}}$
  there is a $<\o_2$-closed $\o_3$-c.c. forcing which kills it. Similarly for each code for a reduction
  from $E_{S^{\o_2}_{\o_1}}$ to $E_{S^{\o_2}_{\o}}$.
  Making an $\o_3$-long iteration, similarly as in
  the end of the proof of Theorem \ref{thm:AntiChainChain}, we can kill all codes for reductions from
  $E_S$ to $E_{S^{\o_2}_{\o_1}}$ and from $E_{S^{\o_2}_{\o_1}}$ to $E_{S^{\o_2}_{\o}}$. 
  Thus, in the extension there are no reductions from $E_{S^{\o_2}_{\o_1}}$
  to $E_{S^{\o_2}_{\o}}$ and no reductions from $E_{S^{\o_2}_{\o}}$ to $E_{S^{\o_2}_{\o_1}}$. (Suppose there 
  is one of a latter kind, $f\colon 2^{\o_2}\to 2^{\o_2}$. Then
  $g(\eta)=f(\eta\cap S)$ is a reduction from $E_{S}$ to $E_{S^{\o_2}_{\o_1}}$.)
\end{proofVOf}

\begin{Def}\label{def:Diamond}
  Let $X,Y$ be subsets of $\k$ and suppose $Y$ consists of ordinals of uncountable cofinality. 
  We say that $X$ \emph{$\diamond$-reflects to} $Y$ if there exists
  a sequence $\la D_\a\ra_{\a\in Y}$ such that
  \begin{myEnumerate}
  \item $D_\a\subset \a$ is stationary in $\a$,
  \item if $Z\subset X$ is stationary, then $\{\a\in Y\mid D_\a=Z\cap \a\}$ is stationary.
  \end{myEnumerate}
\end{Def}

\begin{Thm}\label{thm:Diamond}
  If $X$ $\diamond$-reflects to $Y$, then $E_X\le_c E_Y$.
\end{Thm}
\begin{proof}
  Let $\la D_\alpha\ra_{\alpha \in Y}$ be the sequence of Definition \ref{def:Diamond}. For a
  set $A\subset\k$ define 
  $$f(A)=\{\alpha \in Y | A \cap X \cap D_\alpha\text{ is stationary in }\alpha\}.\eqno(i)$$ 
  We claim that $f$ is a continuous reduction. Clearly $f$ is continuous. Assume that 
  $(A\sd B)\cap X$ is non-stationary. Then there is 
  a cub set $C\subset \k\setminus [(A\sd B)\cap X]$. Now $A\cap X\cap C=B\cap X\cap C$ $(ii)$.  The set
  $C'=\{\a<\k\mid C\cap \a\text{ is unbounded in }\a\}$ is also cub and if $\a\in Y\cap C'$, we have
  that $D_\a\cap C$ is stationary in $\a$. Therefore for $\a\in Y\cap C'$ $(iii)$ we have the following equivalences:
  \begin{eqnarray*}
     \a\in f(A)&\iff&A\cap X\cap D_\a\text{ is stationary }\\
     &\stackrel{(iii)}{\iff}&A\cap X\cap C\cap D_\a\text{ is stationary}\\
     &\stackrel{(ii)}{\iff}&B\cap X\cap C\cap D_\a\text{ is stationary}\\
     &\stackrel{(iii)}{\iff}&B\cap X\cap D_\a\text{ is stationary}\\
     &\stackrel{(i)}{\iff}&\a\in f(B)
  \end{eqnarray*}
  Thus $(f(A)\sd f(B))\cap Y\subset \k\setminus C'$ and is non-stationary.

  Suppose $A\sd B$ is stationary. Then either $A\setminus B$ or $B\setminus A$ is stationary. Without loss of 
  generality suppose the former. Then 
  $$S=\{\a\in Y\mid (A\setminus B)\cap X\cap \a=D_\a\}$$
  is stationary by the definition of the sequence $\la D_\a\ra_{\a\in Y}$. Thus for $\alpha \in S$ we have that
  $A \cap X \cap D_\alpha = A \cap X \cap (A \setminus B) \cap X\cap \alpha = (A \setminus B) \cap X\cap \a$
  is stationary in $\alpha$ and 
  $B \cap X \cap D_\alpha = B \cap X \cap (A \setminus B) \cap X \cap \alpha=\es$
  is not stationary in $\a$. Therefore $(f(A)\sd f(B))\cap Y$ is stationary (as it contains~$S$).
\end{proof}

\begin{Fact}[$\Pii$-reflection]
  Assume that $\k$ is weakly compact.
  If $R$ is any binary predicate on
  $V_\k$ and $\forall A\f$ is some $\Pii$-sentence where $\f$ is a first-order sentence in 
  the language of set theory together with predicates $\{R,A\}$ such that
  $(V_\k,R)\models \forall A\f$, then there exists stationary many $\a<\k$ such that $(V_\a,R\cap V_\a)\models \forall A\f$.
\end{Fact}  

We say that $X$ \emph{strongly reflects to} $Y$ if for all stationary $Z\subset X$ there exist
stationary many $\a\in Y$ with $X\cap \a$ stationary in $\a$. 

\begin{Thm}\label{thm:StrongRefl}
  Suppose $V=L$, $\k$ is weakly compact and that 
  $X\subset \k$ and $Y\subset \reg\k$. If $X$ strongly reflects to $Y$, then $X$ $\diamond$-reflects to $Y$.
\end{Thm}
\begin{proof}
  Define $D_\a$ by induction on $\a\in Y$. For the purpose of the proof also define $C_\a$ for each $\a$ as follows.
  Suppose $(D_\b, C_\b)$ is defined for all $\b<\a$. Let 
  $(D,C)$ be the $L$-least\footnote{The least in the canonical definable ordering on $L$, see \cite{Kunen}.} pair such that
  \begin{myEnumerate}
    \item $C$ is cub subset of $\a$.
    \item $D$ is a stationary subset of $X\cap \a$
    \item for all $\b\in Y\cap C$, $D\cap \b\ne D_\b$
  \end{myEnumerate}
  If there is no such pair then set $D=C=\es$. 
  Then let $D_{\a}=D$ and $C_\a=C$. We claim that the sequence $\la D_\a\ra_{\a\in Y}$ is as needed.
  To show this, let us make a counter assumption: there is a stationary subset $Z$ of $X$ and a cub subset $C$ of $\k$
  such that
  $$C\cap Y\subset \{\a\in Y\mid D_\a\ne Z\cap \a\}.\eqno(\star)$$
  Let $(Z,C)$ be the $L$-least such pair.
  Let $\l>\k$ be regular and let $M$ be an elementary submodel of $L_\l$ such that
  \begin{myEnumerate}
  \item $|M|<\k$,
  \item $\a=M\cap \k\in Y\cap C$,
  \item $Z\cap \a$ is stationary in $\a$,
  \item $\{Z,C,X,Y,\k\}\subset M$
  \end{myEnumerate}
  (2) and (3) are possible by the definition of strong reflection.
  Let $\bar M$ be the Mostowski collapse of $M$ and let $G\colon M\to \bar M$
  be the Mostowski isomorphism. Then 
  $\bar M=L_{\g}$ for some $\g>\a$. Since 
  $\k\cap M=\a$, we have
  \begin{center}
    $G(Z)=Z\cap \a$, $G(C)=C\cap \a$, $G(X)=X\cap \a$, $G(Y)=Y\cap \a$ and $G(\k)=\a$, $(\star\star)$.
  \end{center}
  
  Note that by the definability of the canonical ordering of $L$, the sequence $\la D_{\b}\ra_{\b<\k}$
  is definable.  
  Let $\f(x,y,\a)$ be the formula which says
  \begin{center}
    ``$(x,y)$ is the $L$-least pair such that $x$ is contained in $X \cap \alpha$, $x$ 
      is stationary in $\alpha$, $y$ is cub in $\alpha$ and $x\cap \beta \neq D_\beta$ for all $\beta \in y \cap Y \cap \alpha$.''
  \end{center}
  By the assumption, 
  $$L\models \f(Z,C,\k)\text{, so }M\models \f(Z,C,\k)\text{ and }L_{\g}\models \f(G(Z),G(C),G(\k)).$$ 
  Let us show that this implies $L\models \f(G(Z),G(C),G(\k))$, i.e.
  $L\models \f(Z\cap\a,C\cap\a,\a)$.
  This will be a contradiction because then $D_\a=Z\cap\a$ which contradicts
  the assumptions (2) and $(\star)$ above.

  By the relative absoluteness of being the $L$-least, the relativised formula with parameters $\f^{L_\g}(G(Z),G(C),G(\k))$ says
  \begin{center}
    ``$(G(Z),G(C))$ is the $L$-least pair such that $G(Z)$ is contained in $G(X)$, $G(Z)$ 
      is $(\text{stationary})^{L_{\g}}$ 
      in $G(\k)$, $G(C)$ is cub in $G(\k)$ and $G(Z)\cap \beta \neq D_\beta^{L_\g}$ for all $\beta \in G(C) \cap G(Y) \cap G(\k)$.''
  \end{center}
  Written out this is equivalent to
  \begin{center}
     ``$(Z\cap\a,C\cap\a)$ is the $L$-least pair such that $Z\cap \a$ is contained in $X\cap \a$, $Z\cap \a$ 
      is $(\text{stationary})^{L_{\g}}$ in 
      $\a$, $C\cap\a$ is cub in $\a$ and $Z\cap \beta \neq D_\beta^{L_\g}$ for all $\beta \in C \cap Y \cap \a$.''
  \end{center}
  Note that this is true in $L$. Since $Z\cap \a$ is stationary in $\a$ also in $L$ by (3), 
  it remains to show by induction on $\b\in \a\cap Y$ that $Z\cap \a$ 
      $D_\b^{L_{\g}}=D_\b^L$ and $C_\b^{L_{\g}}=C_\b^L$  and we are done.
  Suppose we have proved this for $\d\in\b\cap Y$ and $\b\in \a\cap Y$. Then $(D_\b^{L_\g},C_\b^{L_\g})$ is 
  \begin{myAlphanumerate}
  \item  $(\text{the least }L\text{-pair})^{L_{\g}}$ such that 
  \item  $(C_\b\text{ is a cub subset of }\b)^{L_{\g}}$,
  \item  $(D_\b\text{ is a stationary subset of }\b)^{L_{\g}}$
  \item  and for all $\d\in Y\cap \b$, $(D_\b\cap\d \ne D_\d)^{L_{\g}}$.
  \item Or there is no such pair and $D_\b=\es$.
  \end{myAlphanumerate}
  The $L$-order is absolute as explained above, so (a) is equivalent
  to (the least $L$-pair)$^{L}$. Being a cub subset of $\a$ is also absolute
  for $L_{\g}$ so (b) is equivalent to  $(C_\b\text{ is a cub subset of }\a)^{L}$.
  All subsets of $\b$ in $L$ are elements of $L_{|\b|^+}$ (see~\cite{Kunen}),
  and since $\a$ is regular and $\b<\a\le\g$, we have $\Po(\b)\subset L_\g$.
  Thus 
  $$(D_\b\text{ is stationary subset of }\b)^{L_\g}\iff (D_\b\text{ is stationary subset of }\b)^{L}.$$
  Finally the statement of (d), 
  $(D_\b\cap\d \ne D_\d)^{L_{\g}}$
  is equivalent to $D_\b\cap\d \ne D_\d^{L_{\g}}$ as it is defining $D_\b$, but by the induction hypothesis
  $D_\d^{L_{\g}}=D_\d^L$, so we are done. For (e), the fact that 
  $$\Po(\b)\subset L_{|\b|^+}\subset L_{\a}\subset L_{\g}$$ 
  as above implies that if there is no such pair in $L_\g$, 
  then there is no such pair in~$L$.
\end{proof}

\begin{proof}[Proof of Theorem \ref{thm:WeaklyCompact1}]
  In the case (A) we will show that $S^\k_\l$ strongly reflects to $\reg(\k)$ in $L$ which suffices 
  by Theorems \ref{thm:Diamond} and \ref{thm:StrongRefl}. For (B) we will assume that $\k$ is a weakly compact cardinal 
  in $L$ and then collapse it to $\o_2$ to get a $\diamond$-sequence which witnesses that $S^{\o_2}_\o$ 
  $\diamond$-reflects to $S^{\o_2}_{\o_1}$ which is sufficient by Theorem \ref{thm:Diamond}.
  In the following we assume: $V=L$ and $\k$ is weakly compact.
    \vspace{5pt}
    
  \noindent {\it(A)}: Let us use $\Pii$-reflection. Let $X\subset S^{\k}_\l$. We want to show that the set 
    $$\{\l\in\reg(\k)\mid X\cap \l\text{ is stationary in }\l\}$$
    is stationary. Let $C\subset \k$ be cub.
    The sentence
    \begin{center}
      ``($X$ is stationary in $\kappa$) $\land$ ($C$ is cub in $\kappa$) $\land$ ($\kappa$ is regular)''      
    \end{center}
    is a $\Pii$-property of $(V_\kappa,X,C)$. By $\Pii$-reflection we get $\delta < \kappa$ such that 
    $(V_\delta,X \cap \delta, C \cap \delta)$ satisfies it. 
    But then $\delta$ is regular, $X \cap \delta$ is stationary and $\delta$ belongs to $C$.
    \vspace{5pt}

    \noindent{\it(B)}: Let $\k$ be weakly compact and let us Levy-collapse $\k$ to $\o_2$ with the following forcing:
     $$\P=\{f\colon \reg\k\to \k^{<\o_1}\mid \ran (f(\mu))\subset \mu,\ |\{\mu\mid f(\mu)\ne\es\}|\le\o\}.$$
     Order $\P$ by $f<g$ if and only if $f(\mu)\subset g(\mu)$ for all $\mu\in\reg(\k)$.
     For all $\mu$ put $\P_\mu=\{f\in\P\mid \sprt f\subset\mu\}$ and $\P^\mu=\{f\in\P\mid \sprt f\subset\k\setminus \mu\}$,
     where $\sprt$ means support, see Section \ref{ssec:Functions} on page~\pageref{ssec:Functions}.
     \begin{claim}{1} For all regular $\mu$, $\o<\mu\le \k$, $\P_\mu$ satisfies the following:
       \begin{myAlphanumerate}
       \item If $\mu>\o_1$, then $\P_\mu$ has the $\mu$-c.c.,
       \item $\P_\mu$ and $\P^\mu$ are $<\o_1$-closed,
       \item $\P=\P_\k\forces \omega_2=\check\k$,
       \item If $\mu<\k$, then $\P\forces\cf(\check\mu)=\o_1$, 
       \item if $p\in\P$, $\s$ a name and $p\forces$ ``$\s$ is cub in $\o_2$'', then there is cub $E\subset \k$
         such that $p\forces \check E \subset \s$.
       \end{myAlphanumerate}
     \end{claim}
     \vspace{-17pt}
     \begin{proof}
       Standard (see for instance \cite{Jech}).%
     \end{proof}
     We want to show that in the generic extension $S^{\o_2}_\o$ $\diamond$-reflects to $S^{\o_2}_{\o_1}$.
     It is sufficient to show that $S^{\o_2}_\o$ $\diamond$-reflects to some stationary $Y\subset S^{\o_2}_{\o_1}$
     by letting $D_\a=\a$ for $\a\notin Y$. In our case $Y=\{\mu\in V[G]\mid (\mu\in\reg(\k))^{V}\}$.
     By (d) of Claim 1, $Y\subset S^{\o_2}_{\o_1}$, $(\reg(\k))^V$ is stationary in $V$ (for instance by $\Pii$-reflection) and by (e)
     it remains stationary in $V[G]$. 

     It is easy to see that $\P\cong\P_\mu\times \P^{\mu}$. Let $G$ be a $\P$-generic over
     (the ground model) $V$. Define 
     $$G_\mu = G\cap \P_\mu.$$
     and 
     $$G^\mu = G\cap \P^\mu.$$
     Then $G_\mu$ is $\P_\mu$-generic over~$V$. 

     Also $G^\mu$ is $\P^\mu$-generic over $V[G_\mu]$ and $V[G] = V[G_\mu][G^\mu]$.

     Let 
     $$E=\{p\in \P\mid (p>q)\land (p_\mu\forces p^\mu\in \dot D)\}$$
     Then $E$ is dense above $q$: If $p>q$ is arbitrary element of $\P$, then
     $q\forces \exists p'>\check p^\mu(p'\in \dot D)$. 
     Thus there exists $q'>q$ with $q'>p_\mu$, $q'\in \P_\mu$ and $p'>p,p'\in \P^{\mu}$ such that
     $q'\forces p'\in \dot D$ and so $(q'\rest \mu)\cup (p'\rest(\k\setminus\mu))$
     is above $p$ and in $E$. So there is $p\in G\cap E$. But then
     $p_\mu\in G_\mu$ and $p^\mu\in G^\mu$ and $p_\mu\forces p^\mu\in \dot D$, 
     so $G^\mu\cap D\ne \es$. Since $D$ was arbitrary,
     this shows that $G^\mu$ is $\P^\mu$-generic over $V[G_\mu]$.
     Clearly $V[G]$ contains both $G_\mu$ and $G^\mu$. On the other hand,
     $G=G_\mu\cup G^\mu$, so $G\in V[G_\mu][G^\mu]$. By minimality of forcing extensions,
     we get $V[G]=V[G_\mu][G^\mu]$.

     For each $\mu\in \reg(\k)\setminus \{\o,\o_1\}$ let 
     $$k_\mu\colon \mu^+\to \{\sigma\mid \sigma\text{ is a nice }\P_\mu\text{ name for a subset of }\mu\}$$
     be a bijection. 
     A nice $\P_\mu$ name for a subset of $\check \mu$ is of the form
     $$\Cup\{\{\check\a\}\times A_{\a}\mid \a\in B\},$$
     where $B\subset \check\mu$ and for each $\a\in B$, $A_\a$ is an antichain in $\P_\mu$.
     By (a) there are no antichains of length $\mu$ in $\P_\mu$ and $|\P_\mu|=\mu$, so there are at most $\mu^{<\mu}=\mu$
     antichains and there are $\mu^+$ subsets $B\subset \mu$, so there indeed exists such a bijection 
     $k_\mu$ (these cardinality facts hold because $V=L$ and $\mu$ is regular). 
     Note that if $\sigma$ is a nice $\P_\mu$-name for a subset of $\check\mu$, then $\sigma\subset V_\mu$.
     
     Let us define 
     $$D_\mu=
     \begin{cases}
       \left[k_\mu\Big([(\cup G)(\mu^+)](0)\Big)\right]_G&\text{ if it is stationary}\\
       \mu &\text{ otherwise.}
     \end{cases}
     $$
     Now $D_\mu$ is defined for all $\mu\in Y$, recall 
     $Y=\{\mu\in V[G]\mid (\mu\in\reg\k)^{V}\}$. 
     We claim that $\la D_{\mu}\ra_{\mu\in Y}$ is the needed $\diamond$-sequence.
     Suppose it is not. Then there is a stationary set $S\subset S^{\o_2}_{\o}$
     and a cub $C\subset \o_2$ such that for all $\a\in C\cap Y$, 
     $D_\a\ne S\cap \a$. By (e) there is a cub set $C_0\subset C$ such that $C_0\in V$. Let $\dot S$
     be a nice name for $S$ and $p'$ such that $p'$ forces that $\dot S$ is stationary. Let us show that
     $$H=\{q\ge p'\mid q\forces D_\mu=\dot S\cap\check\mu\text{ for some }\mu\in C_0\}$$
     is dense above $p'$ which is obviously a contradiction. For that purpose let $p>p'$ be arbitrary and let us show that there
     is $q>p$ in $H$.
     Let us now use $\Pii$-reflection. 
     First let us redefine $\P$. Let $\P^*=\{q\mid \exists r\in \P(r\rest\sprt r=q)\}$.
     Clearly $\P^*\cong \P$ but the advantage is that $\P^*\subset V_\k$ and $\P_\mu^*=\P^*\cap V_\mu$
     where $\P^*_\mu$ is defined as $\P_\mu$.  
     One easily verifies that all the above things (concerning $\P_\mu$, $\P^\mu$ etc.) translate
     between $\P$ and $\P^*$. From now on denote
     $\P^*$ by $\P$.
     Let 
     $$R=(\P\times\{0\})\cup (\dot S\times \{1\})\cup (C_0\times\{2\})\cup (\{p\}\times \{3\})$$
     Then $(V_\k,R)\models \forall A\f$, where $\f$ says:
     ``(if $A$ is closed unbounded and $r>p$ arbitrary, then there exist $q>r$ and $\a$ 
     such that $\a\in A$
     and $q\forces_{\P} \check \a\in \dot S$).''
     So basically $\forall A\f$ says ``$p\forces$ ($\dot S$ is stationary)''. It follows from (e)
     that it is enough to quantify over cub sets in $V$.
     Let us explain why such a formula can be written for $(V_\k,R)$. The sets (classes from the
     viewpoint of $V_\k$) $\P$, $\dot S$ and $C_0$ are coded into $R$, so we can use them as parameters.
     That $r>p$ and $q>r$ and $A$ is closed and unbounded is expressible in first-order as well as $\a\in A$. 
     How do we express $q\forces_{\P} \check \a\in \dot S$?
     The definition of $\check\a$ is recursive in $\a$: 
     $$\check\a=\{(\check\b,1_{\P})\mid \b<\a\}$$
     and is absolute for $V_\k$. Then $q\forces_\P \check\a\in \dot S$ is 
     equivalent to saying that for each $q'>q$ there exists $q''>q'$ with
     $(\check\a,q'')\in \dot S$ and this is expressible in first-order (as we have taken $R$ as a parameter).
     
     By $\Pii$-reflection
     there is $\mu\in C_0$ such that $p\in \P_\mu$ and $(V_{\mu},R)\models \forall A\f$. 
     Note that we may require that $\mu$
     is regular, i.e. $(\check\mu_G\in Y)^{V[G]}$ and such that $\a\in S\cap \mu$ implies
     $(\check\a,\check p)\in \dot S$ for some $p\in \P_\mu$. Let $\dot S_\mu=\dot S\cap V_\mu$.

     Thus
     $p\forces_{\P_\mu}$ ``$\dot S_\mu$ is stationary''.
     Define $q$ as follows: $\dom q=\dom p\cup\{\mu^+\}$, $q\rest\mu=p\rest\mu$ and
     $q(\mu^+)=f$, $\dom f=\{0\}$ and $f(0)=k_{\mu}^{-1}(\dot S_\mu)$.
     Then $q\forces_{\P} \dot S_\mu=D_\mu$ provided that $q\forces_\P$ ``$\dot S_\mu$ is stationary''.
     The latter holds since $\P^\mu$ is $<\o_1$-closed., and does not kill stationarity of $(\dot S_\mu)_{G_\mu}$ so 
     $(\dot S_\mu)_{G_\mu}$ is stationary in $V[G]$ and by the assumption on $\mu$, $(\dot S_\mu)_{G_\mu}=(\dot S_\mu)_{G}$. 
     Finally, it remains to show that in $V[G]$, $(\dot S_\mu)_{G}=S\cap \mu$. But this again follows from the definition~of~$\mu$.
     
     Instead of collapsing $\k$ to $\omega_2$, we could do the same for $\l^{++}$ for any regular $\l<\k$
     and obtain a model in which $E_{S^{\l^{++}}_{\l}}\le_c E_{S^{\l^{++}}_{\l^+}}$.
     \qedhere
 \end{proof}

 \begin{Open}
    Is it consistent that $S^{\omega_2}_{\omega_1}$ Borel reduces to $S^{\omega_2}_\omega$?
 \end{Open}

\subsection{$E_0$ and $E_{S^\k_\l}$}

In the Section \ref{ssec:RedDiffCof} 
above, Theorem \ref{thm:DifferentCofinalitiesNoReduction}, 
we showed that the equivalence relations
of the form $E_{S^\k_\l}$ can form an antichain with respect to $\le_B$.
We will show that under mild set theoretical assumptions, all of them are strictly above
$$E_0=\{(\eta,\xi)\mid \eta^{-1}\{1\}\sd\xi^{-1}\{1\}\text{ is bounded}\}.$$

\begin{Thm}\label{thm:EzeroRedToNonStat}
  Let $\k$ be regular and $S\subset \k$ stationary 
  and suppose that $\diamondsuit_\k(S)$ holds (i.e., 
  $\diamondsuit_\kappa$ holds on the stationary set 
  $S$). Then $E_0$ is Borel 
  reducible to~$E_{S}$. 
\end{Thm}
\begin{proof}
  The proof uses similar ideas than the proof of Theorem \ref{thm:Diamond}.
  Suppose that the $\diamondsuit_\k(S)$ holds and let 
  $\la D_\a\ra_{\a\in S}$ be the  $\diamondsuit_\k(S)$-sequence.
  Define the reduction $f\colon 2^\k\to 2^\k$ by
  $$f(X) =\{\a \in S \mid D_\a\text{ and }X \cap \a\text{ agree on a final segment of }\a\}$$
  If $X,Y$ are $E_0$-equivalent, then $f(X)$, $f(Y)$ are $E_S$-equivalent, because they are in fact even $E_0$-equivalent as is easy to check. 
  If $X,Y$ are not $E_0$-equivalent, then there is a club $C$ of $\a$ where $X$, $Y$ differ cofinally in 
  $\a$; it follows that $f(X)$, $f(Y)$ differ on a stationary subset of $S$, 
  namely the elements $\a$ of $C\cap S$ where $D_\a$ equals $X \cap \a$.
\end{proof}

\begin{Cor}\label{cor:EzeroNotRed} 
  Suppose $\k=\l^+=2^{\l}$. Then $E_0$ is Borel reducible to $E_{S}$ where $S\subset \k\setminus S^\k_{\cf(\l)}$ is stationary.
\end{Cor}
\begin{proof}
  Gregory proved in \cite{Greg} that if $2^\mu=\mu^+=\k$, $\mu$ is regular and $\l<\mu$, then 
  $\diamondsuit_\k(S^\k_\l)$ holds. 
  Shelah extended this result in \cite{Shelah6}
  and proved that if $\k=\l^+=2^\l$
  and $S\subset \k\setminus S^\k_{\cf(\l)}$, then
  $\diamondsuit_\k(S)$ holds. Now apply Theorem~\ref{thm:EzeroRedToNonStat}.
\end{proof}

\begin{Cor}[GCH] 
  Let us assume that $\k$ is a successor cardinal.
  Then in a cofinality 
  and GCH preserving forcing extension,
  there is an embedding
  $$f\colon\la\Po(\k),\subset\ra\to \la\E^{\Sii},\le_B\ra,$$
  where $\E^{\Sii}$ is the set of $\Sii$-equivalence relations
  (see Theorem \ref{thm:AntiChainChain}) such that for all $A\in \Po(\k)$, $E_0$ is strictly below~$f(A)$.
  If $\k$ is not the successor of an $\o$-cofinal cardinal, we may replace $\Sii$ above by Borel*.
\end{Cor}
\begin{proof}
  Suppose first that $\k$ is not the successor of an $\o$-cofinal cardinal.
  By Theorem \ref{thm:AntiChainChain} there is a GCH and cofinality-preserving forcing 
  extension 
  such that there is an embedding
  $$f\colon\la\Po(\k),\subset\ra\to \la\E^{\text{Borel}^*},\le_B\ra.$$
  From the proof of Theorem \ref{thm:AntiChainChain} one sees that $f(A)$ is of the form $E_{S}$ where
  $S\subset S^\k_\o$. Now $E_0$ is reducible to such relations by Corollary \ref{cor:EzeroNotRed}, as GCH continues to hold in the extension. 

  So it suffices to show that $E_{S}\not\le_B E_0$ for stationary $S\subset S^\k_\o$. By the same argument as in Corollary \ref{cor:NSIsNotBorel}, 
  $E_{S}$ is not Borel and by Theorem \ref{thm:EzeroTheorem} $E_0$ is Borel, 
  so by Fact \ref{fact:ReductionNote} $E_{S^\k_\l}$ is not
  reducible to~$E_0$.

  Suppose $\k$ is the successor of an $\o$-cofinal ordinal and $\k>\o_1$. Then, in the proof of Theorem \ref{thm:AntiChainChain} replace $\mu$ by $\o_1$
  and get the same result as above but for relations of the form $E_S$ where $S\subset S^\k_{\o_1}$. 

  The remaining case is $\k=\o_1$. Let $\{S_{\a}\mid \a<\o_1\}$ be a set of pairwise disjoint
  stationary subsets of $\o_1$. Let $\P$ be the forcing given by the proof of Theorem \ref{thm:AntiChainChain}
  such that in the $\P$-generic extension the function $f\colon \la\Po(\o_1),\subset\ra \to \la\E^{\text{Borel}^*},\le_B\ra$
  given by $f(A)=E_{\Cup_{\a\in A} S_\a}$ is an embedding. This forcing preserves stationary sets, so 
  as in the proof of clause \eqref{thm:CUBSETItem4} of Theorem \ref{thm:CUBSET}, we can first force a 
  $\diamondsuit$-sequence which guesses each subset of $\Cup_{\a<\o_1}S_{\a}$ on a set $S$ such that
  $S\cap S_\a$ is stationary for all $\a$. Then by Corollary \ref{cor:EzeroNotRed} $E_0$ is reducible
  to $E_{\Cup_{\a\in A}S_\a}$ for all $A\subset \k$.
\end{proof}

\chapter{Complexity of Isomorphism Relations}\label{chapter:ComplexityofIsomRel}
Let $T$ be a countable complete theory. Let us turn to the question discussed in Section 
\ref{sec:HistMotiv}: ``How is the set theoretic complexity of $\cong_T$ related 
to the stability theoretic properties of $T$?''. The following theorems give some answers. 
As pointed out in Section \ref{sec:HistMotiv}, the assumption that $\k$ is uncountable is crucial in the following theorems. 
For instance the theory of dense linear orderings without end points is unstable, but
$\cong_T$ is an open set in case $\k=\o$, while we show below that for unstable theories
$T$ the set $\cong_T$ cannot be even $\Dii$ when $\k>\o$. Another example introduced by Martin Koerwien in his Ph.D. thesis 
and in \cite{Koe} shows that
there are classifiable shallow theories whose isomorphism is not Borel when $\k=\o$, although we prove below
that the isomorphism of such theories is always Borel, when $\k>\o$.
This justifies in particular the motivation
for studying the space $\k^\k$ for model theoretic purposes: the set theoretic complexity of $\cong_T$
positively correlates with the model theoretic complexity of $T$. 

The following stability theoretical notions will be used:
stable, superstable, DOP, OTOP, shallow, $\l(T)$ and $\k(T)$. Classifiable means superstable with no DOP nor OTOP
and $\l(T)$ is the least cardinal in which $T$ is stable.

The main theme in this section is exposed in the following two theorems:

\begin{Thm}[$\k^{<\k}=\k>\o$]\label{thm:ShallowBorell}
  Assume that $\k$ is not weakly inaccessible and $T$ a complete countable 
  first-order theory. 
  If the isomorphism relation $\cong^\k_T$ is Borel, then $T$ is classifiable and shallow.
  Conversly, if $\k>2^\o$, then if $T$ is classifiable and shallow, then
  $\cong^\k_T$ is Borel.
\end{Thm}

\begin{Thm}[$\k^{<\k}=\k$]\label{thm:DiiisClass}
  Assume that for all $\l<\k$, $\l^\o<\k$ and $\k>\o_1$. 
  Then in $L$ and 
  in the forcing extension after
  adding $\k^+$ Cohen subsets of $\k$ we have: for any theory $T$, 
  $T$ is classifiable if and only if $\cong_T$ is~$\Dii$.
\end{Thm}
The two theorems above are proved in many subtheorems below. Our results are
stronger than those given by \ref{thm:ShallowBorell} and \ref{thm:DiiisClass} (for instance the cardinality
assumption $\k>\o_1$ is needed only in the case where $T$ is superstable with DOP and the stable unsuperstable
case is the only one for which Theorem \ref{thm:DiiisClass} cannot be proved in ZFC). Theorem \ref{thm:ShallowBorell}
follows from Theorems \ref{thm:ClasShalIsBorel}, \ref{thm:ClasNotShalIsNotBorel}. 
Theorem \ref{thm:DiiisClass} follows from Theorems \ref{thm:Dii},
\ref{thm:NotDiiOrNotBorelList}, \ref{thm:StabUnsupstab} 
and items \eqref{thm:CUBSETItem2} and \eqref{thm:CUBSETItem7} of Theorem~\ref{thm:CUBSET}.

\section{Preliminary Results}
The following Theorems \ref{thm:BorelTree} and \ref{thm:NotD1} 
will serve as bridges between the set theoretic complexity and the model theoretic complexity
of an isomorphism relation.

\begin{Thm}[$\k^{<\k}=\k$]\label{thm:BorelTree}
  For a theory $T$, the set $\cong_T$ is Borel if and only if the following holds:
  there exists a $\k^+\o$-tree $t$ such that for all models $\A$ and $\B$ of $T$,
  $\A\cong \B \iff \PlTwo\uparrow \EF^\k_t(\A,\B)$. 
\end{Thm}
\begin{proofV}{Theorem \ref{thm:BorelTree}}
  Recall that we assume $\dom\A=\k$ for all models in the discourse.
  First suppose that there exists a $\k^+\o$-tree $t$ such that for all models $\A$ and $\B$ of $T$,
  $\A\cong \B \iff \PlTwo\uparrow \EF^\k_t(\A,\B)$. 
  Let us show that there exists a $\k^+\o$-tree $u$ which constitutes a Borel code for $\cong_T$ (see Remark \ref{remark:Borel} 
  on page~\pageref{remark:Borel}).

  Let $u$ be the tree of sequences of the form
  $$\la (p_0,A_0), f_0, (p_1,A_1), f_1, \dots,(p_n,A_n),f_n\ra$$
  such that for all $i\le n$
  \begin{myEnumerate}
  \item $(p_i, A_i)$ is a move of player $\PlOne$ in $\EF^\k_t$, i.e. $p_i\in t$ and $A_i\subset \k$ with $|A_i|<\k$,
  \item $f_i$ is a move of player $\PlTwo$ in $\EF^\k_t$, i.e. it is a partial function $\k\to\k$ with $|\dom f_i|,|\ran f_i|<\k$
    and $A_i\subset \dom f_i\cap \ran f_i$
  \item $\la (p_0,A_0), f_0, (p_1,A_1), f_1, \dots,(p_n,A_n),f_n\ra$ is a valid position of the game, 
    i.e. $(p_i)_{i\le n}$ is an initial segment of a branch in 
    $t$ and $A_i\subset A_j$ and $f_i\subset f_j$ whenever $i<j\le n$.
  \end{myEnumerate}
  Order $u$ by end extension.
  The tree $u$ is a $\k^+\o$-tree 
  (because $t$ is and by~(3)).
  
  Let us now define the function
  $$h\colon\{\text{branches of }u\}\to\{\text{basic open sets of }(\k^\k)^2\}.$$
  Let $b\subset u$ be a branch, 
  $$b=\{\es,\la (p_0,A_0)\ra,\la (p_0,A_0),f_0\ra,\dots, \la (p_0,A_0),f_0,\dots,(p_k,A_k),f_k\ra\}.$$
  It corresponds to a unique $\EF$-game between some two structures
  with domains $\k$. In this game the players have chosen some set $A_k=\Cup_{i\le k}A_i\subset\k$
  and some partial function $f_k=\Cup_{i\le k}f_i\colon \k\to\k$. 
  Let $h(b)$ be the set of all pairs $(\eta,\xi)\in (\k^\k)^2$ such that $f_\k\colon \A_{\eta}\rest A_\k \cong \A_{\xi}\rest A_\k$
  is a partial isomorphism. This is clearly an open set: 
  $$(\eta,\xi)\in h(b)\Rightarrow N_{\eta\restl((\sup A_\k)+1)}\times N_{\xi\restl ((\sup A_\k) +1)}\subset h(b).$$
  
  Finally we claim that $\A_\eta\cong\A_\xi \iff \PlTwo\uparrow G(u,h,(\eta,\xi))$. Here
  $G$ is the game as in Definition \ref{def:Borel} of Borel* sets, page \pageref{def:Borel} but played on the
  product $\k^{\k}\times\k^\k$.
  Assume $\A_\eta\cong\A_{\xi}.$ Then $\PlTwo\uparrow \EF^\k_t(\A_\eta,\A_\xi)$. Let $\upsilon$ denote the winning strategy.
  In the game $G(u,h,(\eta,\xi))$, let us define a winning strategy for player $\PlTwo$ as follows.
  By definition, at a particular move, say $n$, $\PlOne$ chooses a sequence 
  $$\la(p_0,A_0),f_0,\dots (p_n,A_n)\ra.$$
  Next $\PlTwo$
  extends it according to $\upsilon$ to 
  $$\la(p_0,A_0),f_0,\dots (p_n,A_n),f_n\ra,$$
  where $f_n=\upsilon((p_0,A_0),\dots, (p_{n},A_n))$.
  Since $\upsilon$ was a winning strategy, it is clear that $f_\k=\Cup_{i<\k} f_i$ is going to be a isomorphism between 
  $\A_\eta\rest A_\k$ and $\A_\xi\rest A_\k$, so $(\eta,\xi)\in h(b)$.
  
  Assume that $\A_\eta\not\cong\A_\xi$. Then by the assumption there is no winning strategy of 
  $\PlTwo$, so player $\PlOne$ can play in such a way 
  that $f_\k=\Cup_{i\le \k} f_i$ is not an isomorphism between $\A_\eta\rest \cup A_i$ and $\A_\xi\rest\cup A_i$, 
  so $(\eta,\xi)$ is not in $h(b)$. This completes the proof of the direction ``$\Leftarrow$''

  Let us prove ``$\Rightarrow$''.
  Suppose $\cong_T$ is Borel and let us show that there is a tree as in the statement of the theorem. 
  We want to use Theorem \ref{thm:BorelIsLkk} and formalize the statement ``$\cong_T$ is definable in $L_{\k^+\k}$''
  by considering the space consisting of pairs of models.
  
  Denote the vocabulary of $\A$ and $\B$ as usual by $L$. 
  Let $P$ be a unary relation symbol not in $L$. We will
  now discuss two distinct vocabularies, $L$ and $L\cup \{P\}$ at the same time, so we have to introduce two distinct codings.
  Fix an $\eta\in 2^{\k}$. Let $\A_\eta$ denote the $L$-structure as defined in
  Definition \ref{def:CodingOfModels} of our usual coding.
  Let $\rho\colon \k\cup\k^{<\o}\to \k$ be a bijection and define $\A^\eta$ to be the model with $\dom\A^\eta=\k$ and
  if $a\in \dom\A^\eta$, then $\A^\eta\models P(a)\iff \eta(\rho(a))=1$ such that if $(a_1,\dots,a_n)\in (\dom\A^\eta)^{n}$,
  then $\A^\eta\models P_n(a_1,\dots,a_n)\iff \eta(\rho(a_1,\dots,a_n))=1$. Note that we are making a distinction here between
  $\k$ and $\k^{\{0\}}$. 

  \begin{claim}{1}
    The set $W=\{\eta\in 2^\k\mid \k=|P^{\A^\eta}|=|\k\setminus P^{\A^\eta}|\}$ is Borel.
  \end{claim}
  \begin{proofVOf}{Claim 1}
    Let us show that the complement is Borel. By symmetry it is sufficient to show that
    $$B=\{\eta\mid \k>|P^{\A^\eta}|\}$$ is Borel.
    Let $I\subset \k$ be a subset of size $<\k$. For $\b\notin I$ define $U(I,\b)$ to be the set
    $$U(I,\b)=\{\eta\mid \eta(\rho(\b))=0\}.$$
    Clearly $U(I,\b)$ is open for all $I$, $\b$.
    Now 
    $$B=\Cup_{I\in [\k]^{<\k}}\Cap_{\b\notin I}U(I,\b).$$
    By the assumption $\k^{<\k}=\k$, this is Borel (in fact a union of closed sets).
  \end{proofVOf}

  Define a mapping $h\colon W\to (2^{\k})^{2}$ as follows. Suppose $\xi\in W$. Let 
  $$r_1\colon \k\to P^{\A^\xi}$$ 
  and  
  $$r_2\colon \k\to \k\setminus P^{\A^\xi}$$
  be the order preserving bijections (note $P^{\A^\eta}\subset\k=\dom \A^\eta$).
  
  Let $\eta_1$ be such that $r_1$ is an isomorphism 
  $$\A_{\eta_1}\to (\A^{\xi}\cap P^{\A^\xi})\rest L$$
  and 
  $\eta_2$ such that $r_2$ is an isomorphism 
  $$\A_{\eta_2}\to (\A^{\xi}\setminus P^{\A^\xi})\rest L.$$
  Clearly $\eta_1$ and $\eta_2$ are unique, so we can define $h(\xi)=(\eta_1,\eta_2)$.

  \begin{claim}{2}
    $h$ is continuous.
  \end{claim}
  \begin{proofVOf}{Claim 2}
     Let $U=N_{p}\times N_q$ be a basic open set of $(2^{\k})^2$, $p,q\in 2^{<\k}$ and let $\xi\in h^{-1}[U]$.
     Let $P^{\A^\xi}=\{\b_i\mid i<\k\}$ be an enumeration such that $\b_i<\b_j\iff i<j$ and similarly
     $\k\setminus P^{\A^\xi}=\{\g_i\mid i<\k\}$. Let $\a=\max\{\b_{\dom p},\g_{\dom q}\}+1$. Then
     $N_{\xi\restl\a}\subset h^{-1}[U]$. Thus arbitrary $\xi$ in $h^{-1}[U]$ have an open neighbourhood
     in $h^{-1}[U]$, so it is open.
  \end{proofVOf}

  Recall our assumption that $E=\{(\eta,\xi)\in 2^{\k}\mid \A_\eta\cong\A_\xi\}$ is Borel. Since
  $h$ is continuous and in particular Borel, this implies that
  $$E'=\{\eta\mid \A_{h_1(\eta)}\cong \A_{h_2(\eta)}\}=h^{-1}E$$
  is Borel in $W$. Because $W$ is itself Borel, $E'$ is Borel in $2^\k$. Additionally, $E'$ is closed under permutations:
  if $\A^\eta$ is isomorphic to $\A^{\xi}$, then $\A^\eta\cap P^{\A^\eta}$ is isomorphic to $\A^\xi\cap P^{\A^\xi}$
  and $\A^\eta\setminus P^{\A^\eta}$ is isomorphic to $\A^\xi\setminus P^{\A^\xi}$, so if $\A^\eta\in E'$, then also 
  $\A^\xi\in E'$ (and note that since $\eta\in W$, also $\xi\in W$). By Theorem \ref{thm:BorelIsLkk}, there is a sentence $\theta$ of 
  $L_{\k^+\k}$ over $L\cup\{P\}$ that defines $E'$. 
  Thus by Theorem \ref{thm:QuantifierRankAndScottRankTheorem} and Remark \ref{rem:LkkLinfk} there is a $\k^+\o$-tree 
  $t$ such that 
  $$\text{ if }\eta\in E'\text{ and }\xi\notin E'\text{, then }\PlTwo\not\wins\EF^\k_t(\A^\eta,\A^\xi).\eqno\bigodot$$ 
  We claim that $t$ is as needed, i.e. for all models $\A,\B$ of $T$
  $$\A\cong \B\iff \PlTwo\wins \EF^\k_t(\A,\B).$$
  Suppose not. Then there are models $\A\not\cong \B$ such that $\PlTwo\wins\EF^\k_t(\A,\B)$. 
  Let $\eta$ and $\xi$ be such that $\A_{h_1(\eta)}=\A_{h_2(\eta)}=\A_{h_1(\xi)}=\A$ and $\A_{h_2(\xi)}=\B$.
  Clearly $\eta\in E'$, but $\xi\notin E'$, so by $\bigodot$ there is no winning strategy of 
  $\PlTwo$ in $\EF^{\k}_t(\A^{\eta},\A^{\xi})$ which is clearly a contradiction, because $\PlTwo$ can apply her
  winning strategies in $\EF^\k_t(\A,\B)$ and $\EF^\k_t(\A,\A)$ to win in $\EF^\k_t(\A^\eta,\A^\xi)$.
\end{proofV}

We will use the following lemma from \cite{MekVaa}:

\begin{Lemma}\label{lem:Diii}
  If $t\subset (\k^{<\k})^2$ is a tree and $\xi\in \k^\k$, denote
  $$t(\xi)=\{p\in \k^{<\k}\mid (p,\xi\rest\dom p)\in t\}$$
  Similarly if $t\in (\k^{<\k})^3$, then
  $$t(\eta,\xi)=\{p\in \k^{<\k}\mid (p,\eta\rest\dom p,\xi\rest\dom p)\in t\}.$$
  Assume that $Z$ is $\Sigma_1^1$. Then $Z$ is $\Delta_1^1$ if and only if for every tree 
  $t\subset (\k^{<\k})^{2}$ such that 
  $$t(\xi)\text{ has a }\k\text{-branch}\iff \xi\in Z$$ 
  there
  exists a $\k^+\k$-tree $t'$ such that $\xi\in Z\iff t(\xi)\not\le t'$.
  (Recall that $t\le t'$ when there exists a strictly order preserving map $t\to t'$)
\end{Lemma}


\begin{Thm}\label{thm:NotD1}
  Let $T$ be a theory and assume that for every
  $\k^+\k$-tree $t$ there exist $(\eta,\xi)\in (2^\k)^2$ such that $\A_\eta,\A_\xi\models T$,
  $\A_\eta\not\cong\A_\xi$ but $\PlTwo\uparrow\EF^\k_t(\A_\eta,\A_\xi)$.
  Then $\cong_T$ is not $\Dii$.
\end{Thm}
\begin{proofV}{Theorem \ref{thm:NotD1}}
  Let us abbreviate some statements:
  \begin{myItemize}
  \item[$A(t)$:] $t\subset (\k^{<\k})^3$ is a tree and for all  $(\eta,\xi)\in (\k^\k)^2$, 
    $$(\eta,\xi)\in \,\cong_T\iff t(\eta,\xi)\text{ contains a }\k\text{-branch }.$$
  \item[$B(t,t')$:] $t\subset (\k^{<\k})^3$ is a $\k^+\k$-tree and for all $(\eta,\xi)\in \k^\k$,
    $$(\eta,\xi)\in \,\cong_T\iff t(\eta,\xi)\not\le t'.$$
  \end{myItemize}
  Now Lemma \ref{lem:Diii} implies that if $\cong_T$ is $\Dii$, then 
  $\forall t [A(t)\rightarrow \exists t'B(t,t')]$. We will show that
  $\exists t [A(t)\land \forall t'\lnot B(t,t')]$, which by Lemma \ref{lem:Diii} suffices to prove the theorem.
  Let us define $t$. In the following, $\nu_\a$, $\eta_\a$ and $\xi_\a$ stand respectively for $\nu\rest\a$,
  $\eta\rest\a$ and $\xi\rest\a$.
  $$t=\{(\nu_\a,\eta_\a,\xi_\a)\mid \a<\k\text{ and }\nu\text{ codes an isomorphism between } 
  \A_\eta \text{ and }\A_\xi\}.$$
  Using Theorem \ref{thm:CodedClosed} it is easy to see that $t$ satisfies $A(t)$.
  Assume now that $t'$ is an arbitrary $\k^+\k$-tree. We will show that $B(t,t')$ does not hold.
  For that purpose
  let $u=\o\times t'$ be the tree defined by the set $\{(n,s)\mid n\in\o,s\in t'\}$
  and the ordering 
  $$(n_0,s_0)<_u(n_1,s_1)\iff \big(s_0<_{t'} s_1 \lor (s_0=s_1\land n_0<_\o n_1)\big).\eqno(1)$$
  This tree $u$ is still a $\k^+\k$-tree, so by the assumption
  of the theorem there is a pair $(\xi_1,\xi_2)$ such that $\A_{\xi_1}$ and $\A_{\xi_2}$ 
  are non-isomorphic, but $\PlTwo\uparrow\EF^\k_{u}(\A_{\xi_1},\A_{\xi_2})$.
  
  It is now sufficient to show that $t(\xi_1,\xi_2)\not\le t'$.
  \vspace{5pt}
  
  \begin{claim}{1}
    There is no order preserving function
    $$\sigma t' \to t',$$ 
    where $\s t'$ is defined in Definition \ref{def:sigmaOfTree}.   
  \end{claim}
  \begin{proofVOf}{Claim 1}
    Assume $g\colon \sigma t' \to t',$ is order preserving. Define $x_0 = g(\es)$ and
    $$x_\a = g(\{y\in t'\mid \exists \b<\a(y\le x_\b)\})\text{ for }0<\a<\k$$
    Then $(x_\a)_{\a<\k}$ contradicts the assumption that $t'$ is a $\k^+\k$-tree.
  \end{proofVOf}
  \begin{claim}{2}
    There is an order preserving function
    $$\sigma t'\to t(\xi_1,\xi_2).$$
  \end{claim}
  \begin{proofVOf}{Claim 2}
    The idea is that players $\PlOne$ and $\PlTwo$ play an $\EF$-game for each
    branch of the tree $t'$ and $\PlTwo$ uses her winning strategy in $\EF^\k_u(\A_{\xi_1},\A_{\xi_2})$to embed that branch into
    the tree of partial isomorphisms. A problem is that the winning strategy gives arbitrary
    partial isomorphisms while we are interested in those which are coded by functions defined
    on page \pageref{page:Closed}. Now the tree $u$ of $(1)$ above becomes useful.
    
    Let $\sigma$ be a winning strategy of player $\PlTwo$ in $\EF^\k_{u}(\A_{\xi_1},\A_{\xi_2})$. 
    Let us define $g\colon \sigma t'\to t(\xi_1,\xi_2)$ recursively. 
    Recall the function $\pi$ from Definition \ref{def:CodingOfModels} and define
    $$C=\{\a\mid \pi[\a^{<\o}]=\a\}.$$
    Clearly $C$ is cub.
    If $s\subset t'$ is an element of $\sigma t'$,
    then we assume
    that $g$ is defined for all $s'<_{\sigma t'} s$ and that $\EF^\k_u$ is played
    up to $(0,\sup s)\in u$. If $s$ does not contain its supremum, then put 
    $g(s)=\Cup_{s'<s} g(s')$. Otherwise let them continue playing the game for $\o$ more moves; 
    at the $n\th$ of these moves player $\PlOne$
    picks $(n,\sup s)$ from $u$ and a $\b<\k$ where $\b$ is an element of $C$ above
    $$\max\{\ran f_{n-1},\dom f_{n-1}\}$$
    where $f_{n-1}$ is the previous move by $\PlTwo$. (If $n=0$, it does not matter what $\PlOne$ does.)
    In that way the function $f=\Cup_{n<\o}f_n$ is a partial isomorphism such that $\dom f=\ran f=\a$ for some ordinal $\a$. 
    It is straightforward to check that such an $f$ is coded by some $\nu_\a\colon \a\to \k$. 
    It is an isomorphism between $\A_{\xi_1}\cap \a$ and $\A_{\xi_2}\cap \a$ and since $\a$ is in $C$,
    there are $\xi_1'$ and $\xi_2'$ such that $\xi_1\rest\a\subset \xi_1'$, $\xi_2\rest\a\subset\xi_2'$ and
    there is an isomorphism $\A_{\xi_1'}\cong \A_{\xi_2'}$ coded by some $\nu$ such that $\nu_\a=\nu\rest\a$.
    Thus $\nu_\a\in t(\xi_1,\xi_2)$ is suitable for setting $g(s)=\nu_\a$.
  \end{proofVOf}
\end{proofV}
\vspace{-15pt}

\section{Classifiable}

Throughout this section $\k$ is a regular cardinal satisfying $\k^{<\k}=\k>\o$.

\begin{Thm}[$\k^{<\k}=\k>2^\o$]\label{thm:ClasShalIsBorel}
  If the theory $T$ is classifiable and shallow, then $\cong_T$ is Borel.
\end{Thm}
\begin{proof} 
  If $T$ is classifiable and shallow, then from \cite{Shelah2} Theorem XIII.1.5 it follows
  that the models of $T$ are characterized by the game $\EF^\k_t$ up to isomorphism, where $t$ is some $\k^+\o$-tree
  (in fact a tree of descending sequences of an ordinal $\a<\k^+$). 
  Hence by Theorem \ref{thm:BorelTree} the isomorphism relation of $T$ is Borel.
\end{proof}

\begin{Thm}\label{thm:ClasNotShalIsNotBorel}
  If the theory $T$ is classifiable but not shallow, then $\cong_T$ is not Borel.
  If $\k$ is not weakly inaccessible and $T$ is not classifiable, then $\cong_T$ is not Borel.
\end{Thm}
\begin{proof}
  If $T$ is classifiable but not shallow, then by \cite{Shelah2} XIII.1.8,
  the $L_{\infty\k}$-Scott heights of models of $T$ of size $\k$ are not bounded by any ordinal $<\k^+$ 
  (see Definition \ref{def:ScottHeight} on page \pageref{def:ScottHeight}).
  Because any $\k^+\o$-tree can be embedded into $t_\a=\{\text{decreasing sequences of }\a\}$ for some $\a$ 
  (see Fact \ref{thm:Embedkappapluskappatree} on page~\pageref{thm:Embedkappapluskappatree}),
  this implies that for any $\k^+\o$-tree $t$ there exists a pair of models $\A,\B$ such that 
  $\A\not\cong\B$ but $\PlTwo\uparrow\EF^\k_t(\A,\B)$. Theorem \ref{thm:BorelTree} now implies that
  the isomorphism relation is not Borel.

  If $T$ is not classifiable $\k$ is not weakly inaccessible, 
  then by \cite{Shelah3} Theorem 0.2 (Main Conclusion), there are non-isomorphic models
  of $T$ of size $\k$ which are $L_{\infty\k}$-equivalent, so the same argument as above,
  using Theorem \ref{thm:BorelTree}, gives that $\cong_T$ is not Borel.
\end{proof}

\begin{Thm}\label{thm:Dii}
  If the theory $T$ is classifiable, then $\cong_T$ is~$\Dii$.
\end{Thm}
\begin{proof}
  Shelah's theorem \cite{Shelah2} XIII.1.4 implies that if a theory $T$ is classifiable, then
  any two models that are $L_{\infty\k}$-equivalent are isomorphic. But $L_{\infty\k}$ equivalence
  is equivalent to $\EF^\k_\o$-equivalence (see Theorem \ref{thm:Karp} on page \pageref{thm:Karp}).
  So in order to prove the theorem it is sufficient to show that if
  for any two models $\A$, $\B$ of the theory $T$ it holds that $\PlTwo \uparrow \EF^\k_\o(\A,\B)\iff \A\cong\B$,
  then the isomorphism relation is $\Dii$.
  The game $\EF^\k_\o$ is a closed game of length $\omega$ and so determined. Hence we have 
  $\PlOne\uparrow \EF^\k_\o(\A,\B)\iff \A\not\cong\B$. 
  By Theorem \ref{thm:EFClosed} the set 
  $$\{(\nu,\eta,\xi)\in (\k^\k)^3\mid \nu \text{ codes a winning strategy for }\PlOne\uparrow \EF^\k_\o(\A_\eta,\A_\xi))\}$$
  is closed and thus $\{(\eta,\xi)\mid \A_\eta\not\cong\A_\xi\}$ is $\Sii$, which further implies that
  $\cong_T$ is $\Dii$ by Corollary~\ref{cor:IsomisS}.
\end{proof}

\section{Unclassifiable}
\subsection{The Unstable, DOP and OTOP Cases}

As before, $\k$ is a regular cardinal satisfying $\k^{<\k}=\k>\o$.
\begin{Thm}\label{thm:NotDiiOrNotBorelList}
  \begin{myEnumerate}
  \item If $T$ is unstable then $\cong_T$ is not~$\Dii$.
  \item If $T$ is stable with OTOP, then $\cong_T$ is not~$\Dii$.
  \item If $T$ is superstable with DOP and $\k>\o_1$, then $\cong_T$ is not~$\Dii$.
  \item If $T$ is stable with DOP and $\l=\cf(\l)=\l(T)+\l^{<\k(T)}\ge\o_1,\k>\l^+$ and for all $\xi<\k$, $\xi^\l<\k$, 
    then $\cong_T$ is not~$\Dii$. (Note that $\k(T)\in \{\o,\o_1\}$.)
  \end{myEnumerate}
\end{Thm}
\begin{proof} 
  For a model $\A$ of size $\k$ of a theory $T$ let us denote by $E(\A)$ the following property: 
  for every $\k^{+}\k$-tree $t$ there is a model $\B$ of $T$ of cardinality
  $\k$ such that $\PlTwo\uparrow \EF^\k_t(\A,\B)$ and $\A\not\cong\B$.

  For (3) we need a result by Hyttinen and Tuuri, Theorem 6.2. from \cite{HytTuu}:
  \begin{Fact}[Superstable with DOP]\label{thm:EFmodels2}
    Let $T$ be a superstable theory with DOP and $\k^{<\k}=\k>\o_1$. Then there exists a model $\A$ of $T$ of cardinality
    $\k$ with the property $E(\A)$.
  \end{Fact}
  
 For (4) we will need a result by Hyttinen and Shelah from~\cite{HytSheC}:
  \begin{Fact}[Stable with DOP]\label{thm:EFmodels1}
    Let $T$ be a stable theory with DOP and $\l=\cf(\l)=\l(T)+\l^{<\k(T)}\ge\o_1$, $\k^{<\k}=\k>\l^+$
    and for all $\xi<\k$, $\xi^{\l}<\k$. Then there is a model $\A$ of $T$ of power $\k$ with the property $E(\A)$.
  \end{Fact}
  For (1) a result by Hyttinen and Tuuri Theorem~4.9 from~\cite{HytTuu}: 
  \begin{Fact}[Unstable]\label{thm:EFmodels2}
    Let $T$ be an unstable theory. Then there exists a model $\A$ of $T$ of cardinality
    $\k$ with the property $E(\A)$.
  \end{Fact}
  And for (2) another result by Hyttinen and Tuuri, Theorem 6.6 in \cite{HytTuu}:
  \begin{Fact}[Stable with OTOP]\label{thm:EFmodels3}
    Suppose $T$ is a stable theory with OTOP. Then there exists a model $\A$ of $T$ of cardinality
    $\k$ with the property $E(\A)$.
  \end{Fact}
  Now (1), (2) and (4) follow immediately from Theorem~\ref{thm:NotD1}. 
\end{proof}

\subsection{Stable Unsuperstable}

We assume $\k^{<\k}=\k>\o$ in all theorems below.

\begin{Thm}\label{thm:StabUnsupstab} Assume that for all $\l<\k$, $\l^\o<\k$.
  \begin{myEnumerate}
    \item  If $T$ is stable unsuperstable, then $\cong_T$ is not Borel. 
    \item If $\k$ is as above and $T$ is stable unsuperstable, 
      then $\cong_T$ is not $\Dii$ in the forcing extension after
      adding $\k^+$ Cohen subsets of~$\k$, or if $V=L$.
  \end{myEnumerate}
\end{Thm}
\begin{proof}
  By Theorem \ref{thm:NSIsReducibleToStableUnsuperstable} on page 
  \pageref{thm:NSIsReducibleToStableUnsuperstable} the relation $E_{S^\k_\o}$ can be reduced to $\cong_T$. 
  The theorem follows now from Corollary \ref{cor:NSIsNotBorel} on page~\pageref{cor:NSIsNotBorel}.
\end{proof}

On the other hand, stable unsuperstable theories sometimes behave nicely to some extent: 

\begin{Lemma}\label{lemma:BorelSt}
  Assume that $T$ is a theory and $t$ a $\k^+\k$-tree such that if $\A$ and $\B$ are models of $T$,
  then $\A\cong \B \iff \PlTwo\uparrow \EF^\k_t(\A,\B)$. Then $\cong$ of $T$ is Borel*.
\end{Lemma}
\begin{proof}
Similar to the proof of Theorem \ref{thm:BorelTree}.
\end{proof}

\begin{Thm}\label{cor:UnsupBorelStar}
  Assume $\k\in I[\k]$ and $\k=\l^+$ (``$\k\in I[\k]$'' is known as the Approachability Property and follows from
  $\l^{<\l}=\l$). Then there exists an unsuperstable theory $T$ whose isomorphism relation is Borel*.
\end{Thm}
\begin{proof}
  In \cite{HytSheA} and \cite{HytSheB} Hyttinen and Shelah show the following
  (Theorem 1.1 of \cite{HytSheB}, but the 
  proof is essentially in~\cite{HytSheA}):
  \begin{myItemize}
  \item[] Suppose $T=((\o^\o,E_i)_{i<\o})$, where $\eta E_i\xi$ if and only if for all $j\le i$, 
    $\eta(j)=\xi(j)$. If $\k\in I[\k]$, $\k=\l^+$ and $\A$ and $\B$ are models of $T$
    of cardinality $\k$, then $\A\cong\B\iff \PlTwo\uparrow \EF^\k_{\l\cdot\o + 2}(\A,\B)$, where
    $+$ and $\cdot$ denote the ordinal sum and product, i.e. $\l\cdot\o + 2$ is just an ordinal.
  \end{myItemize}
  So taking the tree $t$ to be $\l\cdot\o+2$ the claim follows from Lemma~\ref{lemma:BorelSt}.
\end{proof}

\begin{Open}
  We proved that the isomorphism relation of a theory $T$ is Borel if and only if $T$ is classifiable and
  shallow. Is there a 
  connection between the depth of a shallow theory and the Borel degree of its isomorphism relation?
  Is one monotone in the other?
\end{Open}

\begin{Open}
  Can it be proved in ZFC that if $T$ is stable unsuperstable then $\cong_T$ is not $\Dii$?
\end{Open}

\chapter{Reductions}

Recall that in Chapter \ref{chapter:ComplexityofIsomRel} we obtained a provable characterization of
theories which are both classifiable and shallow in terms of the
definability of their isomorphism relations. Without the shallowness
condition we obtained only a consistency result. In this chapter we
improve this to a provable characterization by analyzing isomorphism
relations in terms of Borel reducibility.

Recall the definition of a reduction, Section \ref{def:Reductions} and recall that
if $X\subset \k$ be a stationary subset, we denote by $E_X$ the equivalence relation defined by
$$\forall\eta,\xi\in 2^{\k}(\eta E_X \xi\iff (\eta^{-1}\{1\}\sd\xi^{-1}\{1\})\cap X\text{ is non-stationary}),$$
and by $S^\k_\l$ we mean the ordinals of cofinality $\l$ that are less than $\k$.

The equivalence relations $E_X$ are $\Sii$ ($A E_X B$ if and only if \emph{there exists} 
a cub subset of $\k\setminus (X\cap(A\sd B))$).

Simple conclusions can readily be made from the following observation that
roughly speaking, the set theoretic complexity of a relation does not decrease under reductions:

\begin{Fact}\label{fact:ReductionNote}
  If $E_1$ is a Borel (or $\Dii$) equivalence relation and $E_0$ is an equivalence relation with
  $E_0\le_B E_1$, then $E_0$ is Borel (respectively $\Dii$ if $E_1$ \mbox{is $\Dii$}).  \qed
\end{Fact}

The main theorem of this chapter is:

\begin{Thm}\label{thm:NSembeddable}
  Suppose $\k=\l^+=2^\l>2^\o$ where $\l^{<\l}=\l$. Let
  $T$ be a first-order theory. Then $T$ is classifiable if and only if for all regular $\mu<\k$,
  \mbox{$E_{S^\k_\mu}\not\le_B\,\cong_T\!\!.$}
\end{Thm}

\section{Classifiable Theories}
\label{sec:NSNotToClass}

The following follows from \cite{Shelah2} Theorem XIII.1.4.

\begin{Thm}[\cite{Shelah2}]\label{thm:ClassifiableEF}
  If a first-order theory $T$ is classifiable and $\A$ and $\B$ are non-isomorphic models of $T$ of size $\k$, then
  $\PlOne \uparrow \EF^\k_\o(\A,\B)$.\qed
\end{Thm}

\begin{Thm}[$\k^{<\k}=\k$]\label{thm:NS1}
  If a first-order theory $T$ is classifiable, then for all $\l<\k$ 
  $$E_{S^\k_\l}\not\le_B \,\,\cong_T.$$
\end{Thm}
\begin{proof}
  Let $\NS\in \{E_{S^\k_\l}\mid \l\in\reg(\k)\}$.
  
  Suppose $r\colon 2^\k\to 2^\k$ is a Borel function such that 
  $$\forall \eta,\xi\in 2^{\k}(\A_{r(\eta)}\models T\land \A_{r(\xi)}\models T\land(\eta\,{\NS}\,\xi\iff \A_{r(\eta)}\cong \A_{r(\xi)})).\eqno(\nabla)$$
  
  By Lemma \ref{lem:BaireCont}, let $D$ be an 
  intersection of $\k$-many dense open sets such that $R=r\rest D$ is continuous. $D$ can be coded into a function
  $v\colon \k\times\k\to \k^{<\k}$ such that $D=\Cap_{i<\k}\Cup_{j<\k}N_{v(i,j)}$.
  Since $R$ is continuous, it can also be coded into a single function $u\colon\k^{<\k}\times\k^{<\k}\to \{0,1\}$ such that
  $$R(\eta)=\xi\iff (\forall\a<\k)(\exists\b<\k)[u(\eta\rest\b,\xi\rest\a)=1].$$ 
  (For example define $u(p,q)=1$ if $D\cap N_{p}\subset R^{-1}[N_{q}]$.)
  Let 
  $$\f(\eta,\xi,u,v)=(\forall\a<\k)(\exists\b<\k)[u(\eta\rest\b,\xi\rest\a)=1]\land (\forall i<\k)(\exists j<\k)[\eta\in N_{v(i,j)}].$$
  It is a formula
  of set theory with parameters $u$ and $v$.
  It is easily seen that $\f$ is 
  absolute for transitive elementary submodels $M$ of $H(\k^+)$ containing $\k$, $u$ and $v$ with $(\k^{<\k})^M=\k^{<\k}$.
  Let $\P=2^{<\k}$ be the Cohen forcing. 
  Suppose $M\preccurlyeq H(\k^+)$ is a model as above, i.e. transitive, $\k,u,v\in M$ and $(\k^{<\k})^M=\k^{<\k}$. 
  Note that then $\P\cup \{\P\}\subset M$.
  Then, if $G$ is $\P$-generic over $M$, then $\cup G\in D$ and
  there is $\xi$ such that $\f(\cup G,\xi,u,v)$. By the definition of $\f$ and $u$,
  an initial segment of $\xi$ can be read from an initial segment of $\cup G$. That is why there is
  a nice $\P$-name $\tau$ for a function \mbox{(see \cite{Kunen})} such that
  $$\f(\cup G,\tau_G, u,v)$$
  whenever $G$ is $\P$-generic over $M$. 
                                                 
  Now since the game $\EF^\k_\o$ is determined on all structures, (at least) one of the following holds:
  \begin{myEnumerate}
  \item there is $p$ such that $p\forces \PlTwo\uparrow \EF^\k_\o(\A_\tau,\A_{r(\bar 0)})$
  \item there is $p$ such that $p\forces \PlOne\uparrow \EF^\k_\o(\A_\tau,\A_{r(\bar 0)})$
  \end{myEnumerate}
  where $\bar 0$ is the constant function with value $0$.
  Let us show that both of them lead to a contradiction. 

  Assume (1). Fix a nice $\P$-name $\sigma$ such that 
  $$p\forces\text{``}\sigma\text{ is a winning strategy of }\PlTwo\text{ in }\EF^\k_\o(\A_\tau,\A_{r(\bar 0)})\text{''}$$
  A strategy is a subset of $([\k]^{<\k})^{<\o}\times \k^{<\k}$ (see Definition \ref{def:EF1}), 
  and the forcing does not add elements to that set, 
  so the nice name can be chosen such that all names in $\dom \sigma$
  are standard names for elements that are in $([\k]^{<\k})^{<\o}\times \k^{<\k}\in H(\k^+)$.

  Let $M$ be an elementary submodel of $H(\k^+)$ of size $\k$ such that
  $$\{u,v,\sigma,r(\bar 0),\tau,\P\}\cup (\k+1)\cup M^{<\k}\subset M.$$
  Listing all dense subsets of $\P$ in $M$, it is easy to find a $\P$-generic $G$ over $M$
  which contains $p$ and such that $(\cup G)^{-1}\{1\}$ contains a cub.
  Now in $V$, $\cup G\,\,\,\text{\raisebox{1pt}{$\diagup$}}\!\!\!\!\!\!\!\NS \,\bar 0$. Since $\f(\cup G,\tau_G,u,v)$
  holds, we have by $(\nabla)$: 
  $$\A_{\tau_G}\not\cong\A_{r(\bar 0)}.\eqno(i)$$
  Let us show that $\sigma_G$ is a winning strategy of player $\PlTwo$ in 
  $\EF^\k_\o(\A_{\tau_G},\A_{r(\bar 0)})$ (in $V$) which by Theorem \ref{thm:ClassifiableEF} 
  above is a contradiction with $(i)$.
  
  Let $\mu$ be any strategy of player $\PlOne$ in $\EF^\k_\o(\A_{\tau_G},\A_{r(\bar 0)})$ and 
  let us show that $\sigma_G$ beats it. Consider the play $\sigma_G*\mu$ and assume for a contradiction that
  it is a win for $\PlOne$. This play is well defined, 
  since the moves made by $\mu$ are in the domain of $\sigma_G$ by the note after the definition of $\sigma$, 
  and because $([\k]^{<\k})^{<\o}\times \k^{<\k}\subset M$.

  The play consists of $\o$ moves and is a countable sequence in the set 
  $([\k]^{<\k})\times \k^{<\k}$ (see Definition of EF-games \ref{def:EF1}). Since $\P$ is $<\k$ closed, 
  there is $q_0\in \P$ which decides $\sigma_G*\mu$ (i.e. $\sigma_{G_0}*\mu=\sigma_{G_1}*\mu$ whenever
  $q_0\in G_0\cap G_1$).
  Assume that $G'$ is a $\P$-generic over $V$ with $q_0\in G'$. Then  
  $$(\sigma_{G'}*\mu)^{V[G']}=(\sigma_{G}*\mu)^{V[G']}=(\sigma_G*\mu)^{V}$$ 
  (again, because $\P$ does not add elements of $\k^{<\k}$)  and so 
  $$(\sigma_{G'}*\mu\text{ is a win for }\PlOne)^{V[G']}$$
  But $q_0\forces$ ``$\sigma*\mu\text{ is a win for }\PlTwo$'', because $q_0$ extends $p$ and by the choice of $\sigma$.

  The case (2) is similar, just instead of choosing $\cup G$ such that $(\cup G)^{-1}\{1\}$ contains a cub, choose
  $G$ such that $(\cup G)^{-1}\{0\}$ contains a cub. Then we should have $\A_{\tau_G}\cong \A_{r(\bar 0)}$
  which contradicts (2) by the same absoluteness argument as above.
%
%
%
%
\end{proof}

\section{Unstable and Superstable Theories}
\label{sec:NSToNonClass}
In this section we use Shelah's ideas on how to prove non-structure theorems using Ehrenfeucht-Mostowski models,
see \cite{Shelah3}. We use the definition of Ehrenfeucht-Mostowski models from \cite{HytTuu}, Definition 4.2.

\begin{Def}
  In the following discussion of linear orderings we use the following concepts. 
  \begin{myItemize}
  \item \emph{Coinitiality} or \emph{reverse cofinality} of a linear order $\eta$, denoted $\cf^*(\eta)$
    is the smallest ordinal $\a$ such that there is a map $f\colon \a\to \eta$ 
    which is strictly decreasing and $\ran f$ has no (strict) lower bound in $\eta$.
  \item  If $\eta=\la\eta,<\ra$ is a linear ordering, by $\eta^*$ we denote its mirror image:
    $\eta^*=\la \eta,<^*\ra$ where $x<^* y\iff y<x$.
  \item  Suppose $\l$ is a cardinal. We say that an ordering $\eta$ is $\l$-dense if for all subsets
    $A$ and $B$ of $\eta$ with the properties $\forall a\in A\forall b\in B(a<b)$ and $|A|<\l$
    and $|B|<\l$ there is $x\in \eta$ such that $a<x<b$ for all $a\in A$, $b\in B$. 
    Dense means $\o$-dense.
  \end{myItemize}
\end{Def}

\begin{Thm}\label{thm:NSIsEmbeddableToOTOPDOPUnstable}
  Suppose that $\k=\l^+=2^\l$ such that $\l^{<\l}=\l$.
  If $T$ is unstable or superstable with OTOP, then $E_{S^{\k}_{\l}}\le_c \,\,\cong_T$.
  If additionally $\l\ge 2^\o$, then $E_{S^{\k}_{\l}}\le_c \,\,\cong_T$ holds also for superstable $T$
  with DOP.
\end{Thm}
\vspace{-15pt}
\begin{proofV}{Theorem \ref{thm:NSIsEmbeddableToOTOPDOPUnstable}}
  We will carry out the proof for the case where $T$ is unstable and 
  shall make remarks on how certain steps of the proof should be modified in order this to work for superstable theories with DOP or OTOP.
  First for each $S\subset S^{\k}_{\l}$, let us construct the linear orders $\Phi(S)$ which will serve a fundamental role in the construction.
  The following claim is Lemma 7.17 \mbox{in \cite{HuuHytRau}}:
  \begin{claim}{1}\label{page:claimone}
    For each cardinal $\mu$ of uncountable cofinality there exists a linear ordering $\eta=\eta_\mu$ which satisfies:
    \begin{myEnumerate}
    \item $\eta\cong \eta+\eta$,
    \item for all $\a\le \mu$, $\eta\cong \eta\cdot\a+\eta$,
    \item $\eta\cong \eta\cdot \mu+\eta\cdot \o_1^*$,
    \item $\eta$ is dense,
    \item $|\eta|=\mu$,
    \item $\cf^*(\eta)=\o$.
    \end{myEnumerate}
  \end{claim}
  \vspace{-10pt}  
  \begin{proofVOf}{Claim 1}
    Exactly as in \cite{HuuHytRau}. 
  \end{proofVOf}

  For a set $S\subset S^\k_\l$, define the linear order $\Phi(S)$ as follows:
  $$\Phi(S)=\sum_{i<\k}\tau(i,S),$$
  where $\tau(i,S)=\eta_\l$ if $i\notin S$ and $\tau(i,S)=\eta_\l\cdot \o_1^*$, if $i\in S$.
  Note that $\Phi(S)$ is dense.
  For $\a<\b<\k$ define 
  $$\Phi(S,\a,\b)=\sum_{\a\le i<\b}\tau(i,S).$$  
  (These definitions are also as in \cite{HuuHytRau} although the idea dates back to 
  J. Conway's Ph.D. thesis from the 1960's; they are
  first referred to in \cite{NadelStavi}). From now on denote $\eta=\eta_\l$.
  \begin{claim}{2}
    If $\a\notin S$, then for all $\b\ge\a$ we have $\Phi(S,\a,\b+1)\cong \eta$ and if 
    $\a\in S$, then for all $\b\ge\a$ we have $\Phi(S,\a,\b+1)\cong \eta\cdot \o_1^*$.
  \end{claim}
  \begin{proofVOf}{Claim 2}
    Let us begin by showing the first part, i.e. assume that $\a\notin S$.
    This is also like in \cite{HuuHytRau}. We prove the statement by induction on $\OTP(\b\setminus \a)$. If $\b=\a$, then
    $\Phi(S,\a,\a+1)=\eta$ by the definition of $\Phi$. If $\b=\g+1$ is a successor, then 
    $\b\notin S$, because $S$ contains only limit ordinals, so $\tau(\b,S)=\eta$ and
    $$\Phi(S,\a,\b+1)=\Phi(S,\a,\g+1+1)=\Phi(S,\a,\g+1)+\eta$$
    which by the induction hypothesis and by (1) is isomorphic to $\eta$.
    If $\b\notin S$ is a limit ordinal, then choose a continuous cofinal sequence $s\colon\cf(\b)\to \b$
    such that $s(\g)\notin S$ for all $\g<\cf(\b)$. This is possible since 
    $S$ contains only ordinals of cofinality
    $\l$.
    By the induction hypothesis $\Phi(S,\a,s(0)+1)\cong \eta$, 
    $$\Phi(S,s(\g)+1,s(\g+1)+1)\cong\eta$$ 
    for all successor ordinals $\g<\cf(\b)$, 
    $$\Phi(S,s(\g),s(\g+1)+1)\cong\eta$$ 
    for all limit ordinals $\g<\cf(\b)$ and
    so now 
    $$\Phi(S,\a,\b+1)\cong \eta\cdot \cf(\b)+\eta$$
    which is isomorphic to $\eta$ by (2). If $\b\in S$, then $\cf(\b)=\l$ and we can again choose a cofinal
    sequence $s\colon\l\to \b$ such that $s(\a)$ is not in $S$ for all $\a<\l$. By the induction hypothesis.
    as above,
    $$\Phi(S,\a,\b+1)\cong \eta\cdot \l+\tau(\b,S)$$
    and since $\b\in S$ we have $\tau(\b,S)=\eta\cdot \omega_1^*$, so we have
    $$\Phi(S,\a,\b+1)\cong \eta\cdot \l+\eta\cdot \o_1^*$$
    which by (3) is isomorphic to $\eta$.

    Suppose $\a\in S$. Then $\a+1\notin S$, so by the previous part we have
    $$\Phi(S,\a,\b+1)\cong \tau(\a,S)+\Phi(S,\a+1,\b+1)=\eta\cdot \omega_1^*+\eta=\eta\cdot\omega_1^*.$$
  \end{proofVOf}

  This gives us a way to show that the isomorphism type of $\Phi(S)$ depends only on the $E_{S^\k_\l}$-equivalence class of~$S$:
  \begin{claim}{3}\label{page:claimthree}
    If $S,S'\subset S^{\k}_{\l}$ and $S\sd S'$ is non-stationary, then $\Phi(S)\cong \Phi(S')$.
  \end{claim}
  \begin{proofVOf}{Claim 3}
    Let $C$ be a cub set outside $S\sd S'$. Enumerate it
    $C=\{\a_i\mid i<\k\}$ where $(\a_i)_{i<\k}$ is an increasing and
    continuous sequence. Now $\Phi(S)=\Cup_{i<\k}\Phi(S,\a_i,\a_{i+1})$
    and $\Phi(S')=\Cup_{i<\k}\Phi(S',\a_i,\a_{i+1})$. Note that by the definitions
    these are disjoint unions, so it is enough to show that for all $i<\k$ the orders
    $\Phi(S,\a_i,\a_{i+1})$ and $\Phi(S',\a_i,\a_{i+1})$ are isomorphic. But for all $i<\k$
    $\a_i\in S\iff \a_i\in S'$, so by Claim~2 either
    $$\Phi(S,\a_i,\a_{i+1})\cong \eta\cong \Phi(S',\a_i,\a_{i+1})$$
    (if $\a_i\notin S$) or
    $$\Phi(S,\a_i,\a_{i+1})\cong \eta\cdot\omega^*_1\cong \Phi(S',\a_i,\a_{i+1})$$
    (if $\a_i\in S$).
  \end{proofVOf}

  \begin{Def}\label{def:Kotr}
    $K^{\l}_{tr}$ is the set of $L$-models $\A$ where $L=\{<,\less,(P_\a)_{\a\le \l},h\}$, with the properties
    \begin{myItemize}
    \item $\dom \A\subset I^{\le\l}$ for some linear order $I$.
    \item $\forall x,y\in A(x < y\iff x\subset y)$.
    \item $\forall x\in A(P_\a(x)\iff \lg(x)=\a)$.
    \item $\forall x,y\in A[x\less y\iff \exists z\in A((x,y\in \Succ(z))\land (I\models x < y))]$ 
    \item $h(x,y)$ is the maximal common initial segment of $x$ and $y$.
    \end{myItemize}
  \end{Def}

  For each $S$, define the tree $T(S)\in K^{\l}_{tr}$ by
  \begin{eqnarray*}
    T(S)=\Phi(S)^{<\l}\cup\{\eta\colon \l\to \Phi(S)&\mid &\eta\text{ increasing and }\\
    &&  \cf^*(\Phi(S)\setminus\{x\mid (\exists y\in\ran\eta)(x<y)\})=\o_1\}.
  \end{eqnarray*}
  The relations $<$, $\less$, $P_n$ and $h$ are interpreted in the natural way.
    
  Clearly an isomorphism between $\Phi(S)$ and $\Phi(S')$ induces an isomorphism between $T(S)$ and $T(S')$, thus
  $T(S)\cong T(S')$ if $S\sd S'$ is non-stationary. 

  \begin{claim}{4}
    Suppose $T$ is unstable
    in the vocabulary $v$. Let $T_1$ be $T$ with Skolem functions in
    the Skolemized vocabulary $v_1\supset v$.
    Then there is a function
    $\Po(S^{\k}_{\l})\to \{\A^1\mid \A^1\models T_1, |\A^1|=\k\}$, $S\mapsto \A^1(S)$ which has following properties:
    \begin{myAlphanumerate}
    \item There is a mapping $T(S)\to (\dom\A^1(S))^{n}$ for some $n<\o$, $\eta\mapsto a_\eta$, such that
      $\A^1(S)$ is the Skolem hull of $\{a_\eta\mid \eta\in T(S)\}$, i.e.
      $\{a_{\eta}\mid \eta\in T(S)\}$ is the skeleton of $\A^1(S)$. Denote the skeleton of $\A$ by $\Sk(\A)$.
    \item $\A(S)=\A^{1}(S)\rest v$ is a model of $T$.
    \item $\Sk(\A^{1}(S))$ is indiscernible in $\A^1(S)$, \label{(c)indiscernible}
      i.e. if $\bar\eta,\bar\xi\in T(S)$ and
      $\tp_{\text{q.f.}}(\bar\eta/\es)=\tp_{\text{q.f.}}(\bar\xi/\es)$, then 
      $\tp(a_{\bar\eta}/\es)=\tp(a_{\bar\xi}/\es)$ 
      where
      $a_{\bar\eta}=(a_{\eta_1},\dots,a_{\eta_{\lg\bar \eta}})$. This assignment of types in $\A^1(S)$
      to q.f.-types in $T(S)$ is independent of $S$.
    \item There is a formula $\f\in L_{\o\o}(v)$ 
      such that for all $\eta,\nu\in T(S)$ and $\a<\l$, if $T(S)\models P_{\l}(\eta)\land P_{\a}(\nu)$, 
      then $T(S)\models \eta>\nu$ if and only if 
      $\A(S)\models\f(a_{\eta},a_{\nu})$.
    \end{myAlphanumerate}
  \end{claim}
  \begin{proofVOf}{Claim 4}
    The following is known:
    \begin{myItemize}
    \item[(F1)] Suppose that $T$ is a complete unstable theory. Then for each linear order 
      $\eta$, $T$ has an Ehrenfeucht-Mostowski model $\A$ of vocabulary $v_1$, where $|v_1|=|T|+\o$ and order
      is definable by a first-order formula, such that the template (assignment of types) is independent of
      $\eta$.\footnote{This is from \cite{Shelah4}; there is a sketch of the 
        proof also in \cite{HytTuu}, Theorem 4.7.}
    \end{myItemize}
    It is not hard to see that for every tree $t\in K^\o_{tr}$ we can define a linear order $L(t)$ satisfying the
    following conditions:
    \begin{myEnumerate}
    \item $\dom(L(t))=(\dom t\times\{0\})\cup(\dom t\times\{1\})$,
    \item for all $a\in t$, $(a,0)<_{L(t)}(a,1)$,
    \item if $a,b\in t$, then $a<_t b\iff [(a,0)<_{L(t)}(b,0)]\land [(b,1)<_{L(t)}(a,1)]$,
    \item if $a,b\in t$, then 
      $$(a\not\le b)\land (b\not\le a)\iff [(b,1)<_{L(t)}(a,0)]\lor [(a,1)<_{L(t)}(b,0)].$$
    \end{myEnumerate}
    Now for every $S\subset \k$, by (F1),
    there is an Ehrenfeucht-Mostowski 
    model $\A^1(S)$ for the linear order $L(T(S))$ where order is definable by the formula
    $\psi$ which is in $L_{\infty\o}$.
    Suppose $\bar \eta=(\eta_0,\dots, \eta_n)$ and 
    $\bar \xi=(\xi_0,\dots, \xi_n)$ are sequences in $T(S)$ that have the
    same quantifier free type.
    Then the sequences
    $$\la(\eta_0,0),(\eta_0,1),(\eta_1,0),(\eta_1,1),\dots,(\eta_n,0),(\eta_n,1)\ra$$ 
    and 
    $$\la(\xi_0,0),(\xi_0,1),(\xi_1,0),(\xi_1,1),\dots,(\xi_n,0),(\xi_n,1)\ra$$ 
    have the same quantifier free type in $L(T(S))$. 
    Now let the canonical skeleton of $\A^1(S)$ given by (F1) be $\{a_x\mid x\in L(T(S))\}$.
    Define the $T(S)$-skeleton of $\A^1(S)$ to be the
    set 
    $$\{a_{(\eta,0)}\cat a_{(\eta,1)}\mid \eta\in T(S)\}.$$
    Let us denote $b_\eta=a_{(\eta,0)}\cat a_{(\eta,1)}.$
    This guarantees that (a), (b) and (c) are satisfied. 
    
    For (d) suppose that the order $L(T(S))$ is definable in $\A(S)$ by the formula $\psi(\bar u,\bar c)$, i.e.
    $\A(S)\models \psi(a_x,a_y)\iff x<y$ for $x,y\in L(T(S))$.
    Let $\f(x_0,x_1,y_0,y_1)$ be the formula 
    $$\psi(x_0,y_0)\land \psi(y_1,x_1).$$
    Suppose $\eta,\nu\in T(S)$ are such that $T(S)\models P_\l(\eta)\land P_\a(\nu)$.
    Then
    $$\f((a_{\nu},0),(a_{\nu},1),(a_{\eta},0),(a_{\eta},1))$$
    holds in $\A(S)$ if and only if $\nu<_{T(S)}\eta$. 
  \end{proofVOf}
  \begin{claim}{5}
    Suppose $S\mapsto \A(S)$ is a function as described in Claim~4 with the identical notation. Suppose further
    that $S,S'\subset S^{\k}_{\l}$. Then
    $S\sd S'$ is non-stationary if and only if $\A(S)\cong \A(S')$.
  \end{claim}
  \begin{proofVOf}{Claim 5}
    Suppose $S\sd S'$ is non-stationary. Then by Claim~3 $T(S)\cong T(S')$ which implies $L(T(S))\cong L(T(S'))$
    (defined in the proof of Claim~4) which in turn implies $\A(S)\cong \A(S')$.
    
    Let us now show that if $S\sd S'$ is stationary, then $\A(S)\not\cong \A(S')$.
    Let us make a counter assumption, namely that there is an isomorphism 
    $$f\colon \A(S)\cong \A(S')$$
    and
    that $S\sd S'$ is stationary, and let us deduce a contradiction.
    Without loss of generality we may assume that $S\setminus S'$ is stationary.
    Denote 
    $$X_0=S\setminus S'$$
    
    For all $\a<\k$ define $T^{\a}(S)$ and $T^\a(S')$ by
    $$T^{\a}(S)=\{\eta\in T(S)\mid \ran\eta\subset \Phi(S,0,\b+1)\text{ for some }\b<\a\}$$
    and
    $$T^{\a}(S')=\{\eta\in T(S)\mid \ran\eta\subset \Phi(S',0,\b+1)\text{ for some }\b<\a\}.$$
    Then we have:
    
    \begin{myItemize}
    \item[(i)] if $\a<\b$, then $T^{\a}(S)\subset T^{\b}(S)$
    \item[(ii)] if $\g$ is a limit ordinal, then $T^{\g}(S)=\Cup_{\a<\g}T^{\a}(S)$
    \end{myItemize}
    The same of course holds for $S'$. 
    Note that if $\a\in S\setminus S'$, then
    there is $\eta\in T^{\a}(S)$ cofinal in $\Phi(S,0,\a)$ but there is no such $\eta\in T^\a(S')$ by definition of $\Phi$:
    a cofinal function $\eta$ is added only if $\cf^*(\Phi(S',\a,\k))=\o_1$ which it is not if $\a\notin S'$
    This is the key
    to achieving the contradiction.
    
    But the clauses (i),(ii) are not sufficient to carry out the following argument, because we would like to have
    $|T^\a(S)|<\k$. 
    That is why we want to define a different kind of filtration for $T(S)$, $T(S')$.
    
     For all $\a\in X_0$ fix a function 
     $$\eta^\a_\l\in T(S)\eqno(*)$$ 
     such that $\dom\eta^\a_\l=\l$, for all 
     $\b<\l$, $\eta^\a_\l\rest\b\in T^{\a}(S)$ and $\eta^\a_\l\notin T^{\a}(S)$.

     For arbitrary $A\subset T(S)\cup T(S')$ let $\cl_\Sk(A)$ be the set $X\subset \A(S)\cup \A(S')$
     such that $X\cap \A(S)$ is the Skolem closure of $\{a_{\eta}\mid \eta\in A\cap T(S)\}$ and
     $X\cap \A(S')$ the Skolem closure of $\{a_\eta\mid \eta\in A\cap T(S')\}$. The following is easily verified:

    There exists a $\l$-cub set $C$ and a set $K^\a\subset T^{\a}(S)\cup T^{\a}(S')$ for each $\a\in C$
    such that
    \begin{myItemize}
    \item[(i')] If $\a<\b$, then $K^{\a}\subset K^{\b}$
    \item[(ii')] If $\g$ is a limit ordinal in $C$, then $K^{\g}=\Cup_{\a\in C\cap\g}K^{\a}$
    \item[(iii)] for all $\b<\a$, $\eta^\b_\l\in K^{\a}$. (see $(*)$ above)
    \item[(iv)] $|K^{\a}|=\l$.
    \item[(v)]  $\cl_\Sk(K^{\a})$ is closed under $f\cup f^{-1}$. 
    \item[(vi)] $\{\eta\in T^{\a}(S)\cup T^\a(S')\mid \dom\eta<\l\}\subset K^\a$.
    \item[(vii)] $K^{\a}$ is downward closed.
    \end{myItemize}

    Denote $K^\k=\Cup_{\a<\k}K^\a$. Clearly $K^\k$ is closed under $f\cup f^{-1}$ and so $f$ is an isomorphism between 
    $\A(S)\cap \cl_\Sk(K^{\k})$ and    $\A(S')\cap \cl_\Sk(K^{\k})$.  We will derive a contradiction from this, i.e. we will actually
    show that $\A(S)\cap \cl_\Sk(K^{\k})$ and    $\A(S')\cap \cl_\Sk(K^{\k})$ cannot be isomorphic by $f$.
    Clauses (iii), (v), (vi) and (vii) guarantee that all elements we are going to deal with will be in $K^\k$.

    Let 
    $$X_1=X_0\cap C.$$ 
    For $\a\in X_1$ let us use the following abbreviations:
    \begin{myItemize}
    \item By $\A_\a(S)$ denote the Skolem closure of $\{a_\eta\mid \eta\in K^{\a}\cap T(S)\}$.
    \item By $\A_\a(S')$ denote the Skolem closure of $\{a_\eta\mid \eta\in K^{\a}\cap T(S')\}$. 
    \item $K^{\a}(S)=K^{\a}\cap T(S)$.
    \item $K^{\a}(S')=K^{\a}\cap T(S')$.
    \end{myItemize}

    In the following we will often deal with finite sequences. When defining such a sequence we will use a bar,
    but afterwards we will not use the bar in the notation (e.g. let $a=\bar a$ be a finite sequence...).

    Suppose $\a\in X_1$. Choose 
    $$\xi^\a_\l=\bar\xi^\a_\l\in T(S')\eqno(**)$$ 
    to be such that for some (finite sequence of) 
    terms $\pi=\bar\pi$
    we have 
    \begin{eqnarray*}
    &&f(a_{\eta^\a_\l})=\pi(a_{\xi^\a_\l})\\
    &=&\la \pi_1(a_{\xi^\a_\l(1)},\dots,a_{\xi^\a_\l(\lg(\bar\xi^\a_\l))}),\dots\pi_{\lg \bar\pi}(a_{\xi^\a_\l(1)},
       \dots,a_{\xi^\a_\l(\lg(\xi^\a_\l))})\ra.
    \end{eqnarray*}
    Note that $\xi^\a_\l$ is in $K^\k$ by the definition of $K^\a$'s.

    $$\text{Let us denote by }\eta^\a_\b\text{, the element }\eta^\a_\l\rest \b.\eqno(***)$$

    Let 
    $$\xi^\a_*=\{\nu\in T(S')\mid \exists \xi\in \xi^\a_\l(\nu<\xi)\}.$$
    Also note that $\xi^\a_*\subset K^\b$ for some $\b$. 

    Next define the function $g\colon X_1\to\k$ as follows. Suppose $\a\in X_1$. 
    Let $g(\a)$ be the smallest ordinal $\b$ such that $\xi^\a_*\cap K^{\a}(S')\subset K^{\b}(S')$. We claim that
    $g(\a)<\a$. Clearly $g(\a)\le \a$, so suppose that $g(\a)=\a$. 
    Since $\xi^\a_\l$ is finite, there must be a $\xi^\a_\l(i)\in  \xi^\a_\l$ such that
    for all $\b<\a$ there exists $\g$ such that $\xi^\a_\l(i)\rest \g\in K^\a(S')\setminus K^{\b}(S')$, i.e.
    $\xi^\a_\l(i)$ is cofinal in $\Phi(S',0,\a)$ which it cannot be, because $\a\notin S'$.
    
    Now by Fodor's lemma there exists a stationary set
    $$X_2\subset X_1$$
    and $\g_0$ such that $g[X_2]=\{\g_0\}$. 

    Since there is only $<\k$ many finite sequences in 
    $\A_{\g_0}(S')$, there is a stationary set
    $$X_3\subset X_2$$
    and a finite sequence $\xi=\bar\xi\in K^{\g_0}(S')$
    such that for all $\a\in X_3$ we have $\xi^\a_*\cap K^{\g_0}(S')=\xi_*$ where $\xi_*$ is the set
    $$\xi_*=\{\nu\in T(S')\mid \nu\le\zeta\text{ for some }\zeta\in \bar\xi\}\subset K^{\g_0}(S').$$
    
    Let us fix a (finite sequence of) term(s) $\pi=\bar\pi$ such that the set
    \begin{eqnarray*}
      X_4=\{\a\in X_3&\mid &f(a_{\eta^\a_\l})=\pi(a_{\xi^\a_\l})\}      
    \end{eqnarray*}
    is stationary (see $(*)$). Here $f(\bar a)$ means $\la f(a_1),\dots,f(a_{\lg\bar a})\ra$ and 
    $\bar\pi(\bar b)$ means 
    $$\la \pi_1(b_1,\dots,b_{\lg\bar a}),\dots, \pi_{\lg\pi}(b_1,\dots,b_{\lg\bar a})\ra.$$
    We can find such $\pi$ because there are only countably many such finite sequences of terms.
    
    We claim that in $T(S')$ there are at most $\l$ many quantifier free types over $\xi_*$. 
    All types from now on are quantifier free.
    Let us show that there are at most $\l$ many $1$-types; the general case is left to 
    the reader. To see this, note 
    that a type $p$ over $\xi_*$ is described by the triple 
    $$(\nu_p,\b_p,m_p)\eqno(\star)$$ 
    defined as follows: if $\eta$ satisfies $p$, then 
    $\nu_p$ is the maximal element of $\xi_*$ that is an initial segment of $\eta$,
    $\b_p$ is the level of $\eta$ and $m_p$ 
    tells how many elements of $\xi_*\cap P_{\dom\nu_p+1}$ are there $\less$-below $\eta(\dom\nu_p)$ (recall the vocabulary from
    Definition \ref{def:Kotr}). 
    
    Since $\nu_p\in\xi_*$ and $\xi_*$ is of size $\l$,
    $\b_p\in (\l+1)\cup \{\infty\}$ and $m_p<\o$, there can be at most $\l$ such triples.

    Recall the notations $(*)$, $(**)$ and $(***)$ above.

    We can pick ordinals $\a<\a'$, $\a,\a'\in X_4$, a term $\tau$ and an ordinal $\b<\l$ such that
    $$\eta^{\a'}_{\b}\ne\eta^\a_\b,$$
    $$f(a_{\eta^\a_\b})=\tau(a_{\xi^\a_\b})\text{ and }f(a_{\eta^{\a'}_\b})=\tau(a_{\xi^{\a'}_{\b}})\text{ for some }
    \xi^{\a}_\b,\xi^{\a'}_\b$$
    $$\tp(\xi^\a_\l/\xi_*)=\tp(\xi^{\a'}_\l/\xi_*)$$ 
    and
    $$\tp(\xi^\a_\b/\xi_*)=\tp(\xi^{\a'}_\b/\xi_*).$$ 
   
    We claim that then in fact 
    $$\tp(\xi^\a_\b/(\xi_*\cup \{\xi^{\a'}_l\}))=\tp(\xi^{\a'}_\b/(\xi_*\cup \{\xi^{\a'}\})).$$
    Let us show this. 
    Denote 
    $$p= \tp(\xi^\a_\b /(\xi_*\cup \{\xi^{\a'}_\l\}))$$
    and
    $$p'=\tp(\xi^{\a'}_\b/(\xi_*\cup \{\xi^{\a'}_\l\})).$$
    By the same reasoning as above at $(\star)$ it 
    is sufficient to show that these types $p$ and $p'$ have the same triple of the form $(\star)$.
    Since $\a$ and $\a'$ are in $X_3$ and $X_2$, we have $\xi^{\a'}_*\cap K^{\a'}(S')=\xi_*\subset K^{\g_0}(S')$.
    On the other hand $f\rest\A_{\a'}(S)$ is an isomorphism 
    between $\A_{\a'}(S)$ and $\A_{\a'}(S')$, because $\a$ and $\a'$ are in $X_1$,
    and so $\xi^{\a'}_\b\in K^{\a'}(S')$. Thus $\nu_p=\nu_{p'}\in \xi_*$ and $m_p=m_{p'}$ follows in the same way.
    Clearly $\b_p=\b_{p'}$.

    Now we have: $\xi^\a_\l$ and $\pi$ are such that
    $f(a_{\eta^\a_\l})=\pi(a_{\xi^\a_\l})$ and $\xi^\a_\b$ and $\tau$ are such that $f(a_{\eta^\a_\b})=\tau(a_{\xi^\a_\b})$.
    Similarly for $\a'$. The formula $\f$ is defined in Claim~4.
    
    We know that
    $$\A(S)\models \f(a_{\eta^{\a'}_\l},a_{\eta^{\a'}_\b})$$ 
    and because $f$ is isomorphism, this implies 
    $$\A(S')\models \f(f(a_{\eta^{\a'}_\l}),f(a_{\eta^{\a'}_\b}))$$ 
    which is equivalent to
    $$\A(S')\models \f(\pi(a_{\xi^{\a'}_\l}),\tau(a_{\xi^{\a'}_\b}))$$ 
    (because $\a,\a'$ are in $X_4$).
    Since $T(S')$ is indiscernible in $\A(S')$ and $\xi^{\a'}_\b$ and $\xi^{\a}_\b$ have the same type over
    over $(\xi_*\cup \{\xi^{\a'}_\l\})$, we have
    $$\A(S')\models \f(\pi(a_{\xi^{\a'}_\l}),\tau(a_{\xi^{\a'}_\b}))\iff \f(\pi(a_{\xi^{\a'}_\l}),\tau(a_{\xi^{\a}_\b}))\eqno(*)$$ 
    and so we get
    $$\A(S')\models \f(\pi(a_{\xi^{\a'}_\l}),\tau(a_{\xi^{\a}_\b}))$$ 
    which is equivalent to
    $$\A(S')\models \f(f(a_{\eta^{\a'}_\l}),f(a_{\eta^{\a}_\b}))$$ 
    and this in turn is equivalent to 
    $$\A(S)\models \f(a_{\eta^{\a'}_\l},a_{\eta^{\a}_\b})$$ 
    The latter cannot be true, because the definition of $\b,\a$ and $\a'$ implies that
    $\eta^{\a'}_\b\ne \eta^{\a}_\b$.
  \end{proofVOf}
  
  Thus, the above Claims 1 -- 5 justify the embedding of $E_{S^\k_\l}$ into the isomorphism relation
  on the set of structures that are models for $T$ for unstable $T$. 
  This embedding combined with a suitable coding of models gives a continuous map.\\
  
  \noindent\textbf{DOP and OTOP cases}. The above proof was based on the fact (F1) that for unstable theories there are Ehrenfeucht-Mostowski
  models for any linear order such that the order is definable by a first-order formula $\f$ and is indiscernible relative to $L_{\o\o}$,
  (see (c) on page \pageref{(c)indiscernible}); 
  it is used in $(*)$ above.
  For the OTOP case, we use instead the fact (F2):
  \begin{myItemize}
  \item[(F2)] Suppose that $T$ is a theory with OTOP in a countable vocabulary $v$. 
    Then for each dense linear order $\eta$ we can find a model $\A$ of a countable vocabulary 
    $v_1\supset v$ such that $\A$ is an Ehrenfeucht-Mostowski model of $T$ for $\eta$ where 
    order is definable by an $L_{\o_1\o}$-formula.\footnote{Contained in the 
      proof of Theorem 2.5. of \cite{Shelah1}; see also \cite{HytTuu}, Theorem 6.6.}  
  \end{myItemize}
  Since the order $\Phi(S)$ is dense, it is easy to argue that
  if $T(S)$ is indiscernible relative to $L_{\o\o}$, then it is indiscernible relative to $L_{\infty\o}$
  (define this as in (c) on page \pageref{(c)indiscernible} changing $\tp$ to $\tp_{L_{\infty\o}}$).
  Other parts of the proof remain unchanged, because although the formula $\f$ is not first-order anymore, it is still in $L_{\infty\o}$.
  
  In the DOP case we have the following fact:
  \begin{myItemize}
  \item[(F3)] Let $T$ be a countable superstable theory 
    with DOP of vocabulary $v$. Then there exists a vocabulary $v_1\supset v$,
    $|v_1|=\o_1$, such that for every linear order $\eta$ there exists a $v_1$-model $\A$ which is an Ehrenfeucht-Mostowski
    model of $T$ for $\eta$ where order is definable by an $L_{\o_1\o_1}$-formula.\footnote{This 
      is essentially from \cite{Shelah5} Fact 2.5B; a proof can be found also in \cite{HytTuu} Theorem~6.1.}
  \end{myItemize}
  Now the problem is that $\f$ is in $L_{\infty\o_1}$. By (c) of Claim~4, $T(S)$ is indiscernible in 
  $\A(S)$ relative to $L_{\o\o}$ and by the above relative to $L_{\infty\o}$.  If we could require $\Phi(S)$ to be $\o_1$-dense, we would similarly 
  get indiscernible relative to $L_{\infty\o_1}$. Let us show how to modify the proof in order to do that.
  Recall that in the DOP case,we assume $\l\ge 2^\o$.
  
  In Claim 1 (page \pageref{page:claimone}), we have to replace clauses (3), (4) and (6) by (3'), (4') and (6'):
  \begin{myEnumerate}
  \item[(3')] $\eta\cong \eta\cdot \mu+\eta\cdot \o^*$,
  \item[(4')] $\eta$ is $\o_1$-dense,
  \item[(6')] $\cf^*(\eta)=\o_1$. 
  \end{myEnumerate}
  The proof that such an $\eta$ exists is exactly as the proof of Lemma 7.17 \cite{HuuHytRau} 
  except that instead of putting $\mu=(\o_1)^V$ put $\mu=\o$, build $\theta$-many
  functions with domains being countable initial segments of $\o_1$ 
  instead of finite initial segments of $\o$ and instead of $\Q$ (the countable dense linear order)
  use an $\o_1$-saturated dense linear order -- this order has size $2^\o$ and that is why the assumption $\l\ge 2^\o$ 
  is needed.

  In the definition of $\Phi(S)$ (right after Claim~1), replace $\o_1^*$ by $\o^*$ and $\eta$ by the new $\eta$
  satisfying (3'), (4') and (6') above. Note that $\Phi(S)$ becomes now $\o_1$-dense. In Claim~2 one has to replace
  $\o_1^*$ by $\o^*$. The proof remains similar. In the proof of Claim~3 (page \pageref{page:claimthree}) one has to 
  adjust the use of Claim~2. Then, in the definition of $T(S)$ replace $\o_1$ by $\o$.

  Claim 4 for superstable $T$ with DOP now follows with (c) and (d) modified: instead of indiscernible relative to $L_{\o\o}$,
  demand $L_{\infty\o_1}$ and instead of $\f\in L_{\o\o}$ we have now $\f\in L_{\infty\o_1}$. The proof is unchanged except 
  that the language is replaced by $L_{\infty\o_1}$ everywhere and fact (F1) replaced by (F3) above.

  Everything else in the proof, in particular the proof of Claim~5, remains unchanged modulo some obvious things 
  that are evident from the above explanation.
\end{proofV}
\vspace{-25pt}

\section{Stable Unsuperstable Theories}

In this section we provide a tree construction (Lemma \ref{lem:StoJS}) which is similar to Shelah's construction in \cite{Shelah3}
which he used to obtain (via Ehrenfeucht-Mostowski models) many pairwise non-isomorphic models. Then
using a prime-model construction (proof of Theorem \ref{thm:NSIsReducibleToStableUnsuperstable})
we will obtain the needed result.

\label{sec:NSToNonClass}
\begin{Def}\label{def:Filtration}
  Let $I$ be a tree of size $\k$. Suppose $(I_{\a})_{\a<\k}$ is a collection of subsets of $I$ such that
  \begin{myItemize}
  \item For each $\a<\k$, $I_\a$ is a downward closed subset of $I$
  \item $\Cup_{\a<\k} I_{\a}=I$
  \item If $\a<\b<\k$, then $I_{\a}\subset I_{\b}$
  \item If $\g$ is a limit ordinal, then $I_\g=\Cup_{\a<\g}I_\a$
  \item For each $\a<\k$ the cardinality of $I_\a$ is less than $\k$. 
  \end{myItemize}
  Such a sequence $(I_{\a})_{\a<\k}$ is called $\k$-\emph{filtration} or just \emph{filtration} of $I$. 
\end{Def}

\begin{Def}\label{def:KotrStar}
  Recall $K^\l_{tr}$ from Definition \ref{def:Kotr} on page \pageref{def:Kotr}. Let
  $K^{\l}_{tr*}=\{A\rest L^*\mid A\in K^{\l}_{tr}\}$, where $L^*$ is the vocabulary $\{<\}$.
\end{Def}

\begin{Def}
Suppose $t\in K^{\o}_{tr*}$ is a tree of size $\k$ (i.e. $t\subset \k^{\le\o}$) and let $\I=(I_{\a})_{\a<\k}$ be a filtration of $t$.
Define
$$S_\I(t)=\Big\{\a<\k\mid (\exists\eta\in t)\big[(\dom\eta=\o)\land \forall n<\o (\eta\rest n\in I_{\a})\land(\eta\notin I_\a)\big]\Big\}$$
By $S\sim_{\NS} S'$ we mean that $S\sd S'$ is not $\o$-stationary
\end{Def}

\begin{Lemma}\label{lemma:FiltrEquiv}
  Suppose trees $t_0$ and $t_1$ are isomorphic, and
  $\I=(I_\a)_{\a<\k}$ and $\J=(J_\a)_{\a<\k}$ are $\k$-filtrations of $t_0$ and $t_1$ respectively.
  Then $S_\I(t_0)\sim_{\NS}S_\J(t_1)$.  
\end{Lemma}
\begin{proof}
  Let $f\colon t_0\to t_1$ be an isomorphism. Then 
  $f\I=(f[I_\a])_{\a<\k}$ is a filtration of $t_1$ and 
  $$\a\in S_{\I}(t_0)\iff \a\in S_{f\I}(t_1).\eqno(\star)$$
  Define the set $C=\{\a\mid f[I_\a]=J_\a\}$. Let us show that it is cub. 
  Let $\a\in \k$. Define $\a_0=\a$ and by induction pick $(\a_n)_{n<\o}$ 
  such that $f[I_{\a_n}]\subset J_{\a_{n+1}}$ for odd $n$ and $J_{\a_n}\subset f[I_{\a_{n+1}}]$ for even $n$.
  This is possible by the definition of a $\k$-filtration. Then $\a_\o=\Cup_{n<\o}\a_n\in C$.
  Clearly $C$ is closed and $C\subset \k\setminus S_{f\I}(t_1)\sd S_\J(t_1),$ so now by $(\star)$
  \begin{equation*}
    S_{\I}(t_0)=S_{f\I}(t_1)\sim_{\NS} S_{\J}(t_1).\qedhere    
  \end{equation*}
\end{proof}

\begin{Lemma}\label{lem:StoJS}
  Suppose for $\l<\k$, $\l^\o<\k$ and $\k^{<\k}=\k$.
  There exists a function $J\colon \Po(\k)\to K^\o_{tr*}$ such that 
  \begin{myItemize}
  \item $\forall S\subset\k (|J(S)|=\k)$.
  \item If $S\subset \k$ and $\I$ is a $\k$ filtration of $J(S)$, then $S_\I(J(S))\sim_\NS S$.
  \item If $S_0\sim_\NS S_1$, then $J(S_0)\cong J(S_1)$.
  \end{myItemize}
\end{Lemma}
\begin{proofV}{Lemma \ref{lem:StoJS}}  
  Let $S\subset S^\k_\o$ and let us define a preliminary tree $I(S)$
  as follows. For each $\a\in S$ let $C_\a$ be the set of all strictly increasing cofinal functions 
  $\eta\colon\o\to\a$.
  Let $I(S)=\underline{[\k]}^{<\o}\cup \Cup_{\a\in S}C_\a$ where $\underline{[\k]}^{<\o}$ is the set of strictly increasing functions
  from finite ordinals to $\k$.

  For ordinals $\a<\b\le \k$ and $i<\o$ we adopt the notation:
  \begin{myItemize}
  \item $[\a,\b]=\{\g\mid \a\le \g\le\b\}$
  \item $[\a,\b)=\{\g\mid \a\le \g <\b\}$
  \item $\ruszh(\a,\b,i)=\Cup_{i\le j\le\o}\{\eta\colon [i,j)\to [\a,\b)\mid \eta\text{ strictly increasing}\}$
  \end{myItemize}

  For each $\a,\b<\k$ let us define the sets $P^{\a,\b}_{\g}$, for $\g<\k$ as follows. 
  If $\a=\b=\g=0$, then $P^{0,0}_{0}=I(S)$. Otherwise let $\{P^{\a,\b}_\g\mid \g<\k\}$
  enumerate all downward closed subsets of $\ruszh(\a,\b,i)$ for all $i$, i.e. 
  $$\{P^{\a,\b}_\g\mid \g<\k\}=\Cup_{i<\o}\Po(\ruszh(\a,\b,i))\cap \{A\mid A\text{ is closed under inital segments}\}.$$
  Define 
  $$\rusch(P^{\a,\b}_{\g})$$ 
  to be the natural number $i$ such that
  $P^{\a,\b}_{\g}\subset \ruszh(\a,\b,i)$.
  The enumeration is possible, because by our assumption $\k^{<\k}=\k$ we have
  \begin{eqnarray*}
  \Big|\Cup_{i<\o}\Po(\ruszh(\a,\b,i))\Big|&\le& \o\times |\Po(\ruszh(0,\b,0))|\\    
  &\le&\o\times |\Po(\b^\o)|\\
  &=&\o\times 2^{\b^{\o}}\\
  &\le&\o\times \k\\
  &=&\k
  \end{eqnarray*}
  Let $S\subset \k$ be a set and define $J(S)$ to be the set of all $\eta\colon s\to \o\times \k^4$ such
  that $s\le\o$ and the following conditions are met for all $i,j<s$:
  \begin{myEnumerate}
  \item \label{J0}$\eta$ is strictly increasing with respect to the lexicographical order on $\o\times\k^4$.
  \item \label{J1}$\eta_1(i)\le \eta_1(i+1)\le \eta_1(i)+1$
  \item \label{J3}$\eta_1(i)=0\rightarrow \eta_{2}(i)=\eta_{3}(i)=\eta_{4}(i)=0$
  \item \label{J4}$\eta_1(i)<\eta_1(i+1)\rightarrow \eta_2(i+1)\ge \eta_3(i)+\eta_4(i)$
  \item \label{J5}$\eta_1(i)=\eta_1(i+1)\rightarrow (\forall k\in\{2,3,4\})(\eta_k(i)=\eta_k(i+1))$
  \item \label{J6}if for some $k<\o$, $[i,j)=\eta_1^{-1}\{k\}$, then\\
    $\eta_5\rest[i,j)\in P^{\eta_2(i),\eta_3(i)}_{\eta_4(i)}$
  \item \label{J7}if $s=\o$, then either \\
    $(\exists m<\o)(\forall k<\o)(k>m\rightarrow \eta_1(k)=\eta_1(k+1))$ \\
    or\\
    $\sup\ran\eta_5\in S.$
  \item   Order $J(S)$ by inclusion.\\
  \end{myEnumerate}

  Note that it follows from the definition of $P^{\a,\b}_\g$ and the conditions \eqref{J6} and
  \eqref{J4} that for all $i<j<\dom \eta,$ $\eta\in J(S)$:\\

  \begin{myEnumerate}\setcounter{enumi}{8}
  \item \label{J2} $i<j\rightarrow \eta_5(i)<\eta_5(j)$.\\
  \end{myEnumerate}

  For each $\a<\k$ let 
  $$J^{\a}(S)=\{\eta\in J(S)\mid \ran\eta\subset\o\times(\b+1)^4\text{ for some }\b<\a\}.$$
  Then $(J^{\a}(S))_{\a<\k}$ is a $\k$-filtration of $J(S)$ (see Claim~2 below). For the first
  item of the lemma, clearly $|J(S)|=\k$.

  Let us observe that if $\eta\in J(S)$ and $\ran\eta_1=\o$, then 
  $$\sup\ran\eta_4\le\sup\ran\eta_2=\sup\ran\eta_3=\sup\ran\eta_5\eqno(\#)$$
  and if in addition to that, $\eta\rest k\in J^{\a}(S)$ for all $k$ and $\eta\notin J^{\a}(S)$ or if $\ran\eta_1=\{0\}$, then 
  $$\sup\ran\eta_5=\a.\eqno(\circledast)$$\label{circledast}
  To see $(\#)$ suppose $\ran\eta_1=\o$.  By \eqref{J2}, $(\eta_5(i))_{i<\o}$ is an
  increasing sequence. By \eqref{J6} $\sup\ran\eta_3\ge\sup\ran\eta_5\ge\sup\ran\eta_2$.
  By \eqref{J4}, $\sup\ran\eta_2\ge\sup\ran\eta_3$ and again by \eqref{J4} $\sup\ran\eta_2\ge \sup\ran\eta_4$.
  Inequality $\sup\ran\eta_5\le\a$ is an immediate consequence of the definition of $J^{\a}(S)$,
  so $(\circledast)$ follows now from the assumption that $\eta\notin J^{\a}(S)$.

  \begin{claim}{1}
    Suppose $\xi\in J^{\a}(S)$ and $\eta\in J(S)$. Then
    if $\dom \xi<\o$, $\xi\subsetneq \eta$ and 
    $(\forall k\in \dom\eta\setminus\dom\xi)\big(\eta_1(k)=\xi_1(\max\dom\xi)\land \eta_1(k)>0\big),$
    then
    $\eta\in J^{\a}(S)$.
  \end{claim}
  \begin{proofVOf}{Claim 1}
    Suppose $\xi,\eta\in J^{\a}(S)$ are as in the assumption. Let us define 
    $\b_2=\xi_2(\max\dom\xi)$, $\b_3=\xi_2(\max\dom\xi)$, and
    $\b_4=\xi_4(\max\dom\xi)$. Because $\xi\in J^{\a}(S)$, there is $\b$ such that
    $\b_2,\b_3,\b_4<\b+1$ and $\b<\a$.
    Now by \eqref{J5} $\eta_2(k)=\b_2$, $\eta_3(k)=\b_3$ and $\eta_4(k)=\b_4$, for all $k\in \dom\eta\setminus\dom\xi$. 
    Then by \eqref{J6} for all 
    $k\in\dom\eta\setminus\dom\xi$ we have that 
    $\b_2<\eta_5(k)<\b_3<\b+1$. Since $\xi\in J^{\a}(S)$,
    also $\b_4<\b+1$, so $\eta\in J^{\a}(S)$. 
  \end{proofVOf}
  \begin{claim}{2}
    $|J(S)|=\k$, $(J^\a(S))_{\a<\k}$ is a $\k$-filtration of $J(S)$
    and if $S\subset \k$ and $\I$ is a $\k$-filtration of $J(S)$, then $S_\I(J(S))\sim_\NS S$.
  \end{claim}
  \begin{proofVOf}{Claim 2}
    For all $\a\le\k$,
    $J^\a(S)\subset (\o\times\a^{4})^{\le\o}$, so by the cardinality assumption of the lemma, 
    the cardinality of $J^\a(S)$ is $<\k$ if $\a<\k$  ($J^\k(S)=J(S)$).
    Clearly $\a<\b$ implies $J^{\a}(S)\subset J^{\b}(S)$. Continuity is verified by
    \begin{eqnarray*}
      \Cup_{\a<\g}J^{\a}(S)&=&\{\eta\in J(S)\mid\exists\a<\g,\exists\b<\a(\ran\eta\subset \o\times (\b+1)^4)\}\\
      &=&\{\eta\in J(S)\mid\exists\b<\cup\g(\ran\eta\subset \o\times (\b+1)^4)\}
    \end{eqnarray*}
    which equals $J^\g(S)$ if $\g$ is a limit ordinal.
    By Lemma \ref{lemma:FiltrEquiv} it is enough to show $S_\I(J(S))\sim_\NS S$ for
    $\I=(J^\a(S))_{\a<\k}$, and we will show that if $\I=(J^\a(S))_{\a<\k}$, then 
    in fact $S_\I(J(S))=S$.
    
    Suppose $\a\in S_{\I}(J(S))$. Then there is $\eta\in J(S)$, $\dom\eta=\o$, such that
    $\eta\rest k\in J^{\a}(S)$ for all $k<\o$ but $\eta\notin J^{\a}(S)$. Thus there is no $\b<\a$
    such that $\ran\eta\subset \o\times (\b+1)^4$ but on the other hand for all $k<\o$ there is $\b$
    such that $\ran \eta\rest k\subset \o\times (\b+1)^4$. By \eqref{J5} and \eqref{J6} this implies that
    either $\ran\eta_1=\o$ or $\ran\eta_1=\{0\}$. By ($\circledast$) on page \pageref{circledast} it now follows that
    $\sup\ran\eta_5=\a$ and by \eqref{J7}, $\a\in S$. 

    Suppose then that $\a\in S$. Let us show that $\a\in S_{\I}(J(S))$. Fix a function $\eta_\a\colon \o\to\k$ with 
    $\sup\ran\eta_\a=\a$. Then $\eta_\a\in I(S)$ and the function $\eta$ such that
    $\eta(n)=(0,0,0,0,\eta_\a(n))$ is as required. 
    (Recall that $P^{0,0}_0=I(S)$ in the definition of $J(S)$).
  \end{proofVOf}
  \begin{claim}{3}
    Suppose $S\sim_{\NS}S'$. Then $J(S)\cong J(S')$.  
  \end{claim}
  \begin{proofVOf}{Claim 3}
    Let $C\subset \k\setminus (S\sd S')$ be the cub set which exists by the assumption.
    By induction on $i<\k$ we will define $\a_i$ and $F_{\a_i}$ such that
    \begin{myAlphanumerate}
    \item If $i<j<\k,$ then $\a_i<\a_j$ and $F_{\a_i}\subset F_{\a_j}$.
    \item If $i$ is a successor, then $\a_i$ is a successor and if $i$ is limit, then $\a_i\in C$.
    \item If $\g$ is a limit ordinal, then $\a_\g=\sup_{i<\g}\a_i$,
    \item $F_{\a_i}$ is a partial isomorphism $J(S)\to J(S')$
    \item Suppose that $i=\g+n$, where $\g$ is a limit ordinal or 0 and $n<\o$ is even. Then
      $\dom F_{\a_i}=J^{\a_i}(S)$ (e1). If 
      also $n>0$ and $(\eta_k)_{k<\o}$ is an increasing sequence in $J^{\a_i}(S)$
      such that $\eta=\Cup_{k<\o}\eta_k\notin J(S)$, then $\Cup_{k<\o} F_{\a_i}(\eta_k)\notin J(S')$ (e2).     
    \item If $i=\g+n$, where $\g$ is a limit ordinal or 0 and $n<\o$ is odd, then
      $\ran F_{\a_i}=J^{\a_i}(S')$ (f1).
      Further, if $(\eta_k)_{k<\o}$ is an increasing sequence in $J^{\a_i}(S')$
      such that $\eta=\Cup_{k<\o}\eta_k\notin J(S')$, then $\Cup_{k<\o} F^{-1}_{\a_i}(\eta_k)\notin J(S)$ (f3).
    \item If $\dom \xi<\o$, $\xi\in \dom F_{\a_i}$,
      $\eta\rest\dom\xi=\xi$ and $(\forall k\ge\dom\xi)\big(\eta_1(k)=\xi_1(\max\dom\xi)\land \eta_1(k)>0\big)$, then
      $\eta\in \dom F_{\a_i}$. Similarly for $\ran F_{\a_i}$
    \item If $\xi\in\dom F_{\a_i}$ and $k<\dom \xi$, then $\xi\rest k\in \dom F_{\a_i}$.
    \item For all $\eta\in \dom F_{\a_i}$, $\dom\eta=\dom (F_{\a_i}(\eta))$
    \end{myAlphanumerate}
    
    \noindent\textbf{The first step.} 
    The first step and the successor steps are similar, but the first step is easier. Thus we give it separately
    in order to simplify the readability.
    Let us start with $i=0$. Let $\a_0=\b+1$, for arbitrary $\b\in C$. Let us denote by 
    $$\rusi(\a)$$
    the ordinal that is  order isomorphic to $(\o\times\a^4,<_{\text{lex}})$.
    Let $\g$ be such that there is an isomorphism $h\colon P^{0,\rusi(\a_0)}_{\g}\cong J^{\a_0}(S)$ 
    and such that $\rusch(P^{0,\a_0}_{\g})=0$. Such exists by \eqref{J0}.
    Suppose that $\eta\in J^{\a_0}(S)$. Note that because $P^{0,\a_0}_{\g}$ and $J^{\a_0}(S)$ are closed 
    under initial segments and by the definitions of $\rusch$ and $P^{\a,\b}_\g$, we have $\dom h^{-1}(\eta)=\dom \eta$,
    Define $\xi=F_{\a_0}(\eta)$ such that $\dom\xi=\dom\eta$ 
    and for all $k<\dom \xi$
    \begin{myItemize}
      \item $\xi_1(k)=1$
      \item $\xi_2(k)=0$
      \item $\xi_3(k)=\rusi(\a_0)$
      \item $\xi_4(k)=\g$
      \item $\xi_5(k)=h^{-1}(\eta)(k)$
    \end{myItemize}
    Let us check that $\xi\in J(S')$. Conditions \eqref{J0}-\eqref{J5} and \eqref{J7} are satisfied because
    $\xi_k$ is constant for all $k\in \{1,2,3,4\}$, $\xi_1(i)\ne 0$ for all $i$ and $\xi_5$ is increasing. For \eqref{J6}, if
    $\xi_1^{-1}\{k\}$ is empty, the condition is verified since each $P^{\a,\b}_\g$ is closed under initial segments 
    and contains the empty function. If it is non-empty, then $k=1$ and in that case $\xi_1^{-1}\{k\}=[0,\o)$ and
    by the argument above ($\dom h^{-1}(\eta)=\dom \eta=\dom\xi$) we have
    $\xi_5=h^{-1}(\eta)\in P^{0,\rusi(\a_0)}_\g=P^{\xi_2(0),\xi_3(0)}_{\xi_4(0)}$, so the condition is satisfied.
    
    Let us check whether all the conditions (a)-(i) are met. In (a), (b), (c), (e2) and (f) there is nothing to check.
    (d) holds, because $h$ is an isomorphism. (e1) and (i) are immediate from the definition. 
    Both $J^{\a_0}(S)$ and $P^{0,\rusi(\a_0)}_\g$
    are closed under initial segments, so (h) follows, because $\dom F_{\a_0}=J^{\a_0}(S)$ and 
    $\ran F_{\a_0}=\{1\}\times \{0\}\times \{\rusi(\a_0)\}\times \{\g\}\times P^{0,\a_0}_{\g}$.
    Claim~1 implies (g) for $\dom F_{\a_0}$. Suppose $\xi\in \ran F_{\a_0}$ and $\eta\in J(S')$ are as in the assumption of (g).
    Then $\eta_1(i)=\xi_1(i)=1$ for all $i<\dom \eta$. By \eqref{J5} it follows that 
    $\eta_2(i)=\xi_2(i)=0$,
    $\eta_3(i)=\xi_3(i)=\rusi(\a_0)$ and
    $\eta_4(i)=\xi_4(i)=\g$
    for all $i<\dom \eta$, so by \eqref{J6} $\eta_5\in P^{0,\rusi(\a_0)}_\g$ and since $h$ is an isomorphism,
    $\eta\in\ran F_{\a_0}$.\\

    \noindent\textbf{Odd successor step.}
    We want to handle odd case but not the even case first, because the most important case is the successor of a limit ordinal, see $(\iota\iota\iota)$
    below. Except that, the even case is similar to the odd case.

    Suppose that $j<\k$ is a successor ordinal. Then there exist $\b_j$ and $n_j$ such that
    $j=\b_j+n_j$ and $\b$ is a limit ordinal or $0$. Suppose that $n_j$ is odd and
    that $\a_l$ and $F_{\a_l}$ are defined for all $l<j$ such that the conditions (a)--(i) and \eqref{J0}--\eqref{J2}
    hold for $l<j$.
    
    Let $\a_j=\b+1$ where $\b$ is such that $\b\in C$, $\ran F_{\a_{j-1}}\subset J^{\b}(S')$, $\b>\a_{j-1}$.
    For convenience define $\xi(-1)=(0,0,0,0,0)$ for all $\xi\in J(S)\cup J(S')$.
    Suppose $\eta\in \ran F_{\a_{j-1}}$ has finite domain $\dom\eta=m<\o$ and denote $\xi=F^{-1}_{\a_{j-1}}(\eta)$. 
    Fix $\g_\eta$ to be such that $\rusch(P^{\a,\b}_{\g_\eta})=m$ and
    such that there is an isomorphism 
    $h_\eta\colon P^{\a,\b}_{\g_\eta}\to W,$
    where 
    $$W=\{\zeta\mid \dom\zeta=[m,s), m<s\le\o, 
       \eta^{\frown}\la m,\zeta(m)\ra\notin 
       \ran F_{\a_{j-1}}, \eta^{\frown}\zeta\in J^{\a_j}(S')\},$$
    $\a=\xi_3(m-1)+\xi_4(m-1)$ and $\b=\a+\rusi(\a_j)$ (defined in the beginning of the First step).

    We will define $F_{\a_{j}}$ so that its range is $J^{\a_{j}}(S')$ and instead of $F_{\a_j}$ we will
    define its inverse. 
    So let $\eta\in J^{\a_j}(S')$. We have three cases:
    \begin{myItemize}
    \item[($\iota$)] $\eta\in \ran F_{\a_{j-1}}$,
    \item[($\iota\iota$)] $\exists m<\dom\eta(\eta\rest m\in \ran F_{\a_{j-1}}\land \eta\rest(m+1)\notin F_{\a_{j-1}})$,
    \item[($\iota\iota\iota$)] $\forall m<\dom\eta(\eta\rest(m+1)\in \ran F_{\a_{j-1}}\land \eta\notin \ran F_{\a_{j-1}})$.
    \end{myItemize}
    Let us define $\xi=F^{-1}_{\a_j}(\eta)$ such that $\dom\xi=\dom\eta$. If ($\iota$) holds, define
    $\xi(n)=F^{-1}_{\a_{j-1}}(\eta)(n)$ for all $n<\dom\eta$. Clearly $\xi\in J(S)$ by the induction hypothesis.
    Suppose that ($\iota\iota$) holds
    and let $m$ witness this. For all $n<\dom \xi$ let
    \begin{myItemize}
    \item If $n<m$, then $\xi(n)=F^{-1}_{\a_{j-1}}(\eta\rest m)(n)$.
    \item Suppose $n\ge m$. Let 
      \begin{myEnumerate}
      \item[$\cdot$] $\xi_1(n)=\xi_1(m-1)+1$
      \item[$\cdot$] $\xi_2(n)=\xi_3(m-1)+\xi_4(m-1)$
      \item[$\cdot$] $\xi_3(n)=\xi_2(m)+\rusi(\a_j)$
      \item[$\cdot$] $\xi_4(n)=\g_{\eta\restl m}$
      \item[$\cdot$] $\xi_5(n)=h_{\eta\restl m}^{-1}(\eta)(n)$.
      \end{myEnumerate}
    \end{myItemize}
    Next we should check that $\xi\in J(S)$; let us check items \eqref{J0} and \eqref{J6}, the rest are left to the reader.
    \begin{myItemize}
    \item[\eqref{J0}] By the induction hypothesis $\xi\rest m$ is increasing. Next,
      $\xi_1(m)=\xi_1(m-1)+1$, so $\xi(m-1)<_{\text{lex}}\xi(m)$. If $m\le n_1<n_2$,
      then $\xi_k(n_1)=\xi_{k}(n_2)$ for all $k\in\{1,2,3,4\}$ and $\xi_5$ is increasing.
    \item[\eqref{J6}] Suppose that $[i,j)=\xi_1^{-1}\{k\}$. Since $\xi_1\rest [m,\o)$ is constant,
      either  $j<m$, when we are done by the induction hypothesis, or $i=m$ and $j=\o$. In that case
      one verifies that $\eta\rest[m,\o)\in W=\ran h_{\eta\restl m}$ and then, imitating
      the corresponding argument in the first step, that 
      $$\xi_5\rest [m,\o)=h_{\eta\restl m}^{-1}(\eta\rest [m,\o))$$
      and hence in $\dom h_{\eta\restl m}=P^{\xi_2(m),\xi_3(m)}_{\xi_4(m)}$.
    \end{myItemize}   
    
    Suppose finally that ($\iota\iota\iota$) holds. Then $\dom\eta$ must be $\o$ since otherwise the condition ($\iota\iota\iota$)
    is simply contradictory 
    (because $\eta\rest(\dom\eta-1+1)=\eta$ (except for the case $\dom\eta=0$, but then condition ($\iota$) holds and we are done)).
    By (g), we have $\ran\eta_1=\o$, because otherwise we had $\eta\in \ran F_{\a_{j-1}}$.
    Let $F^{-1}_{\a_j}(\eta)=\xi=\Cup_{n<\o}F^{-1}_{\a_{j-1}}(\eta\rest n)$. 

    Let us check that it is in $J(S)$. Conditions \eqref{J0}--\eqref{J6}
    are satisfied by $\xi$, because they are satisfied by all its initial segments. Let us check \eqref{J7}.

    First of all $\xi$ cannot be in $J^{\a_{j-1}}(S)$, since
    otherwise, by (d) and (i), 
    $$F_{\a_{j-1}}(\xi)=\Cup_{n<\o}F_{\a_{j-1}}(\xi\rest n)=\Cup _{n<\o}\eta\rest n=\eta$$ 
    were again in $\ran F_{\a_{j-1}}$.
    If $j-1$ is a successor ordinal, then we are done: by (b) $\a_{j-1}$ is a successor and 
    we assumed $\eta\in J(S')$, so by
    (e2) we have $\xi\in J(S)$. Thus we can assume that $j-1$ is a limit ordinal. Then by (b),
    $\a_{j-1}$ is a limit ordinal in $C$ and by (a), (e) and (f),
    $\ran F_{\a_{j-1}}=J^{\a_{j-1}}(S')$ and $\dom F_{\a_{j-1}}=J^{\a_{j-1}}(S)$. This implies 
    that $\ran\eta\not\subset \o\times \b^4$ for any $\b<\a_{j-1}$ 
    and by ($\circledast$) on page \pageref{circledast} we must have $\sup\ran\eta_5=\a_{j-1}$ which gives $\a_{j-1}\in S'$ by \eqref{J7}. 
    Since $\a_{j-1}\in C\subset \k\setminus S\sd S'$, we have
    $\a_{j-1}\in S$. Again by $(\circledast)$ and that $\dom F_{\a_{j-1}}=J^{\a_{j-1}}(S)$ by (e1), 
    we have $\sup\ran\xi_5=\a_{j-1}$, thus $\xi$ satisfies the condition \eqref{J7}. 

    Let us check whether all the conditions (a)-(i) are met. (a), (b), (c) are common to the cases ($\iota$), ($\iota\iota$) 
    and ($\iota\iota\iota$) in the definition of $F_{\a_j}^{-1}$ and are easy to verify. 
    Let us sketch a proof for (d); the rest is left to the reader.
    \begin{myAlphanumerate}
    \item[(d)] Let $\eta_1,\eta_2\in \ran F_{\a_{j}}$ and let us show that 
      $$\eta_1\subsetneq\eta_2\iff F^{-1}_{\a_j}(\eta_1)\subsetneq F^{-1}_{\a_j}(\eta_2).$$
      The case where both $\eta_1$ and $\eta_2$ satisfy $(\iota\iota)$ is the interesting one (implies all the others).
      
      So suppose $\eta_1,\eta_2\in (\iota\iota)$. Then there exist $m_1$ and $m_2$ as described in the statement of ($\iota\iota$).
      Let us show that $m_1=m_2$. We have $\eta_1\rest (m_1+1)= \eta_2\rest (m_1+1)$ and $\eta_1\rest (m_1+1)\notin \ran F_{\a_{j-1}}$,
      so $m_2\le m_1$. If $m_2\le m_1$, then $m_2<\dom \eta_1$, since $m_1<\dom \eta_1$. 
      Thus if $m_2\le m_1$, then $\eta_1\rest (m_2+1)=\eta_2\rest (m_2+1)\notin \ran F_{\a_{j-1}}$, which implies $m_2=m_1$.
      According to the definition of $F^{-1}_{\a_j}(\eta_i)(k)$ for $k<\dom \eta_1$, 
      $F^{-1}_{\a_j}(\eta_i)(k)$ depends only on 
      $m_i$ and $\eta\rest m_i$ for $i\in \{1,2\}$. Since $m_1=m_2$ and $\eta_1\rest m_1=\eta_2\rest m_2$,
      we have $F^{-1}_{\a_j}(\eta_1)(k)=F^{-1}_{\a_j}(\eta_2)(k)$ for all $k<\dom\eta_1$.

      Let us now assume that $\eta_1\not\subset \eta_2$. Then take the smallest $n\in \dom\eta_1\cap\dom\eta_2$ such that
      $\eta_1(n)\ne \eta_2(n)$. It is now easy to show that $F^{-1}_{\a_j}(\eta_1)(n)\ne F^{-1}_{\a_j}(\eta_2)(n)$ by the construction.
    \end{myAlphanumerate}
    \noindent\textbf{Even successor step.} Namely the one where $j=\b+n$ and $n$ is even.
    But this case goes exactly as the above completed step, except that we start with $\dom F_{\a_j}=J^{\a_j}(S)$
    where $\a_j$ is big enough successor of an element of $C$ such that $J^{\a_j}(S)$ contains $\ran F_{\a_{j-1}}$
    and define $\xi=F_{\a_j}(\eta)$. Instead of (e) we use (f) as the induction hypothesis. 
    This step is easier since one does not need to care about the successors of limit ordinals.

    \noindent\textbf{Limit step.}
    Assume that $j$ is a limit ordinal. Then let $\a_j=\Cup_{i<j}\a_i$ and $F_{\a_j}=\Cup_{i<j}F_{\a_i}$.
    Since $\a_i$ are successors of ordinals in $C$, $\a_j\in C$, so (b) is satisfied. 
    Since each $F_{\a_i}$ is an isomorphism, 
    also their union is, so (d) is satisfied. 
    Because conditions (e), (f) and (i) hold for $i<j$, the conditions (e) and (i)
    hold for $j$. (f) is satisfied because the premise is not true.
    (a) and (c) are clearly satisfied. Also (g) and (h) are satisfied by Claim~1 since now $\dom F_{\a_j}=J^{\a_j}(S)$
    and $\ran F_{\a_j}=J^{\a_j}(S')$ (this is because (a), (e) and (f) hold for $i<j$).

    \noindent\textbf{Finally} $F=\Cup_{i<\k}F_{\a_i}$ is an isomorphism between $J(S)$ and $J(S')$.
  \end{proofVOf}
\end{proofV}

\begin{Thm}\label{thm:NSIsReducibleToStableUnsuperstable}
  Suppose $\k$ is such that $\k^{<\k}=\k$ and 
  for all $\l<\k$, $\l^\o<\k$ and that $T$ is a stable unsuperstable theory.
  Then $E_{S^\k_\o}\le_c\,\, \cong_T$.
\end{Thm}
\begin{proof}
  For $\eta\in 2^\k$ let $J_\eta=J(\eta^{-1}\{1\})$ where the function $J$ is as in Lemma \ref{lem:StoJS} above.
  For notational convenience, we assume that $J_\eta$ is a downward closed subtree of $\k^{\le\o}$.
  Since $T$ is stable unsuperstable, for all $\eta$ and 
  $t\in J_\eta$, there are finite sequences $a_t=a_t^\eta$ in the monster model
  such that 
  \newcommand{\da}{\mathop{\downarrow}}
  \begin{myEnumerate}
  \item If $\dom(t)=\o$ and $n<\o$ then 
    $$a_t\,\,\,\,\,\,\,\,\,\,\,\,\not\!\!\!\!\!\!\!\!\!\!\!\!\da_{\displaystyle\mathop{\cup}_{m<n}\!a_t\rest m}a_{t\restl n}.$$
  \item For all downward closed subtrees $X,Y\subset J_\eta$, 
    $$\Cup_{t\in X}a_t\da_{\displaystyle\!\!\mathop{\cup}_{t\in X\cap Y}\!a_t}\,\Cup_{t\in Y}a_t$$
  \item For all downward closed subtrees $X\subset J_\eta$ and 
    $Y\subset J_{\eta'}$ the following holds:
    If $f\colon X\to Y$ is an isomorphism, then there is an automorphism $F$ of the monster model such that
    for all $t\in X$, $F(a_t^\eta)=a_{f(t)}^{\eta'}$
  \end{myEnumerate}
  Then we can find an $F^f_\o$-construction 
  $$(\Cup_{t\in J_\eta}a_t,(b_i,B_i)_{i<\k})$$
  (here $(t(b/C),D)\in F^f_\o$ if $D\subset C$ is finite and $b\da_D C$, see \cite{Shelah2})
  such that 
  \begin{myItemize}
  \item[$(\star)$] for all $\a<\k$, $c$ and finite $B\subset \Cup_{t\in J_\eta}a_t\cup \Cup_{i<\a}b_i$ there is
    $\a<\b<\k$ such that $B_\beta=B$ and
    $$\stp(b_\beta/B)=\stp(c/B).$$
  \end{myItemize}
  Then 
  $$M_\eta=\Cup_{t\in J_\eta}a_t\cup \Cup_{i<\k}b_i\models T.$$
  Without loss of generality we may assume that the trees $J_\eta$ 
  and the $F^f_\o$-constructions for $M_\eta$ are chosen coherently enough such that
  one can find a code $\xi_\eta$ for (the isomorphism type of) $M_\eta$ so that 
  $\eta\mapsto \xi_\eta$ is continuous. Thus we are left to show that
  $\eta E_{S^\k_\o} \eta'\iff M_{\eta}\cong M_{\eta'}$
  
  \begin{myItemize}
  \item[``$\Rightarrow$''] Assume $J_\eta\cong J_{\eta'}$. By (3) it is enough to show that
    $F^f_\o$-construction of length $\k$ satisfying $(\star)$ are unique up to isomorphism over
    $\Cup_{t\in J_\eta}a_t$. But $(\star)$ guarantees that the proof of the uniqueness of $F$-primary models from
    \cite{Shelah2} works here.
  \item[``$\Leftarrow$''] Suppose $F\colon M_\eta\to M_{\eta'}$ is an isomorphism and for a contradiction suppose
    $(\eta,\eta')\notin E_{S^\k_\o}$. Let $(J^\a_\eta)_{\a<\k}$ be a filtration of $J_\eta$ and $(J^\a_{\eta'})_{\a<\k}$
    be a filtration of $J_{\eta'}$ (see Definition \ref{def:Filtration} above). For $\a<\k$, let 
    $$M^\a_\eta=\Cup_{t\in J^{\a}_\eta}a_t\cup\Cup_{i<\a}b_i$$
    and similarly for $\eta'$:
    $$M^\a_{\eta'}=\Cup_{t\in J^{\a}_{\eta'}}a_t\cup\Cup_{i<\a}b_i.$$
    Let $C$ be the cub set of those $\a<\k$ such that $F\rest M^\a_\eta$ is onto $M^\a_{\eta'}$ and 
    for all $i<\a$, $B_i\subset M^\a_\eta$ and $B'_i\subset M^\a_{\eta'}$,
    where $(\Cup_{t\in J_{\eta'}},(b'_i,B'_i)_{i<b})$ is in the construction of $M_{\eta'}$.
    Then we can find $\a\in \lim C$ such that in $J_\eta$ there is
    $t^*$ satisfying (a)--(c) below, but in $J_{\eta'}$ there is no such $t^*.$:
    \begin{myAlphanumerate}
    \item $\dom(t^*)=\o$,
    \item $t^*\notin J^{\a}_\eta$,
    \item for all $\b<\a$ there is $n<\o$ such that $t^*\rest n\in J^\a_\eta\setminus J^\b_\eta,$
    \end{myAlphanumerate}
  \end{myItemize}
  Note that
  \begin{myItemize}
  \item[$(\star\star)$] if $\a\in C$ and $c\in M^\a_\eta$, there is a finite $D\subset \Cup_{t\in J^\a_\eta}a_t$
    such that $(t(c,\Cup_{t\in J_\eta}a_t),D)\in F^f_\o$,
  \end{myItemize}
  Let $c=F(a_{t^*})$. By the construction we cat find finite $D\subset M^\a_{\eta'}$, and $X\subset J_{\eta'}$ such that
  $$\Big(t(c,M^\a_{\eta'}\cup \Cup_{t\in J_{\eta'}}a^{\eta'}_t),D\cup \Cup_{t\in X}a^{\eta'}_{t}\Big)\in F^f_\o.$$
  But then there is $\b\in C$, $\b<\a$, such that $D\subset M^\b_{\eta'}$ and if $u\le t$ for some $t\in X$, then
  $u\in J^\b_{\eta'}$ (since in $J_{\eta'}$ there is no element like $t^*$ is in $J_\eta$).
  But then using $(\star\star)$ and (2), it is easy to see that 
  $$c\da_{M^\b_{\eta'}}M^\a_{\eta'}.$$
  On the other hand, using (1), (2), $(\star\star)$ and the choice of $t^*$ one can see that 
  $\displaystyle a_{t^*}\,\,\,\not\!\!\!\da_{M^\b_\eta}M^a_\eta$,
  a contradiction.
\end{proof}

\begin{Open}
  If $\k=\l^+$, $\l$ regular and uncountable, 
  does equality modulo $\l$-non-sationary ideal, $E_{S^\k_\l}$, Borel reduce to $T$ for all stable unsuperstable~$T$?
\end{Open}

\chapter{Further Research}

In this chapter we merely list all the questions that also appear in the text:

\begin{Open}
  Is it consistent that Borel* is a proper subclass of $\Sii$, or even equals $\Dii$? Is it consistent that
  all the inclusions are proper at the same time: $\Dii\subsetneq {\rm Borel}^*\subsetneq \Sii$?
\end{Open}

\begin{Open}
  Does the direction left to right of Theorem \ref{thm:BorelIsLkk} hold without the assumption
  $\k^{<\k}=\k$?
\end{Open}

\begin{Open}
  Under what conditions on $\kappa$ does the conclusion of Theorem \ref{thm:VaughtConjecture} hold?
\end{Open}

\begin{Open}
  Is the Silver Dichotomy for uncountable $\k$ consistent?
\end{Open}

\begin{Open}
  Can there be two equivalence relations, 
  $E_1$ and $E_2$ on $2^{\k}$, $\k>\o$ such that $E_1$ and $E_2$ are Borel and incomparable,
  i.e. $E_1\not\le_B E_2$ and $E_2\not\le_B E_1$? 
\end{Open}

 \begin{Open}
    Is it consistent that $S^{\omega_2}_{\omega_1}$ Borel reduces to $S^{\omega_2}_\omega$?
 \end{Open}

\begin{Open}
  We proved that the isomorphism relation of a theory $T$ is Borel if and only if $T$ is classifiable and
  shallow. Is there a 
  connection between the depth of a shallow theory and the Borel degree of its isomorphism relation?
  Is one monotone in the other?
\end{Open}

\begin{Open}
  Can it be proved in ZFC that if $T$ is stable unsuperstable then $\cong_T$ is not $\Dii$?
\end{Open}

\begin{Open}
  If $\k=\l^+$, $\l$ regular and uncountable,  does equality modulo $\l$-non-sationary ideal, $E_{S^\k_\l}$, Borel reduce to $T$ for all stable 
  unsuperstable~$T$?
\end{Open}

\newcommand{\dlo}{\text{dlo}}
\newcommand{\gr}{\text{gr}}

\begin{Open}
  Let $T_{\dlo}$ be the theory of dense linear orderings without end points and $T_{\gr}$ the theory of random graphs.
  Does the isomorphism relation of $T_{\gr}$ Borel reduce to $T_{\dlo}$, i.e. $\cong_{T_{\gr}}\le_B \cong_{T_{\dlo}}$?
\end{Open}

\end{document}